\documentclass[12pt]{article}
\usepackage{mathrsfs}
\usepackage[active]{srcltx}
\usepackage[curve]{xypic}
\usepackage{amsthm}
\usepackage{amsfonts}
\usepackage{amscd}
\usepackage{amsmath}
\usepackage{amssymb}
\usepackage{color}
\usepackage{eqnarray}
\usepackage{array}
\usepackage{epsfig}
\usepackage{verbatim}
\usepackage{psfrag}
\usepackage{graphicx,caption2}
\usepackage[usenames,dvipsnames]{pstricks}
\usepackage{epsfig}
\usepackage{pst-grad} 
\usepackage{pst-plot} 

\input epsf.tex

\newcommand{\mylabel}[1]{\label{#1}}

\newtheorem{theorem}{Theorem}[section]
\newtheorem{lemma}[theorem]{Lemma}
\newtheorem{proposition}[theorem]{Proposition}
\newtheorem{corollary}[theorem]{Corollary}

\newcommand{\REFEQN}[1] { \begin{equation}\mylabel{#1} }
\newcommand{\ENDEQN}{\end{equation}}
\newcommand{\REFTHM}[1] { \begin{theorem}\mylabel{#1} }
\newcommand{\ENDTHM}{\end{theorem}}
\newcommand{\REFPROP}[1]{\begin{proposition}\mylabel{#1} }
\newcommand{\ENDPROP}{\end{proposition} }
\newcommand{\REFLEM}[1]{\begin{lemma}\mylabel{#1} }
\newcommand{\ENDLEM}{\end{lemma} }
\newcommand{\REFCOR}[1]{\begin{corollary}\mylabel{#1} }
\newcommand{\ENDCOR}{\end{corollary} }

\def\Ga{\Gamma}

\def\la{\lambda}

\def\wt{\widetilde}
\def\wh{\widehat}

\def\CC{\mbox{$\mathbb C$}}
\def\DD{\mbox{$\mathbb D$}}
\def\AA{\mbox{$\mathbb A$}}

\def\NN{\mbox{$\mathbb N$}}

\def\TT{\mbox{$\mathbb T$}}
\def\ZZ{\mbox{$\mathbb Z$}}

\def\BBB{{\mathcal B}}

\def\FFF{{\mathcal F}}

\def\III{{\mathcal I}}
\def\JJJ{{\mathcal J}}

\def\MMM{{\mathcal M}}

\def\OOO{{\mathcal O}}
\def\PPP{{\mathcal P}}
\def\QQQ{{\mathcal Q}}
\def\RRR{{\mathcal R}}
\def\SSS{{\mathcal S}}
\def\TTT{{\mathcal T}}
\def\UUU{{\mathcal U}}
\def\VVV{{\mathcal V}}
\def\WWW{{\mathcal W}}

\def\fm{{\mathfrak M}}

\def\fw{{\mathfrak W}}

\def\sA{{\mathscr A}}
\def\sC{{\mathscr C}}
\def\sD{{\mathscr D}}

\def\sP{{\mathscr P}}
\def\sQ{{\mathscr Q}}
\def\sR{{\mathscr R}}

\def\cbar{\overline{\CC}}
\def\md{\!\!\!\!\mod\!}
\def\beginp{{\noindent\em Proof.} }
\def\smm{{\smallsetminus}}
\def\ds{\displaystyle }
\def\height{\text{\rm Height}}
\def\width{\text{\rm Width}}
\def\area{\text{\rm Area}}

\begin{document}

\title{Hyperbolic-parabolic deformations \\
of rational maps
\footnote{2010 Mathematics Subject
Classification: 37F20, 37F45}}
\author{Guizhen CUI
\thanks{supported by NSFC grant no. 11125106.}
\and Lei TAN }
\date{\today}
\maketitle
\begin{abstract}
We develop a Thurston-like theory to characterize geometrically finite rational maps, then apply it to study pinching and plumbing deformations of rational maps. We show that in certain conditions the pinching path converges uniformly and the quasiconformal conjugacy converges uniformly to a semi-conjugacy from the original map to the limit. Conversely, every geometrically finite rational map with parabolic points is the landing point of a pinching path for any prescribed plumbing combinatorics.
\end{abstract}
\thispagestyle{empty}

\newpage
\pagenumbering{roman}
\tableofcontents

\newpage
\pagenumbering{arabic}

\section{Introduction}

One of the most important results in complex dynamics is Thurston's topological characterization of rational maps. It states that a post-critically finite branched covering of the $2$-sphere with a hyperbolic orbifold is combinatorially equivalent to a rational map if and only if it has no Thurston obstructions \cite{DH2}. This condition has been proved to be necessary for geometrically finite rational maps \cite{Mc1}. C. McMullen \cite{Bi} proposed extending Thurston's Theorem to this case.

This problem has been solved for the sub-hyperbolic case \cite{CT1, JZ}. In the first part of this work, we solve this problem in the presence of parabolic points. The main tool is pinching.

In general, pinching provides a path of quasiconformal deformations of a hyperbolic Riemann surface, whose Beltrami differential is supported in a neighborhood of a finite disjoint union of simple closed geodesics. Along the deformations the lengths of these geodesics shrink to zero. The limit is a stable curve in the Deligne-Mumford compactification.

One may perform such deformations on Riemann surfaces occurring as the quotient space of a Kleinian group. See \cite{Mas} for related results. Pinching has also been applied in order to study the parameter space of the dynamics of rational maps \cite{Mak}, where the simple closed geodesics on the quotient space of the rational map are chosen such that their lifts to the dynamical space are simple closed curves.

For our purpose, creating parabolic points from attracting points, we have to choose the simple closed geodesics on the quotient space such that their lifts to the dynamical space are arcs which join attracting periodic points to the Julia set. Such a choice makes the control of the distortion of quasiconformal conjugacy more difficult.

In this work, we first study pinching on specific simple closed geodesics on the quotient spaces of rational maps. We will call this type of pinching {\it simple pinching}. By using the length-area method to control the distortion of quasiconformal conjugacy, we prove that the pinching path is convergent. This result was generalized to other cases \cite{H3}.

Plumbing is a surgery on a nodded Riemann surface which is like an inverse to pinching: it replaces pairs of cusp neighborhoods by thin annuli. We will apply this surgery to create attracting points from parabolic points in the dynamics of rational maps.

Begin with a semi-rational map with parabolic points. We start by constructing a sub-hyperbolic semi-rational map by plumbing. We prove that if the original map has neither a Thurston obstruction nor a {\it connecting arc}, which can be viewed as a degenerate Thurston obstruction, then the resulting semi-rational map has no Thurston obstruction. Now applying the characterization theorem obtained in \cite{CT1, JZ} and the above result on simple pinching, we obtain a rational map that is {\it c-equivalent} to the original semi-rational map.

Simple pinching can also be used to study a conjecture proposed by L. Goldberg and J.Milnor \cite{GM}, which states that if a polynomial has a parabolic cycle, then its immediate basins can be converted to be attracting by a small perturbation without changing the topology of the Julia set.  We prove this conjecture in the setting of geometrically finite rational maps. Refer \cite{H1,H2,H3} and \cite{K1,K2} for proofs of the same conjecture in various settings, but with different methods.

The second part of this work is devoted to general pinching of rational maps. Unfortunately, the distortion control for simple pinching is not valid in the general case. We start by studying plumbing instead of pinching. We develop a new distortion control for univalent maps, then show that any plumbing of a geometrically finite rational map can be realized as a pinching path of a geometrically finite rational map converging to the original map. Applying the characterization theorem obtained in the first part, we prove that a pinching path is convergent if it satisfies the {\it non-separating} condition. This condition has been proved to be necessary \cite{T2}. As a by-product, we show that the quasiconformal deformation of a geometrically finite rational map supported on a non-separating multi-annulus is bounded in its moduli space.

The above result is an analogue of Maskit's theorem \cite{Mas} for producing cusps. It is known that cusps are dense in the space of boundary groups \cite{Mc3}. It has been asked whether cusps are dense in the boundary of a hyperbolic component of rational maps \cite{Mc4}. Our results about pinching may be helpful in solving this problem.

\vskip 0.24cm

{\bf Main results}. We now give definitions and statements. Let $F$ be a branched covering of the Riemann sphere $\cbar$ with degree $\deg F\ge 2$. Denote by $\Omega_F =\{z\in\cbar:\ \deg_z F>1\}$ the {\bf critical set} of $F$, denote by
$$
\PPP_F =\overline{\ds\bigcup_{n>0}F^n(\Omega_F)}
$$
the {\bf post-critical set} of $F$ and denote by $\PPP'_F$ the set of {\bf accumulation points} of $\PPP_F$. The map $F$ is called {\bf post-critically finite}\/ if $\PPP_F$ is finite, and {\bf geometrically finite}\/ if $\PPP'_F$ is finite.

Let $f$ be a rational map of $\cbar$ with $\deg f\ge 2$. Denote by $\FFF_f$ the Fatou set of $f$ and by $\JJJ_f$ the Julia set of $f$; refer to \cite{B} or \cite{Mi} for the definitions. In the literature, the geometric finiteness of a rational map $f$ is defined by the condition that $\PPP_f\cap\JJJ_f$ is finite. It turns out that the two definitions are equivalent when $f$ is a rational map \cite{CJ}.

\vskip 0.24cm

A geometrically finite branched covering $F$ is called a {\bf semi-rational map}\/ if

(1) $F$ is holomorphic in a neighborhood of $\PPP'_F$,

(2) each cycle in $\PPP'_F$ is attracting or super-attracting or parabolic, and

(3) any attracting petal at a parabolic periodic point in $\PPP'_F$ contains points of $\PPP_F$.

A semi-rational map $F$ is called {\bf sub-hyperbolic}\/ if $\PPP'_F$ contains no parabolic cycles.

See \cite{Mi} or \S2.6 for the definitions of attracting petals and attracting flowers. The condition (3) implies that any attracting petal at a parabolic periodic point in $\PPP'_F$ contains infinitely many points of $\PPP_F$.

\vskip 0.24cm

Let $F$ be a semi-rational map. An open set $\UUU\subset\cbar$ is called a {\bf fundamental set}\/ of $F$ if $\UUU\subset F^{-1}(\UUU)$ and $\UUU$ contains every attracting and super-attracting cycle in $\PPP'_F$ and an attracting flower at each parabolic cycle in $\PPP'_F$.

A fundamental set could be the empty set if $F$ is post-critically finite. It is contained in the Fatou set if $F$ is a rational map.

\vskip 0.24cm

Two semi-rational maps $F$ and $G$ are called {\bf c-equivalent}\/ if there exist a pair of orientation-preserving homeomorphisms $(\phi, \psi)$ of $\cbar$ and a fundamental set $\UUU$ of $F$ such that:

(a) $\phi\circ F=G\circ\psi$,

(b) $\phi$ is holomorphic in $\UUU$,

(c) $\psi=\phi$ in $\UUU\cup\PPP_F$ and $\psi$ is isotopic to $\phi$ rel $\UUU\cup\PPP_F$.

\vskip 0.24cm

Let $G$ be a semi-rational map with parabolic cycles in $\PPP'_G$. An open arc $\beta\subset\cbar\smm\PPP_G$ which joins two points $z_0, z_1\in\PPP'_G$ is called a {\bf connecting arc}\/ if:

(i) either $z_0\neq z_1$, or $z_0=z_1$ and both components of $\cbar\smm\overline{\beta}$ contains points of $\PPP_G$,

(ii) $\beta$ is disjoint from a fundamental set of $G$, and

(iii) $\beta$ is isotopic to a component $\wt\beta$ of $G^{-p}(\beta)$ rel $\PPP_G$ for some integer $p>0$, i.e. there exists an isotopy $H: [0,1]\times\cbar\to\cbar$ with $H(0,\cdot)=\text{\rm id}$ and $H(t, \cdot)=\text{\rm id}$ on $\PPP_G$ for $t\in [0,1]$ such that $H(1, \beta)=\wt\beta$.

\REFTHM{unicity}{\bf (Unicity)} Two c-equivalent geometrically finite rational maps with infinite post-critical sets are holomorphic conjugate in the isotopy class of the c-equivalence.
\ENDTHM

\REFTHM{existence}{\bf (Existence)} A semi-rational map with infinite post-critical set is c-equivalent to a rational map if and only if it has neither Thurston obstructions nor connecting arcs.
\ENDTHM

See \cite{DH2, Mc2} or \S4 below for the definition of Thurston obstructions. In order to prove Theorem \ref{existence},  we first establish the following two results which have independent interest. See \S5 for the precise definition of pinching.

\REFTHM{s-pinching}{\bf (Simple pinching)} Let $f$ be a geometrically finite rational map and let $f_t=\phi_t\circ f\circ\phi_t^{-1}$ $(t\ge 0)$ be a simple pinching path starting from $f=f_0$. Then the following properties hold:

(a) $f_t$ converges uniformly to a geometrically finite rational map $g$ as $t\to\infty$.

(b) $\phi_t$ converges uniformly to a continuous onto map $\varphi$ of $\cbar$ as $t\to\infty$.

(c) $\varphi\circ f=g\circ\varphi$ and $\varphi:\, \JJJ_f\to\JJJ_g$ is a homeomorphism.
\ENDTHM

\REFTHM{s-plumbing}{\bf (Simple plumbing)} Any geometrically finite rational map with parabolic cycles is the limit of a simple pinching path starting from a sub-hyperbolic rational map.
\ENDTHM

Denote by $\fm_d$ the space of holomorphic conjugate classes of rational maps of degree $d\ge 2$. For $[f]\in\fm_d$, define $\fm[f]\subset\fm_d$ by
$[g]\in\fm[f]$ if $g$ is quasiconformally conjugate to $f$.

Refer to \S5 for the definition of {\bf non-separating}.

\REFTHM{pinching}{\bf (Pinching)} Let $f$ be a geometrically finite rational map. Let $f_t=\phi_t\circ f\circ\phi_t^{-1}$ ($t\ge 0$) be a pinching path starting from $f=f_0$ and supported on a non-separating multi-annulus. Then the following properties hold:

(a) $f_t$ converges uniformly to a geometrically finite rational map $g$ as $t\to\infty$.

(b) $\phi_t$ converges uniformly to a continuous onto map $\varphi$ of $\cbar$ as $t\to\infty$.

(c) $\varphi(\JJJ_f)=\JJJ_g$.

(d) $\fm[g]\subset\partial\fm[f]$.
\ENDTHM

The pinching deformation can be reversed via a plumbing surgery on the limit rational map $g$. The complete set of possibilities for plumbing can be encoded by a finite set of combinatorial data -- plumbing combinatorics; see \S10.1 for the definition.

\REFTHM{plumbing}{\bf (Plumbing)} Let $g$ be a geometrically rational map with parabolic cycles and let $\sigma$ be a plumbing combinatorics of $g$. Then there exist a geometrically finite rational map $f$ and a pinching path $f_t=\phi_t\circ f\circ\phi_t^{-1}$ $(t\ge 0)$ starting from $f=f_0$ such that $f_t$ is a plumbing of $g$ along $\sigma$ and $\{f_t\}$ converges uniformly to $g$ as $t\to\infty$.
\ENDTHM

For each $[f]\in\fm_d$ and any non-separating multi-annulus $\sA$ in its quotient space, define $\fm[f, \sA]\subset\fm[f]$ by $[g]\in\fm[f, \sA]$ if there exists a quasiconformal map $\phi$ of $\cbar$ such that $g\circ\phi=\phi\circ f$ and the Beltrami differential of $\phi$ is supported on $\pi_f^{-1}(\sA)$.
In general, $\fm[f]$ need not have compact closure in $\fm_d$. However, we have:

\REFTHM{boundedness}{\bf (Boundedness)} Let $f$ be a geometrically rational map and let $\sA$ be a non-separating multi-annulus. Then $\fm[f, \sA]$ is compactly contained in $\fm_d$ and any rational map in the closure of $\fm[f, \sA]$ is geometrically finite.
\ENDTHM

\vskip 0.24cm

{\bf Outline of the paper}.
In \S2, we recall some basic results which will be used in the sequel. Most of them are known except for some lemmas whose proofs are not difficult. In \S3 we prove Theorem \ref{unicity} by considering local conjugacy at parabolic points and the boundary dilatation of a c-equivalence. In \S4 we begin by recalling Thurston's Theorem. Then we show that  geometrically finite rational maps have no connecting arcs, which is the necessary part of Theorem \ref{existence}.  We also show that Thurston's algorithm is convergent. In \S5 we give the definition of the pinching path of a rational map through a favorable model. The proof of Theorem \ref{s-pinching} is given in \S6.

In \S7 we prove Theorem \ref{s-plumbing} and complete the proof of Theorem \ref{existence}. The strategy is to make a detour to sub-hyperbolic semi-rational maps via plumbing and pinching: Starting with a semi-rational map $G$, we first construct a sub-hyperbolic semi-rational map $F$ from $G$ by simple plumbing. Then we show that $F$ has no Thurston obstructions and thus is c-equivalent to a rational map $f$. Finally we show that the simple pinching limit of $f$ is a rational map c-equivalent to $G$.

In the last three sections, we study general pinching and plumbing. In \S8 we define a new type of distortion for univalent maps and give a control for it using a property of the domains of the univalent maps. Then we apply this distortion control to the dynamics of rational maps. The proofs of Theorems \ref{pinching}, \ref{plumbing} and \ref{boundedness} are given in the last two sections.

The proof of Theorem \ref{pinching} is different from the proof of Theorem \ref{s-pinching} and is quite involved. This is because the distortion control in the proof of Theorem \ref{s-pinching} cannot be applied to general pinching. We start with a rational map $f$ and its pinching path; instead of showing directly that the pinching path is convergent, we first construct a semi-rational map $G$ from the limit of truncated quasiconformal maps whose convergence is easy to prove. Then we show that the map $G$ has neither Thurston obstructions nor connecting arcs, and hence is c-equivalent to a rational map $g$ by Theorem \ref{existence}. This provides us the candidate limit map $g$ of our pinching path. Now using a similar strategy as in the proof of Theorem \ref{existence}, we plumb $g$ and then pinch. We have to check that we get exactly the same pinching path, and that it converges uniformly to $g$ by the distortion control established in Theorem \ref{distortion}.

\vskip 0.24cm

{\bf Notation}. The following notation and conventions will be used in this paper.
\begin{align*}
& \DD(z_0, r)=\{z\in\CC: |z-z_0|<r\}\quad\text{for } z_0\in\CC\text{ and } r>0.  \\
& \DD(r)=\DD(0,r)\quad\text{and}\quad \DD=\DD(0,1).  \\
& \DD^*(r)=\DD(r)\smm\{0\}\quad\text{and}\quad \DD^*=\DD^*(1). \\
& \AA(z_0; r_1, r_2)=\{z\in\CC: r_1<|z-z_0|<r_2\}\quad\text{for $z_0\in\CC$ and $0<r_1<r_2$}.  \\
& \AA(r_1, r_2)=\{z\in\CC: r_1<|z|<r_2\}\quad\text{for } 0<r_1<r_2  \\
& \AA(r)=\{z\in\CC: 1/r<|z|<r\}\quad\text{for } r>1.
\end{align*}

$\cbar=\CC\cup\{\infty\}$ is the Riemann sphere and $\CC^*=\CC\smm\{0\}$.

$U\Subset V$ if $\overline{U}\subset V$ for open sets $U, V\subset\cbar$.

A {\bf disk} is a Jordan domain in $\cbar$.

An annulus $A_1$ is contained {\bf essentially} in another annulus $A_2$ if $A_1\subset A_2$ and $A_1$ separates the two boundary components of $A_2$.

\section{Preliminaries}

\subsection{Modulus and extremal length}

Let $A\subset\cbar$ be an annulus such that each of its two complementary components contains at least two points. Then there exist a constant $r>1$ and a conformal map $\chi_A:\, A\to\AA(r)$, where $r$ is unique and $\chi_A$ is unique up to post-composition of a rotation and an inversion $z\mapsto 1/z$.

$\bullet$ The {\bf modulus} of $A$ is defined by $\md A=(\log r)/\pi$.

$\bullet$ $\chi_A^{-1}\{z: \arg z=\theta\}$ is a {\bf vertical line} in $A$ for $0\le\theta<2\pi$,

$\bullet$ $\chi_A^{-1}\{z: |z|=\rho\}$ is a {\bf horizontal circle} in $A$ for $1/r<\rho<r$ and

$\bullet$ $e(A)=\chi_A^{-1}\{z: |z|=1\}$ is the {\bf equator} of $A$.

\vskip 0.24cm

Refer to \cite{Mc1} Theorem 2.1 for the next lemma (the statement is slightly different but with the same proof).

\REFLEM{mod1} Let $A\subset\CC$ be an annulus with $\md A >\frac{5\log 2}{2\pi}$ and let $z_0$ be a point in the bounded component of $\CC\smm A$. Then there exists an annulus $\AA(z_0; r_1, r_2)$ contained essentially in $A$ such that
$$
\md \AA(z_0; r_1, r_2) \ge\md A -\frac{5\log 2}{2\pi}.
$$
\ENDLEM

The modulus of an annulus is related to extremal length as follows. Let $\rho(z)$ be a non-negative Borel measurable function on $A$ satisfying
$$
0<\area(\rho,A)=\iint_A\rho^2(z)dxdy<\infty.
$$
The {\bf $\rho$-length} of a locally rectifiable arc $\alpha\subset A$ is:
$$
L(\rho,\alpha)=\int_{\alpha}\rho(z)|dz|.
$$
Let $\height(\rho, A)$ be the infimum of $L(\rho,\alpha)$ over all locally rectifiable arcs $\alpha\subset A$ which join the two components of $\cbar\smm A$. Let $\width(\rho, A)$ be the infimum of $L(\rho, \gamma)$ over all locally rectifiable Jordan curves $\gamma\subset A$ which separate the two components of $\cbar\smm A$. The following classical inequalities (refer to \cite{L}) will be used several times in this paper.

\REFLEM{mod2}
$$
\frac{\height(\rho, A)^2}{\area(\rho,A)}\le \md A
\le\frac{\area(\rho,A)}{\width(\rho, A)^2}.
$$
Both equalities hold for $\rho(z)=|(\log\chi_A)'(z)|$, which is called an {\bf extremal metric} of $A$.
\ENDLEM

A {\bf (topological) quadrilateral}\/ $Q=Q(\alpha, \alpha')$ is a disk $Q$ in $\cbar$ together with a pair of open arcs $\alpha, \alpha'\subset\partial Q$  that have disjoint closures. We will call the two arcs $\alpha$ and $\alpha'$ {\bf the horizontal sides} of $Q$, and the two arcs $\partial Q\smm(\overline{\alpha}\cup \overline{\alpha'})$ the {\bf vertical sides} of $Q$. There exist a unique constant $b>0$ and a conformal map
$$
\chi_Q:\ Q\to R_b=\{z=x+iy:\, 0\!<\!x\!<\!1, 0\!<\!y\!<\!b\}
$$
such that the continuous extension of $\chi_Q$ maps the two horizontal sides $(\alpha, \alpha')$ onto the two horizontal sides $(0,1)$ and $(ib, 1+ib)$ of the rectangle $R_b$.

$\bullet$ The {\bf modulus} of $Q$ is defined by $\md Q=b$.

$\bullet$ $\chi_Q^{-1}\{x+iy:\ x=x_0\}$ is a {\bf vertical line} of $Q$ for $0<x_0<1$ and

$\bullet$ $\chi_Q^{-1}\{x+iy: y=y_0\}$ is a {\bf horizontal line} of $Q$ for $0<y_0<b$.

\vskip 0.24cm

Lemma \ref{mod2} also holds for the quadrilateral $Q$, where $\height(\rho, A)$ and $\width(\rho, A)$ are replaced by the following: $\height(\rho, Q)$ is the infimum of the $\rho$-length of all locally rectifiable arcs in $Q$ which join the two horizontal sides of $Q$, whereas $\width(\rho, Q)$ is the infimum of the $\rho$-length of all locally rectifiable arcs in $Q$ which join the two vertical sides of $Q$. The following two lemmas will be used in the proof of Lemma \ref{nested-annuli}.

\begin{lemma}\label{mod3}{\bf (From quadrilaterals to an annulus)}
Let $Q_i(\alpha_i, \alpha'_i)$ ($1\le i\le n$) be quadrilaterals such that $\cup_{i=1}^n Q_i=A$ is an annulus and
$$
\bigcup_{i=1}^n(\partial Q_i\smm\alpha'_i)=\partial_{+}A=\bigcup_{i=1}^n\overline{\alpha_i},\eqno (2.1)
$$
where $\partial_{+}A$ is one of the two components of $\partial A$ (see Figure 1(a)). Then
$$
\dfrac 1{\md A }\le\ds\sum^n_{i=1}\frac 1{\md Q_i}.
$$
\end{lemma}

\beginp
Set $\rho_i(z)=\frac{|\chi'_{Q_i}(z)|}{\md Q_i}$ on $Q_i$ and $\rho_i(z)=0$ elsewhere for $1\le i\le n$. Then
$$
\area(\rho_i,Q_i)=\frac1{\md Q_i}\quad\text{and}\quad \height(\rho_i, Q_i)=1.
$$
Set $\rho(z)=\max\{\rho_i(z)\}$. Then
$$
\area(\rho,A)\le\sum_{i=1}^n\frac{1}{\md Q_i}.
$$
By assumption (2.1), the other component of $\partial A$ is contained in $\cup\alpha'_i$. Let $\delta: (0,1)\to A$ be a locally rectifiable arc in $A$ which joins the two components of $\partial A$. By the second equation of condition (2.1), one end of $\delta$ must be contained in $\alpha_j$ for some $0\le j\le n$. By the first equation of the condition, either $\delta$ is totally contained in $Q_j$ and hence the other endpoint of $\delta$ lands on $\alpha'_j$, or $\delta$ intersects $\alpha'_j$.  In both cases, we have $L(\rho, \delta)\ge L(\rho_j, \delta)\ge 1$. So $\height(\rho, A)\ge 1$. Now the lemma follows from Lemma \ref{mod2}.
\qed

\begin{figure}[htbp]
\begin{center}
\includegraphics[scale=0.9]{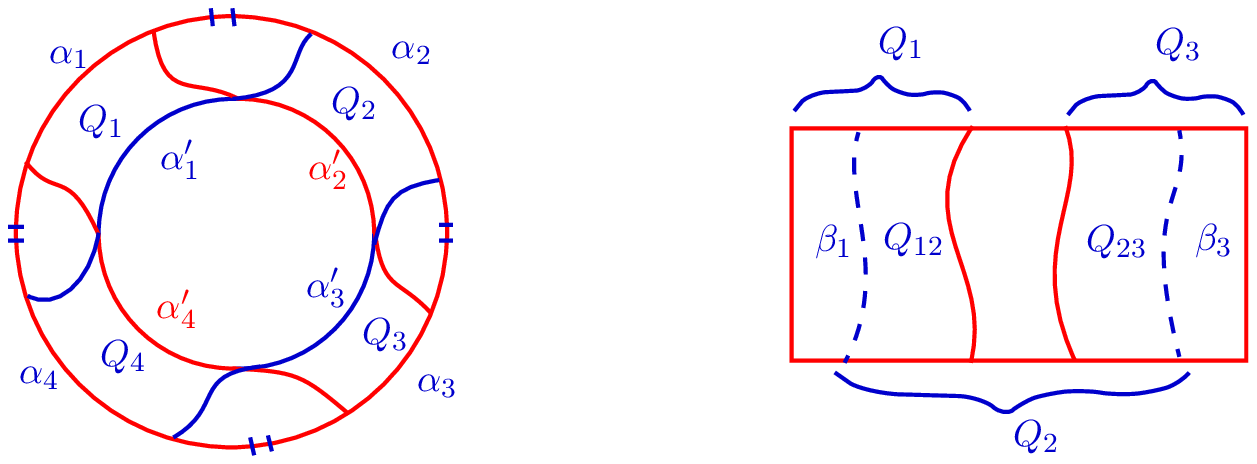}
\end{center}
\begin{center}{\sf Figure 1. (a) Four quadrilaterals form an annulus. (b) Sub-quadrilaterals.}
\end{center}
\end{figure}

\vskip 0.24cm

A quadrilateral $Q_0(\alpha_0, \alpha'_0)$ is a {\bf sub-quadrilateral}\/ of a quadrilateral $Q(\alpha, \alpha')$ if $Q_0\subset Q$ and $(\alpha_0\cup\alpha'_0)\subset(\alpha\cup\alpha')$.

Let $Q$ be a quadrilateral and let $Q_1, Q_3\subset Q$ be disjoint sub-quadrilaterals of $Q$ such that each of the two vertical sides of $Q$ is a vertical side of $Q_1$ or $Q_3$. Let $\beta_1$ and $\beta_3$ be vertical lines in $Q_1$ and $Q_3$, respectively. Then we get a sub-quadrilateral $Q_2$ of $Q$ such that $(\beta_1, \beta_3)$ are vertical sides of $Q_2$. Let $Q_{12}=Q_1\cap Q_2$, which is a sub-quadrilateral of $Q_1$ and $Q_2$.  Let $Q_{23}=Q_3\cap Q_2$, which is a sub-quadrilateral of $Q_3$ and $Q_2$ (see Figure 2(b)). Let
$$
M=\max\{\!\!\md Q_{12}, \md Q_{23}, 1\}.
$$

\REFLEM{mod4} {\bf (Three overlapping quadrilaterals)}
$$
\dfrac 1{\md Q }\le\sum_{i=1}^3\frac{M^2}{\md Q_i}.
$$
\ENDLEM

\beginp
Set $\rho_i(z)=\frac{|\chi'_{Q_i}(z)|}{\md Q_i}$ on $Q_i$ and $\rho_i(z)=0$ otherwise for $i=1,2,3$. Then
$$
\area(\rho_i,Q_i)=\dfrac1{\md Q_i}\quad\text{and}\quad \height(\rho_i, Q_i)=1.
$$
Set $\rho(z)=\max\{\rho_i(z):\, i=1,2,3\}$. Then $Q=E_1\cup E_2\cup E_3$, where $E_i$ is defined by
$$
E_i=\{z\in Q: \, \rho(z)=\rho_i(z)\},\ i=1,2,3.
$$
Since $\rho(z)>0$ for $z\in Q$, we have $E_i\subset Q_i$ and
\begin{align*}
\area(\rho, Q) & \le  \area(\rho, E_1)+\area(\rho,E_2)+\area(\rho, E_3) \\
& =\area(\rho_1, E_1)+\area(\rho_2,E_2)+\area(\rho_3, E_3) \\
& \le \area(\rho_1, Q_1)+\area(\rho_2, Q_2)+\area(\rho_3, Q_3) \\
& = \frac{1}{\md Q_1}+\frac{1}{\md Q_2}+\frac{1}{\md Q_3} \ .
\end{align*}

Since $\beta_1$ is a vertical line in $Q_1$, the quadrilateral $Q_{12}$ becomes a rectangle in the rectangle model of $Q_1$. Therefore $\rho_1$ restricted to $Q_{12}$ is an extremal metric of $Q_{12}$. Thus
$$
\area(\rho_1,Q_{12})=\width(\rho_1, Q_{12})=\frac1{\md Q_{12}}.
$$
Similarly,
$$
\area(\rho_3,Q_{23})=\width(\rho_3, Q_{23})=\frac1{\md Q_{23}}.
$$

For any arc $\alpha$ in $Q$ which joins the two horizontal sides of $Q$, either $\alpha$ is contained in $Q_i$ for some $i\in \{1,2,3\}$ and joins the two horizontal sides of $Q_i$, or $\alpha$ intersects the two vertical sides of either $Q_{12}$ or $Q_{23}$. In the former case,
$$
L(\rho, \alpha)\ge L(\rho_i, \alpha)\ge\height(\rho_i, Q_i)=1.
$$
In the latter case, suppose $\alpha$ intersects the two vertical sides of $Q_{12}$; then there is a sub-arc $\alpha'$ of $\alpha$ which stays in $Q_{12}$ and joins the two vertical sides of $Q_{12}$. Thus
$$
L(\rho,\alpha)\ge L(\rho, \alpha')\ge L(\rho_1, \alpha')\ge \width(\rho_1, Q_{12})= \frac{1}{\md Q_{12}}\ge\frac{1}{M}.
$$
In summary, we have $L(\rho, \alpha)\ge 1/M$ and hence $\height(\rho, Q)\ge 1/M$. Applying Lemma \ref{mod2}, we obtain
$$
\frac{1}{\md Q}\le\frac{\area(\rho, Q)}{\height(\rho, Q)^2}\le \sum_{i=1}^3\frac{M^2}{\md Q_i}.
$$
\qed

\REFTHM{continuity}
Let $\{\phi_t\}$ ($t\ge 0$) be a family of homeomorphisms of $\cbar$ which converges uniformly to a homeomorphism $\phi$ of $\cbar$ as $t\to\infty$.
Let $A\subset\cbar$ be an annulus with $\md(A)<\infty$. Then
$$
\lim_{t\to\infty}\md\phi_t(A)=\md\phi(A).
$$
\ENDTHM

\beginp
We may assume that $\phi$ is the identity by considering the sequence $\{\phi_t\circ\phi^{-1}\}$. Then $\{\phi_t(A)\}$ converges to $A$ as $t\to\infty$ in the sense of Carath\'{e}odory, i.e.:

(1) any compact subset $E\subset A$ is contained in $\phi_t(A)$ once $t$ is large enough, since $\{\phi_t\}$ converges uniformly to the identity as $t\to\infty$, and

(2) if $U\subset\cbar$ is an open set and $t_0\ge 0$ is a constant such that $U\subset\phi_t(A)$ for $t\ge t_0$, then $U\subset A$.

In fact, for any point $z\in U$, since $\phi_t(z)$ converges to the point $z$ as $t\to\infty$, there exists a constant $t_1\ge 0$ such that $\phi_t(z)\in U$ foe $t\ge t_1$. Thus $\phi_t(z)\in U\subset\phi_t(A)$ for $t\ge\max\{t_0, t_1\}$. So $z\in A$. Therefore $U\subset A$.

Each of the two components $E_1, E_2$ of $\cbar\smm A$ contains more than one point since $\md(A)<\infty$. Since $\{\phi_t\}$ converges uniformly to the identity, for $i=1,2$,
$$
\lim_{t\to\infty}\text{\rm diam}_s\phi_t(E_i)=\text{\rm diam}_s E_i>0.
$$
Therefore,
$$
\limsup_{t\to\infty}\md\phi_t(A)<\infty
$$
by Teichm\"{u}ller's Theorem (refer to \cite{A}, Theorem 4-7). Let
$$
\chi_t:\,\phi_t(A)\to B_t:=\{z:\, 1/r_t<|z|<r_t\}
$$
be a conformal map, where $(\log r_t)/\pi=\md\phi_t(A)$. Then there exists a constant $t_0\ge 0$ such that the family $\{\chi_t\}$ ($t>t_0$) is uniformly bounded.

For any sequence $t_n\in (t_0,\infty)$ ($n\in\NN$) with $t_n\to\infty$ as $n\to\infty$, by the Motel principle, there exists a subsequence (also denoted by $\{t_n\}$) and a holomorphic map $\chi$ on $A$ such that for any compact set $E\subset A$, $\chi_{t_n}$ converges uniformly to $\chi$ on $E$.

By  Carath\'{e}odory's Theorem (refer to \cite{P}, Theorem 1.8), either the map $\chi$ is constant or it is univalent, and in the latter case, $\{B_{t_n}\}$ converges to $\chi(A)$ as $n\to\infty$, i.e.

(a) any compact subset $E\subset\chi(A)$ is contained in $B_{t_n}$ once $n$ is large enough, and

(b) if $U\subset\cbar$ is an open set and $n_0\ge 0$ is an integer such that $U\subset B_{t_n}$ for $n\ge n_0$, then $U\subset\chi(A)$.

Obviously, the map $\chi$ is not constant since $\chi$ fixes the unit circle. Thus $\{B_{t_n}\}$ converges to $\chi(A)$ as $n\to\infty$. It follows that
$\chi(A)=\{z: r<|z|<1/r\}$ with $r=\lim_{n\to\infty}r_{t_n}$. Therefore,
$$
\lim_{n\to\infty}\md\phi_{t_n}(A)=\lim_{n\to\infty}\md B_{t_n}=\md\chi(A)=\md A.
$$
\qed

\subsection{Quotient maps}

A continuum $E\subset\cbar$ is {\bf full} if $\cbar\smm E$ is connected. By a {\bf quotient map}\/ of $\cbar$ we mean a continuous onto map $q$ of $\cbar$ such that for any point $w\in\cbar$, $q^{-1}(w)$ is either a single point or a full continuum.

We will call two quotient maps $q_0$ and $q_1$ of $\cbar$ {\bf isotopic} rel a closed subset $E\subset\cbar$ if there exists a continuous map $H: [0,1]\times\cbar\to\cbar$ such that $H(0,\cdot)=q_0$, $H(1,\cdot)=q_1$, $q_t=H(t,\cdot)$ is a quotient map of $\cbar$ for $t\in [0,1]$ and $q_t^{-1}(w)=q_0^{-1}(w)$ for $w\in q_0(E)$ and $t\in [0,1]$. Refer to \cite{CPT} or \cite{PM} for the following lemma.

\REFLEM{quotient} Let $q$ be a continuous onto map of $\cbar$. The following conditions are equivalent.

(i) The map $q$ is a quotient map.

(ii) $q^{-1}(E)$ is a full continuum if $E\subset\cbar$ is a full continuum.

(iii) $q^{-1}(E)$ is a continuum if $E\subset\cbar$ is a continuum.

(iv) $q^{-1}(U)$ is an $n$-connected domain if $U\subset\cbar$ is an $n$-connected domain ($n\ge 1)$.

(v) There exists a sequence of homeomorphisms $\{\phi_n\}$ of $\cbar$ such that $\phi_n$ converges uniformly to $q$ as $n\to\infty$.
\ENDLEM

The following lemma will be used  in the proof of Lemma \ref{normal-q}.

\REFLEM{equi-p} Let $\sQ$ be a family of quotient maps of $\cbar$. Then $\sQ$ is equicontinuous if, for any point $w_0\in\CC$ and any disk $U\subset\cbar$ with $w_0\in U$, there exist a constant $\delta(w_0)>0$ and a disk $V\ni w_0$ with $\overline{V}\subset U$, such that for any $q\in\sQ$, the spherical distance
$$
\text{\rm dist}_s(q^{-1}(\partial U), q^{-1}(V))>\delta(w_0).
$$
\ENDLEM

\beginp If $\sQ$ is not equicontinuous, then for any $\varepsilon>0$, there exists a sequence $\{q_n\}$ in $\sQ$ and a sequence of pairs $(z_n, z'_n)$ of points in $\cbar$ such that $\text{\rm dist}_s(z_n, z'_n)\to 0$ as $n\to\infty$ but
$$
\text{\rm dist}_s(q_n(z_n), q_n(z'_n))>\varepsilon.\eqno{(2.2)}
$$
Passing to a subsequence, we assume that $\{q_n(z_n)\}$ converges to a point $w_0$. Let $U=\DD(w_0, \varepsilon/2)$. By hypothesis, there exist a constant $\delta(w_0)>0$ and a disk $V\ni w_0$ with $\overline{V}\subset U$, such that
$$
\text{\rm dist}_s(q_n^{-1}(\partial U), q_n^{-1}(V))>\delta(w_0).
$$

When $n$ is large enough, $q_n(z_n)\in V$. Thus $q_n(z'_n)\notin U$ by (2.2). Therefore,
$$
\delta(w_0)<\text{\rm dist}(q_n^{-1}(\partial U), q_n^{-1}(V))\le\text{\rm dist}(z_n,z'_n).
$$
This contradicts the fact that $\text{\rm dist}_s(z_n, z'_n)\to 0$ as $n\to\infty$.
\qed

\subsection{Convergence of rational map sequences}

If a sequence of rational maps $\{f_n\}$ is uniformly convergent on $\cbar$, then it converges to a rational map $g$. Moreover, $\deg(f_n)=\deg(g)$ once $n$ is large enough. The next lemma will be used in \S6.2 and \S9.2.

\REFLEM{rational-s} Let $\{f_n\}$ be a sequence of rational maps with constant degree $d\ge 1$. Suppose that $U\!\subset\!\cbar$ is a non-empty open set and $\{f_n\}$ converges uniformly to a map $g$ on $U$ as $n\to\infty$. Then $g$ is a rational map and $\deg g\le d$. Moreover, $\deg g=d$ implies that $\{f_n\}$ converges uniformly to $g$ on $\cbar$ as $n\to\infty$.
\ENDLEM

\beginp
By composing M\"{o}bius transformations, we may assume that $\infty\in U$ and $f_n(\infty)\to 1$. Thus, once $n$ is large enough, the function $f_n$ has the form
$$
f_n(z)=k_n\dfrac{(z-a_{1,n})\cdots (z-a_{d,n})}{(z-b_{1,n})\cdots (z-b_{d,n})},
$$
and $k_n\to 1$ as $n\to\infty$. Since $\{f_n\}$ converges uniformly in a neighborhood of $\infty$ and $f_n(\infty)\to 1$, both $\{a_{i,n}\}$ and $\{b_{j,n}\}$ are bounded in $\CC$. Passing to a subsequence $\{f_{n_k}\}$, we have
$$
(a_{1,{n_k}},\cdots, a_{d, {n_k}}; b_{1,{n_k}},\cdots, b_{d,{n_k}})\to (a_1,\cdots,a_d; b_1,\cdots, b_d)\text{ as }n_k\to\infty.
$$
If $a_i\ne b_j$ for $1\le i,j\le d$, let
$$
g_1(z)=\dfrac{(z-a_1)\cdots(z-a_d)}{(z-b_1)\cdots(z-b_d)};
$$
then $\{f_{n_k}\}$ converges uniformly to $g_1$ on $\cbar$. Since $\{f_n\}$ converges uniformly to $g$ on $U$, we have $g=g_1$ on $U$, and $\{f_n\}$ converges uniformly to $g$ on $\cbar$.

Otherwise, suppose $a_i\ne b_j$ for $1\le i,j\le d_0$ and $a_k=b_k$ for $d_0<k\le d$. Let
$$
g_1(z)=\dfrac{(z-a_1)\cdots(z-a_{d_0})}{(z-b_1)\cdots(z-b_{d_0})}.
$$
Then $\{f_{n_k}\}$ converges uniformly to $g_1$ on any compact subset of $\cbar\smm\{a_{d_0+1},\cdots, a_d\}$. So $g=g_1$ on $U$ and hence $g$ is a rational map with $\deg g=d_0<d$.
\qed

\subsection{Shrinking Lemma}

Let $f$ be a rational map with $\deg f>1$. The following lemma is well-known (refer to \cite{LM}) and will be used several times throughout this paper.

\REFLEM{shrinking} Let $U\subset\cbar$ be a simply-connected domain disjoint from $\PPP_f$ and rotation domains. For any domain $D\Subset U$ and any integer $n\ge 1$, denote by $C_n$ the maximum of the spherical diameters of all the components of $f^{-n}(D)$.  Then $C_n\to 0$ as $n\to\infty$.
\ENDLEM

\beginp If $\# f^{-1}(\PPP_f)<3$, then $f$ is holomorphically conjugate to the power map $z\mapsto z^{\pm d}$. The lemma is easy to check in this case.

Now we assume that $\# f^{-1}(\PPP_f)\ge 3$. If $C_n\not\to 0$ as $n\to\infty$, there exists a sequence $\{n_k\}$ in $\NN$ with $n_k\to\infty$ as $k\to\infty$, and there exists a component $D_k$ of $f^{-n_k}(D)$, such that the spherical diameter
$$
\text{\rm diam}_s D_k\to C>0\text{ as }k\to\infty.
$$
Denote by $U_k$ the component of $f^{-n_k}(U)$ containing $D_k$ and by $g_k$ the inverse map of the univalent map $f^{n_k}:\, U_k\to U$. Then $\{g_k\}$ is a normal family on $U$ since $g_k(U)=U_k$ is disjoint from $f^{-1}(\PPP_f)$. Thus there is a subsequence of $\{g_k\}$, which we will also denote by $\{g_k\}$, which converges uniformly to a holomorphic function $g$ on $D$. Moreover, $g$ is not a constant since $\text{\rm diam}_s g_k(D)\to C>0$ as $k\to\infty$. Thus $g: D\to g(D)$ is a univalent map.

We claim that $g(D)\subset\FFF_f$. Otherwise, there is a domain $V\Subset g(D)$ such that $V\cap\JJJ_f\ne\emptyset$. So $f^n(V)$ covers $\cbar$ with at most two exceptional points once $n$ is large enough. This contradicts the fact that $f^n(V)\subset D$ for infinitely many $n\in\NN$. Thus $g(D)\subset\FFF_f$.

Since $U$ is disjoint from rotation domains, for any domain $W\Subset g(D)$, $f^n(W)$ converges to a periodic orbit in $\PPP_f$. On the other hand, once $k$ is large enough, $W\Subset g_k(D)$ and $f^{n_k}|_W=g^{-k}|_W$ converges uniformly to the univalent map $g^{-1}$ on $W$. This is a contradiction. Thus $C_n\to 0$ as $n\to\infty$.
\qed

\subsection{Parabolic points}

By a {\bf parabolic fixed point} we mean $(g,y)$ where $g$ is a holomorphic map from a neighborhood of $y\in\cbar$ into $\cbar$ with $g(y)=y$, such that if $y\neq\infty$, $g'(y)=e^{2\pi i\frac{p}{q}}$ where $p$ and $q$ are co-prime positive integers with $p\le q$, and
$$
g^{q}(z)-y=(z-y)(1+c_{kq}(z-y)^{kq}+\cdots)
$$
with $c_{kq}\neq 0$ for some integer $k\ge 1$; or if $y=\infty$, $f(z)=1/g(1/z)$ satisfies the above conditions at the origin. Its rotation number is $p/q$ and its multiplicity is $kq+1$. Refer to \cite{Mi}, \S10, for the following results, with a little modification.

\vskip 0.24cm

{\bf Attracting/repelling petals and flowers}.
Let $(g,y)$ be a parabolic fixed point with rotation number $p/q$ and multiplicity $kq+1$. Let $N\subset\cbar$ be a neighborhood of the point $y$ such that $g$ is injective on $N$. Suppose that $V_i\Subset N$ $(i=1,\cdots,kq)$ are pairwise disjoint disks such that their union $\VVV$ satisfies the following conditions:

(a) $g(\overline{\VVV})\subset\VVV\cup\{y\}$.

(b) $\{g^n(z)\}\to y$ as $n\to\infty$ uniformly on any compact subset of $\VVV$.

(c) If $\{g^n(z)\}\to y$ as $n\to\infty$ for $z\in N\smm\{y\}$, then $g^n(z)\in\VVV$ once $n$ is large enough.

We will call these domains $V_i$ {\bf attracting petals}\/ and their union $\VVV$ an {\bf attracting flower} of $(g,y)$. A {\bf repelling petal} and a {\bf repelling flower}\/ of $g$ are defined as an attracting petal and an attracting flower of $g^{-1}$, respectively.

\REFTHM{petal}{\bf (Leau-Fatou Flower Theorem)} Let $(g,y)$ be a parabolic fixed point. Then there exist an attracting flower $\VVV$ and a repelling flower $\VVV'$ such that $\VVV\cup\VVV'\cup\{y\}$ is a neighborhood of the point $y$.
\ENDTHM

\REFTHM{cylinder}{\bf (Cylinder Theorem)} Let $\VVV$ be an attracting flower of the parabolic fixed point $(g,y)$. Then the quotient space $\VVV/\langle g\rangle$ is a disjoint union of $k$ cylinders.
\ENDTHM

Denote $\sC_a=\VVV/\langle g\rangle$. We will call it the {\bf attracting cylinder} of $(g,y)$. The {\bf repelling cylinder} of $(g,y)$ is defined to be the quotient space $\sC_r=\VVV'/\langle g\rangle$ of a repelling flower $\VVV'$. Then $\sC_r$ is also a disjoint union of $k$ cylinders. The following corollary is easy to prove from the fact that every univalent map from a cylinder to itself is surjective.

\REFCOR{cylinder1}
Let $\VVV\Subset N$ be a disjoint union of $kq$ disks satisfying conditions (a) and (b) in the definition of attracting flowers. If the quotient space $\VVV/\langle g\rangle$ is a disjoint union of $k$ cylinders, then $\VVV$ is an attracting flower.
\ENDCOR

The existence of an attracting flower is a complete characterization of a parabolic fixed point, by the following lemma.

\REFLEM{parabolic} Let $N\subset\CC$ be a domain and let $g:\, N\to\CC$ be a univalent map with a fixed point $y\in N$. Suppose that $\VVV\Subset N$ is a finite disjoint union of disks with $y\notin\VVV$ satisfying the following conditions:

(a) $g(\overline{\VVV})\subset\VVV\cup\{y\}$.

(b) $\{g^n(z)\}\to y$ as $n\to\infty$ uniformly on any compact subset of $\VVV$.

(c) If $\{g^n(z)\}\to x$ as $n\to\infty$ for $z\in N\smm\{x\}$, then $g^n(z)\in\VVV$ once $n$ is large enough.

Then $(g, y)$ is a parabolic fixed point.
\ENDLEM

\beginp There exist a component $V$ of $\VVV$ and an integer $p\ge 1$ such that $\overline{g^p(V)}\subset V\cup\{y\}$, by (a). Pick a point $z\in\partial V\smm\{y\}$ and a Jordan arc $\delta_0$ in $V\smm\overline{g^p(V)}$ joining the point $z$ with $g^p(z)$. Let $\delta_k=g^{kp}(\delta_0)$ for $k\ge 1$. Then $\delta_k$ converges to the point $y$ as $k\to\infty$, by (b). Let $\gamma=\cup_{k\ge 1}\delta_k\cup\{g^{kp}(z)\}$. Then $\gamma$ is an arc joining the point $g^{p}(z)$ with the point $y$, and $g^p(\gamma)\subset\gamma$. Thus by the Snail lemma (refer to \cite{Mi}) the fixed point $y$ is not a Cremer point. Obviously, it is neither repelling nor Siegel. Consequently, it is parabolic or attracting.

If the fixed point is attracting, then there exists a disk $D\Subset N$ with $y\in D$ such that $g(D)\Subset D$. Let $W_0=\VVV\cap g(D)$. Then $g(W_0)\subset W_0$. Let $W_1=g^{-1}(W_0)\cap g(D)$. Then $g(W_1)=W_0\subset W_1$. Inductively, let $W_{n+1}=g^{-1}(W_n)\cap g(D)$ for $n\ge 1$. Then $g(W_{n+1})=W_{n}\subset W_{n+1}$.

Let $E_n=\partial g(D)\smm g^{-1}(W_n)$ for $n\ge 1$. Then $E_n$ is closed and $E_{n+1}\subset E_n$. In particular, $E_n$ is non-empty. Otherwise, we would have $\partial g(D)\subset g^{-1}(W_n)$. It follows that $g(\partial g(D))\subset W_n$. Since $g(W_n)\subset W_{n-1}$ for all $n\ge 1$, we obtain  $g^{n+1}(\partial g(D))\subset W_0\subset\VVV$. Since each component of $\VVV$ is simply connected and $g(D)\Subset D$, we get $g^{n+2}(D)\subset\VVV$.  This is a contradiction since $y\in g^{n+2}(D)$ but $y\not\in\VVV$.

Set $E_{\infty}=\cap_{n\ge 1}E_n$. It is non-empty. Pick a point $z\in E_{\infty}$. Then $z\in\partial g(D)$ and hence $g^n(z)\in g(D)$ for all $n\ge 1$, and
$\{g^n(z)\}\to y$ as $n\to\infty$. On the other hand, $z\not\in g^{-1}(W_n)$ for all $n\ge 1$. Thus $g(z)\in g(D)\smm W_n$. We claim that $g^n(z)\not\in W_0$ for all $n\ge 1$. Otherwise, we would have $g^{n-1}(z)\in W_1$ by the definition of $W_1$. Inductively, $g(z)\in W_{n-1}$ and hence $g(z)\in W_n$. This is a contradiction. By the claim and the definition of $W_0$, $g^n(z)\notin\VVV$ for all $n\ge 1$. This contradicts condition (c). Therefore $(g, y)$ is a parabolic fixed point.
\qed

\vskip 0.24cm

{\bf Sepals, calyxes and horn maps}. Let $(g,y)$ be a parabolic fixed point with rotation number $p/q$ and multiplicity $kq+1$. Let $N\subset\cbar$ be a neighborhood of the point $y$ such that $g$ is injective on $N$.

\REFTHM{sepal}{\bf (Sepal)} There exist $2kq$ disjoint disks $W_j\Subset N$ such that $g$ is a conformal map from their union $\WWW$ onto itself and $\{g^n(z)\}\to y$ as $n\to\infty$ uniformly on any compact subset of $\WWW$.
\ENDTHM

We will call these domains $W_i$ {\bf sepals}\/ and their union $\WWW$ a {\bf calyx}. There are exactly two sepals of $\WWW$ intersecting with an attracting petal; we call them a {\bf left sepal}\/ and a {\bf right sepal}, as viewed from the parabolic fixed point to the attracting petal.

Denote by $\pi_a$ (respectively $\pi_r$) the natural projection from an attracting (repelling) flower to the attracting (repelling) cylinder. For the sake of simplicity, we do not specify the domain of the projection $\pi_a$ or $\pi_r$. This does not cause confusion since if a point is contained in two different attracting flowers (or repelling flowers), then its projection is independent of the choice of flower. Consequently, both $\pi_a$ and $\pi_r$ are well-defined on a calyx.

Let $\WWW$ be a calyx. Then both $\pi_a(\WWW)$ and $\pi_r(\WWW)$ are the disjoint union of $2k$ one-punctured disks. Define
$$
\Upsilon_g:\pi_r(\WWW)\to\pi_a(\WWW)\quad\text{by}\quad\Upsilon_g(\pi_r(z))=\pi_a(z)\text{ for }z\in\WWW.
$$
This is a well-defined conformal map and it is called a {\bf horn map}\/ of $(g,y)$.

\vskip 0.24cm

{\bf Regular flowers}.
Let $(g,y)$ be a parabolic fixed point with rotation number $p/q$ and multiplicity $kq+1$. Let $\VVV$ be an attracting flower of $(g,y)$. Then the projection $\pi_a(\partial\VVV\smm\{y\})$ is a disjoint union of $kq$ arcs and each of the $k$ attracting cylinders contains $q$ such arcs. Each arc lands on punctures at both ends, or else the limit set of each arc is disjoint from punctures, or else the limit sets of the arcs are complicated sets. On the other hand, $\cap_{n=1}^{\infty}g^n(\VVV)$ is empty, or else contains a calyx, or else is a complicated set. To avoid the complexity, we will need a further requirement for attracting flowers in this paper.

An attracting flower $\VVV$ is called {\bf regular} if $\cap_{n=1}^{\infty}g^n(\VVV)$ is empty, or equivalently, each arc in $\pi_a(\partial\VVV\smm\{y\})$ lands on punctures at both ends. A {\bf regular} repelling flower is defined similarly. The following proposition will be used in \S3:

\REFPROP{regular}
(1) Any attracting flower contains a regular attracting flower.

(2) Let $\alpha_n\subset\sC_a$ $(1\le n\le kq)$ be pairwise disjoint arcs connecting two punctures such that each of the $k$ attracting cylinders contains $q$ arcs. Then there exists a regular attracting flower $\VVV$ such that $\pi_a(\partial\VVV\smm\{y\})=\cup_{n=1}^{kq}\alpha_n$.

(3) For any attracting flower $\VVV$ of $(g,y)$, there exists a regular repelling flower $\VVV'$ of $(g,y)$ such that $\VVV\cup\{y\}\cup\VVV'$ is a neighborhood of the point $y$.
\ENDPROP

\beginp
We only prove the proposition in case that $kq=1$. The proof for general case has no essential difficulty.

(1) Let $\VVV$ be an attracting flower of the parabolic fixed point $(g,y)$ and let $\pi_a: \VVV\to\sC_a=\CC^*$ be the natural projection to the cylinder. Pick a bi-infinite sequence $\{w_n\}$ in $\CC^*$ ($n\in\ZZ$) such that $|w_n|<|w_{n+1}|$, $|w_n|\to 0$ as $n\to -\infty$ and $|w_n|\to\infty$ as $n\to\infty$. Let $\beta_n\subset\CC^*$ be a round circle with center zero and radius $|w_n|$ for $n\in\ZZ$. Then $w_n\in\beta_n$.

By condition (c) in the definition of attracting flower, for each point $w\in\beta_n$, there exists a point $z\in\VVV$ such that $\pi_a(z)=w$. Moreover, there exists a disk $D_z\subset\VVV$ with $z\in D_z$ such that $\pi_a$ is injective on $D_z$. Set $U_w=\pi_a(D_z)$. The sets $U_w$ form an open cover of $\beta_n$. Thus there is a finite open sub-cover. This fact shows that there exists an arc $\gamma_n\subset\VVV$ joining a point $z_n\in\pi_a^{-1}(w_n)$ with the parabolic fixed point $y$, such that $\pi_a(\gamma_n)=\beta_n$, and hence $g(\gamma_n)\subset\gamma_n$.

Let $\delta_n\subset\CC^*$ be an arc in the annulus bounded by $\beta_n$ and $\beta_{n+1}$ which joins the point $w_n$ with $w_{n+1}$. As above, there exists an integer $i_n\ge 0$ and an arc $\wt\delta_n\subset\VVV$ which joins the point $g^{i_n}(z_n)$ with $g^{i_n}(z_{n+1})$ such that $\pi_a(\wt\delta_n)=\delta_n$. This implies that the domain bounded by $\wt\delta_n$, $\gamma_n$ and $\gamma_{n+1}$ is contained in $\VVV$.

Start from the point $z_0\in\VVV$. From the existence of $\wt\delta_n$, we know that there exists an integer $k_1\ge 0$ such that the Fatou line segment $\alpha_1$ (a line segment in Fatou coordinates), which joins the point $z_0$ with $g^{k_1}(z_1)$, is contained in $\VVV$. Inductively, there exists a sequence of non-negative integers $\{k_n\}$ ($n\ge 1$) such that the Fatou line segment $\alpha_n$ which joins the point $g^{k_{n-1}}(z_{n-1})$ with $g^{k_n}(z_n)$ is contained in $\VVV$.

Repeating this argument for $n\le -1$ and setting $k_0=0$, we get a sequence of non-negative integers $\{k_n\}$ ($n\le 0$) such that the Fatou line segment $\alpha_n$ which joins the point $g^{k_{n+1}}(z_{n+1})$ with $g^{k_n}(z_n)$ is contained in $\VVV$.

These arcs $\{\alpha_n\}$ for $n\neq 0$ are pairwise disjoint. The union of them together with their endpoints forms an arc $\alpha$ whose endpoints are both the parabolic fixed point $y$. Thus $\alpha\cup\{y\}$ bounds a disk in $\VVV$, which is a regular attracting flower.

\vskip 0.24cm

(2) We assume that $kq=1$ for simplicity. Let $\WWW$ be a calyx of the parabolic fixed point $(g,y)$. Then $\pi_a(\WWW)$ is the disjoint union of two once-punctured disks.  The arc $\alpha_1$ can be cut into three arcs $\gamma_0, \gamma_1$ and $\gamma_2$ such that both $\gamma_0$ and $\gamma_1$ are contained in $\pi_a(\WWW)$ and the endpoints of $\gamma_2$ stay on $\alpha$.

By a similar argument as in the proof of (1), there exists an arc $\wt\gamma_2$ in an attracting flower of $(g,y)$ such that $\pi_a(\wt\gamma_2)=\gamma_2$. Let $\wt\gamma_i\subset\WWW$ be a lift of $\gamma_i$ that has a common endpoint with $\wt\gamma_2$ for $i=0,1$. Then the union of $\wt\gamma_i$ ($i=0,1,2$) together with their endpoints bounds a regular attracting flower $\VVV$ and $\pi_a(\partial\VVV\smm\{y\})=\alpha_1$.

\vskip 0.24cm

(3) Pick a calyx $\WWW$ of $(g,y)$ such that $\partial\WWW\smm\{y\}$ is two horizontal lines in the Fatou coordinate. Then there exists a regular repelling flower $\VVV'_1$ of $(g,y)$ such that $\partial\VVV'_1\cap\WWW$ is two vertical lines in the Fatou coordinate. Denote by $z_0$ and $z'_0$ the two endpoints of these two lines on $\partial\WWW\smm\{y\}$.

Pick a sequence of calyxes $\{\WWW_n\}$ ($n\ge 1$) in $\WWW$ such that $\overline{\WWW_{n+1}}\subset\WWW_n\cup\{y\}$,  $\partial\WWW_n\smm\{y\}$ is two horizontal lines in the Fatou coordinate and $\{\WWW_n\}$ converges to the point $y$.  Then there exists a pair of points $(z_n, z'_n)$ in $\partial\WWW_n\smm\{y\}$ for each $n\ge 1$ such that $z_n$ is contained in the same sepal as $z_0$, $z_n'$ is contained in the same sepal as $z_0'$, and
the Fatou line segments $\delta_n, \delta'_n$ which join $z_n$ with $z_{n+1}$, and $z'_n$ with $z'_{n+1}$, respectively, are contained in $\VVV$. Since $\{\WWW_n\}$ converges to the point $y$, both $\{z_n\}$ and $\{z'_n\}$ converge to $y$ as $n\to\infty$.

Let $\delta_0, \delta'_0$ be the Fatou line segments which join $z_0$ with $z_1$, and $z'_0$ with $z'_1$, respectively. Let $\alpha$ be the union of $\delta_n$ and $\delta'_n$ together with their endpoints for $n\ge 0$. Then $\alpha\cup\{y\}$ bounds a regular repelling flower $\VVV'$ and $\VVV\cup\{y\}\cup\VVV'$ is a neighborhood of the point $y$.
\qed

\subsection{Quotient space of rational maps and periodic arcs}

Let $f$ be a rational map with attracting or parabolic domains. Define $\wt\sR_f\subset\FFF_f$ by $z\in\wt\sR_f$ if its forward orbit $\{f^n(z)\}$ is infinite, disjoint from $\PPP_f$ and contained in parabolic or attracting (but not superattracting) cycles. Define the grand orbit equivalence relation by $z_1\sim z_2$ if $f^n(z_1)= f^m(z_2)$ for integers $n,m>0$. Then the quotient space $\sR_f=\wt\sR_f/\sim$ has only finitely many components. Each of them is either a punctured torus with at least one puncture (but only finitely many punctures) corresponding to an attracting basin, or a punctured sphere with at least three puncture (but only finitely many punctures) corresponding to a parabolic basin (refer to \cite{McS}). We will call $\sR_f$ the {\bf (punctured) quotient space}\/ of $f$ and denote by $\pi_f:\,\wt\sR_f\to\sR_f$ the natural projection.

\vskip 0.24cm

Let $\beta$ be an open arc in the attracting or parabolic domains of $f$ with $\beta\cap\PPP_f=\emptyset$. We call $\beta$ a {\bf periodic arc}\/ if $\beta$ coincides with a component of $f^{-p}(\beta)$ for some integer $p\ge 1$, or an {\bf eventually periodic arc}\/ if $f^k(\beta)$ is a periodic arc for some integer $k\ge 0$. In that case the projection $\pi_f(\beta)$ is a simple closed curve in $\sR_f$.

Conversely, let $\gamma\subset\sR_f$ be a Jordan curve. Then either each component of $\pi_f^{-1}(\gamma)$ is a Jordan curve or each component of $\pi_f^{-1}(\gamma)$ is an eventually periodic arc.

\REFLEM{arc} A periodic arc lands at both ends.
\ENDLEM

\beginp
Let $\beta$ be a periodic arc of a rational map $f$ with period $p\ge 1$. Then $g:=f^p|_{\beta}$ is a homeomorphism on $\beta$. Pick a point $y\in\beta$ and denote by $\alpha_n$ the closed arcs in $\beta$ with endpoints $g^{n}(y)$ and $g^{n+2}(y)$, for each integer $n$. Then $g(\alpha_n)=\alpha_{n+1}$ and $\cup_{n=-\infty}^{\infty}\alpha_n=\beta$.

Since $\beta$ is contained in an attracting or parabolic periodic domain, $\alpha_n$ converges to the attracting or parabolic point $x_0$ in the basin as $n\to\infty$. Thus $\beta$ lands on $x_0$ from one direction.

Denote by $r(\beta)$ the limit set of the other end of $\beta$. It is connected and $f^p(r(\beta))=r(\beta)$. For any point $x\in r(\beta)$, let $\{x_n\}$ ($n\ge 1$) be a sequence of points in $\beta$ converging to $x$. Then for each point $x_n$, there exists an integer $k_n$ such that both $x_n$ and $g(x_n)$ are contained in $\alpha_{k_n}$. Moreover $k_n\to -\infty$ as $n\to\infty$ since $\{x_n\}$ converges to the point $x\in r(\beta)$. By Lemma \ref{shrinking}, $\text{\rm diam}(\alpha_{k})\to 0$ as $k\to -\infty$. Thus $\text{\rm dist}(x_n, g(x_n))\to 0$ as $n\to\infty$. So $f^p(x)=x$. Therefore $r(\beta)$ can only be a single point.
\qed

\vskip 0.24cm

Let $\beta$ be a periodic arc with period $p\ge 1$. Denote by $a(\beta)$ and $r(\beta)$ the limit points of the forward and backward orbits on $\beta$ under $f^p$, respectively. We call $a(\beta)$ the {\bf attracting end}\/ and $r(\beta)$ the {\bf repelling end}\/ of $\beta$. They may coincide if both of them are parabolic.

If $\beta_1$ is a pre-periodic arc, let $k>0$ be an integer such that $f^k(\beta_1)$ is periodic. We will denote by $a(\beta_1)$ and $r(\beta_1)$ the two endpoints of $\beta_1$ so that $f^k(a(\beta_1)) =a(f^k(\beta_1))$ and $f^k(r(\beta_1))=r(f^k(\beta_1))$. Both $a(\beta_1)$ and $r(\beta_1)$ are pre-periodic.

\section{Unicity}

We will prove Theorem \ref{unicity} in this section. It is known that a parabolic fixed point has infinitely many analytic invariants, e.g. the power series representations of the horn maps. This fact makes our attempt to analytically conjugate two parabolic points a delicate issue. Our approach is to modify  c-equivalency to be a local quasiconformal conjugacy with a small distortion.

\subsection{A lemma about quasiconformal maps}

Let $\phi: R\to R'$ be a quasiconformal homeomorphism between Riemann surfaces. We denote by
$$
\mu(\phi)=\mu_{\phi}(z)\dfrac{d\bar z}{dz}=\dfrac{\partial_{\bar z}\phi}{\partial_z\phi},\quad
K_{\phi}(z)=\dfrac{1+|\mu_{\phi}(z)|}{1-|\mu_{\phi}(z)|}\quad\text{and}\quad K(\phi)=\|K_{\phi}\|_{\infty}
$$
the Beltrami differential, the dilatation and the maximal dilatation of $\phi$,  respectively.

A quasiconformal map $\phi$ is {\bf extremal} if $K(\phi)\le K(\psi)$ for all the quasiconformal maps $\psi$ isotopic to $\phi$ rel the boundary. There always exists an extremal quasiconformal map in the isotopy class of a quasiconformal map.

A quasiconformal map $\phi$ is called a {\bf Teichm\"{u}ller map} associated with an integrable holomorphic quadratic differential $\omega(z)dz^2$ if
$$
\mu_{\phi}(z)=k\frac{\overline{\omega(z)}}{|\omega(z)|}
$$
for some constant $0<k<1$. It is known that a Teichm\"{u}ller map is the unique extremal quasiconformal map in its isotopy class. Refer to \cite{Str} or \cite{RS} for the following theorems.

\REFTHM{Strebel}{\bf (Strebel)}
Let $\phi_0: R\to R'$ be an extremal quasiconformal map between open Riemann surfaces with $K(\phi_0)>1$. If there exists a quasiconformal map $\phi$ isotopic to $\phi_0$ rel the boundary such that $K_{\phi}(z)<K(\phi_0)$ in some neighborhood of $\partial R$, then $\phi_0$ is a Teichm\"{u}ller map associated with an integrable holomorphic quadratic differential and hence is the unique extremal quasiconformal map in its isotopy class.
\ENDTHM

\REFTHM{MI}{\bf (Main inequality)} Let $\phi, \psi: R\to R'$ be quasiconformal maps between open Riemann surfaces which are isotopic rel the boundary. Let $\omega(z)dz^2$ be an integrable holomorphic quadratic differential on $R$. Then
$$
\|\omega\|:=\int_R|\omega(z)|dxdy\le\int_R|\omega(z)|\dfrac{\left|1-\mu_{\phi}\dfrac{\omega(z)}
{|\omega(z)|}\right|^2}{1-|\mu_{\phi}|^2}K_{\psi^{-1}}\circ\phi(z)dxdy.
$$
\ENDTHM

From these known results we obtain:

\REFLEM{qc} Let $\phi$ be a quasiconformal map of $\CC^*$. Then for any $\epsilon>0$ and any $r_0>1$, there exist a constant $r>r_0$ and a quasiconformal map $\psi$ of $\CC^*$  such that:

(1) $\psi=\phi$ on $E(r):=\cbar\smm\AA(r)$,

(2) $\psi|_{\AA(r)}$ is isotopic to $\phi|_{\AA(r)}$ rel the boundary, and

(3) $K(\psi)<K(\phi|_{E(r_0)})+\epsilon$.
\ENDLEM

\beginp
For any $r>r_0$ , let $\psi_r: \AA(r)\to\phi(\AA(r))$ be an extremal quasiconformal map isotopic to $\phi|_{\AA(r)}$ rel the boundary. If
$$
K(\psi_r)<K(\phi|_{E(r_0)})+\epsilon
$$
for some $r>r_0$, set $\psi=\psi_r$ on $\AA(r)$ and $\psi=\phi$ on $E(r)$; then $\psi$ satisfies the conditions.

Now we assume
$$
K(\psi_r)\ge K(\phi|_{E(r_0)})+\epsilon
$$
for $r>r_0$. By Theorem \ref{Strebel}, $\psi_r$ is a Teichm\"{u}ller map associated with an integrable holomorphic quadratic differential $\omega_r(z)dz^2$. We may assume $\|\omega_r\|=1$. Then $\omega_r$ converges to zero uniformly on any compact subset of $\CC^*$ as $r\to\infty$ since there is no non-zero integrable holomorphic quadratic differential on $\CC^*$.

Let $U_r=\psi_r^{-1}\circ\phi(\AA(r_0))$. Then there exists a compact subset $V\subset\CC^*$ such that $U_r\subset V$ for all $r>r_0$, since $K(\psi_r)\le K(\phi)$. Since $\omega_r$ converges uniformly to zero on $V$ as $r\to\infty$, we obtain
$$
\int_{U_r}|\omega_r|dxdy < \dfrac{\epsilon}{K(\phi)}
$$
when $r$ is large enough. Applying Theorem \ref{MI} for $\phi|_{\AA(r)}$, $\psi_r$ and $\omega_r$, we get

\begin{align*}
1=\|\omega\| & =\int_{\AA(r)}|\omega_r|dxdy \\
& \le\int_{\AA(r)}|\omega_r|\dfrac{\left|1-\mu_{\psi_r}\dfrac{\omega_r}{|\omega_r|}\right|^2}{1-|\mu_{\psi_r}|^2}
           K_{\phi^{-1}}\circ\psi_r dxdy \\
& =\int_{\AA(r)}\dfrac{|\omega_r|}{K(\psi_r)}K_{\phi^{-1}}\circ\psi_r dxdy.
\end{align*}
Thus
\begin{align*}
K(\psi_r) & \le\int_{\AA(r)}|\omega_r|K_{\phi^{-1}}\circ\psi_r dxdy \\
          & \le\int_{U_r}|\omega_r|K_{\phi^{-1}}\circ\psi_r dxdy
           +\ds\int_{\AA(r)\smm U_r}|\omega_r|K_{\phi^{-1}}\circ\psi_r dxdy.
\end{align*}
Note that $\psi_r(\AA(r)\smm U_r)\subset\phi(E(r_0))$. Therefore when $r$ is large enough,
$$
K(\psi_r) <\epsilon+K(\phi|_{E(r_0)}).
$$
This is a contradiction.
\qed

\subsection{Local conjugacy between parabolic points}

\REFLEM {p2l}{\bf (From petal conjugacy to local conjugacy)} Let $(f,x)$ and $(g, y)$ be parabolic fixed points. Let $\phi: \VVV(f)\to\VVV(g)$ be a $K$-quasiconformal conjugacy between regular attracting flowers of $(f,x)$ and $(g, y)$. Then for any $\epsilon>0$, there exist a neighborhood $N$ of the point $x$ with $\VVV(f)\subset N$ and a $(K+\epsilon)$-quasiconformal map $\phi_0$ on $N\cup f(N)$ such that $\phi_0=\phi$ on $\VVV(f)$ and $\phi_0\circ f=g\circ\phi_0$ on $N$.
\ENDLEM

\beginp
Let $\WWW(f)$ be a calyx at the parabolic fixed point $(f,x)$ such that each arc of $\partial\WWW(f)\smm\{x\}$ intersects $\partial\VVV(f)$ at exactly one point. Once the calyx $\WWW(f)$ is small enough, $\phi(\WWW(f)\cap\VVV(f))$ is contained in a calyx $\WWW_1(g)$ of $(g,y)$. Thus the union of the backward orbit of $\phi(\WWW(f)\cap\VVV(f))$ under $g|_{\WWW_1(g)}$ forms a calyx $\WWW(g)$ of $(g,y)$. The map $\phi$ can be extended to a quasiconformal map from $\WWW(f)$ to $\WWW(g)$ by the equation  $\phi\circ f=g\circ\phi$.

By proposition \ref{regular} (2), there exists a regular repelling flower $\VVV'(f)$ of $(f,x)$, disjoint from $\VVV(f)$, such that each arc of $\partial\VVV'(f)\smm\{x\}$ intersects $\partial\WWW(f)$ at exactly two points and there exists an integer $k\ge 1$ such that
$$
f^k(\partial\VVV'(f)\cap\WWW(f))=\partial\VVV(f)\cap\WWW(f).
$$
Similarly, there exists a regular repelling flower $\VVV'(g)$ of $(g, y)$, disjoint from $\VVV(g)$, such that each arc of $\partial\VVV'(g)\smm\{y\}$ intersects $\partial\WWW(g)$ at exactly two points and there exists an integer $k_1\ge 1$ such that
$$
g^{k_1}(\partial\VVV'(g)\cap\WWW(g))=\partial\VVV(g)\cap\WWW(g).
$$
We may assume that $k_1=k$. Otherwise, if $k_1>k$, we may use $f^{k-k_1}(\VVV'(f))$ to replace $\VVV'(f)$.

Denote by $(\pi_{a,f}, \pi_{r,f})$ and $(\pi_{a,g}, \pi_{r,g})$ the projections to attracting and repelling cylinders $(\sC_a(f), \sC_{r}(f))$ of $(f,x)$ and attracting and repelling cylinders $(\sC_{a}(g), \sC_{r}(g))$ of $(g,y)$, respectively. Denote by
$$
\Upsilon_f:\pi_{r,f}(\WWW(f))\to\pi_{a,f}(\WWW(f))\quad\text{and}\quad \Upsilon_g:\pi_{r,g}(\WWW(g))\to\pi_{a,g}(\WWW(g))
$$
the horn maps. Let $\Phi:\, \sC_{a}(f)\to\sC_{a}(g)$ be the projection of $\phi$ to the attracting cylinders. It is well-defined since $\phi$ is a conjugacy. Composing with the horn maps $\Upsilon_f$ and $\Upsilon_g$, we get a $K$-quasiconformal map
$$
\Psi:=\Upsilon_g^{-1}\circ\Phi\circ\Upsilon_f: \pi_{r,f}(\WWW(f))\to\pi_{r,g}(\WWW(g)).
$$
By the choices of $\VVV'(f)$ and $\VVV'(g)$, we have
$$
\Psi\circ\pi_{r,f}(\partial\VVV'(f)\cap\WWW(f))=\pi_{r,g}(\partial\VVV'(g)\cap\WWW(g)).
$$
Note that the boundaries of $\pi_{r,f}(\WWW(f))$ and $\pi_{r,g}(\WWW(g))$ are disjoint unions of simple closed curves. Thus the map $\Psi$ can be extended continuously to a homeomorphism $\Psi:\, \sC_{r}(f)\to\sC_{r}(g)$ such that
$$
\Psi\circ\pi_{r,f}(\partial\VVV'(f)\smm\{x\})=\pi_{r,g}(\partial\VVV'(g)\smm\{y\}).
$$

\begin{figure}[htbp]
\begin{center}
\includegraphics[scale=0.8]{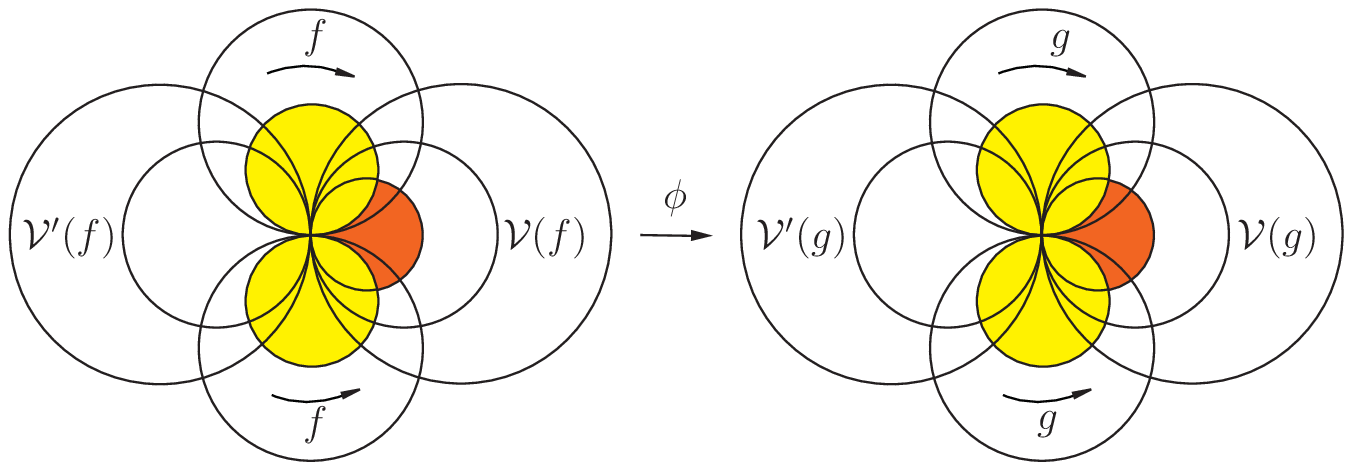}
\end{center}
\begin{center}{\sf Figure 3. From petal to local conjugacy.}
\end{center}
\end{figure}

By Lemma \ref{qc}, for any $\epsilon>0$, there exist a smaller calyx $\WWW_0(f)\subset\WWW(f)$ and a ($K+\epsilon$)-quasiconformal map $\Psi_0:\, \sC_{r}(f)\to\sC_{r}(g)$ such that $\Psi_0=\Psi$ on $\pi_{r,f}(\WWW_0(f))$, and restricted to the complement of $\pi_{r,f}(\WWW_0(f))$, $\Psi_0$ is isotopic to $\Psi$ modulo the boundary. Let $\psi_0$ be the lift of $\Psi_0$ through the isotopy from $\Psi_0$ to $\Psi$. It is well-defined on $\VVV'_n(f)=f^{-n}(\VVV'(f))$ for some integer $n\ge 0$, and $\psi_0(\VVV'_{n}(f))\subset\VVV'(g)$.

Because $\Psi_0=\Psi$ on $\pi_{r,f}(\WWW_0(f))$, and restricted to the complement of $\pi_{r,f}(\WWW_0(f))$, $\Psi_0$ is isotopic to $\Psi$ rel the boundary, we have $\psi_0=\phi$ on $\VVV'_n(f)\cap\WWW_0(f)$. Set $\phi_0=\phi$ on $\VVV(f)\cup\WWW_0(f)$ and $\phi_0=\psi$ on $\VVV'_{n}(f)$. Then $\phi_0$ is a ($K+\epsilon$)-quasiconformal map and $g\circ\phi_0=\phi_0\circ f$ on $\VVV(f)\cup\WWW(f)\cup\VVV'_{n+1}(f)$.
\qed

\vskip 0.24cm

The following result, which was first proved by McMullen in \cite{Mc2}, is a corollary of the previous lemma.

\REFCOR {1+e} Two parabolic fixed points with same rotation number and multiplicity are locally $(1+\epsilon)$-quasiconformal conjugate for any $\epsilon>0$.
\ENDCOR

\subsection{Proof of Theorem \ref{unicity}}

Let $F$ and $G$ be semi-rational maps and let $(\phi, \psi)$ be a c-equivalence between them on a fundamental set $\UUU$ of $F$, i.e., $(\phi, \psi)$ is a pair of orientation-preserving homeomorphisms of $\cbar$ such that:

(a) $\phi\circ F=G\circ\psi$,

(b) $\phi$ is holomorphic in $\UUU$, and

(c) $\psi=\phi$ in $\UUU\cup\PPP_F$ and $\psi$ is isotopic to $\phi$ rel $\UUU\cup\PPP_F$.

\vskip 0.24cm

The relation c-equivalence is an equivalence relation by the following Lemma.

\begin{lemma}\label{equivalent}
$\phi(\UUU)$ is a fundamental set of $G$.
\end{lemma}

\beginp
Obviously, $\phi(\UUU)\subset G^{-1}(\phi(\UUU))$ and $\phi(\UUU)$ contains all the attracting or super-attracting points in $\PPP'_G$.  We only need to prove that if $x\in\PPP'_F$ is a parabolic periodic point of $F$ with period $p\ge 1$ and $\VVV\subset\UUU$ is an attracting flower of $F^p$ at the point $x$, then $\phi(x)$ is a parabolic periodic point of $G$ and $\phi(\VVV)$ is an attracting flower of $G^p$ at the point $\phi(x)$.

By Theorem \ref{cylinder}, the quotient space $\VVV/\langle F^p\rangle$ is a disjoint union of cylinders. If $\phi(x)$ is attracting, then $\phi$ induces a holomorphic injection from the cylinders into a torus. This is impossible. So $\phi(x)$ is parabolic. Since each attracting petal of $G^p$ at the point $\phi(x)$ contains points of $\PPP_G$, the multiplicities of the parabolic fixed points $(G^p, \phi(x))$ and $(F^p, x)$ are equal. The quotient space $\phi(\VVV)/\langle G^p\rangle$ is also a disjoint union of cylinders. By Corollary \ref{cylinder1}, $\phi(\VVV)$ is an attracting flower of $G^p$ at the point $\phi(x)$.
\qed

\vskip 0.24cm

The following totally topological lemma will be used when we deal with c-equivalence. It implies that c-equivalence may be defined with a weaker condition;  condition (c) can be replaced by: $\psi=\phi$ in $\UUU\cup\PPP_F$ and $\psi$ is isotopic to $\phi$ rel $\PPP_F$.

\begin{lemma}\label{isotopy}
Let $E_1, E_2\subset$ be closed subsets of $\cbar$. Let $\theta$ be an orientation-preserving homeomorphism of $\cbar$ isotopic to the identity rel $E_1$ with $\theta=\text{\rm id}$ on $E_2$. Assume that

(1) $(E_1\smm E_2)$ is a finite set,

(2) $E_2$ has only finitely many components and each of them contains points of $E_1$, and

(3) Each component of $E_2$ is a either a closed Jordan disk, or finitely many closed Jordan disks intersecting at a single point.

Then $\theta$ is isotopic to the identity rel $E_1\cup E_2$.
\end{lemma}

\beginp
Assume $E_2\neq\emptyset$. Otherwise the lemma is trivial. Let $H: I\times\cbar\to\cbar$, where $I=[0,1]$, be an isotopy such that $H(0,\cdot)=\text{\rm id}$, $H(1, \cdot)=\theta$ and $H(t, z)=z$ for $z\in E_1$ and $t\in I$. Pick one point of $E_1$ in each component of $E_2$, and denote by $E_0$ the set of them together with all points of $E_1$ outside of $E_2$. Then $E_0\subset E_1$ and $E_0\cup E_2=E_1\cup E_2$.

Consider the path $H(I, z)\subset\cbar\smm E_0$ for each point $z\in U:=\cbar\smm(E_0\cup E_2)$. Its two endpoints $(z,\theta(z))$ are contained in $U$. Since each component of $E_2$ contains exactly one point of $E_0$, there exists a path $\beta(z)\subset U$ connecting $(z, \theta(z))$ which is homotopic to $H(I, z)$ in $\cbar\smm E_0$, and such paths are homotopic to each other in $U$.

Let $\gamma(z)\subset U$ be the unique geodesic under the Poincar\'{e} metric on $U$ connecting $(z,\theta(z))$ and homotopic to $\beta(z)$ in $U$. Define $\theta_t(z)$ by $\theta_0(z)=z$, $\theta_1(z)=\theta(z)$ and $\theta_t(z)\in\gamma(z)$ with
$$
\frac{L(\theta(0,t))}{L(\theta(t,1))}=\frac{t}{1-t},
$$
where $L(\cdot)$ denotes the length under Poincar\'{e} metric. Then $\theta_t$ is a homotopy in $U$ connecting the identity with $\theta$ and $\theta_t=\text{\rm id}$ on the boundary of $U$ for $t\in I$.

By Theorem 1.12 in \cite{BD}, two homotopic homeomorphisms between compact surfaces with finitely many punctures and holes are isotopic. Thus there exists an isotopy $\wt\theta_t$ on $U$ connecting the identity with $\theta$ such that $\wt\theta_t=\text{\rm id}$ on $\partial U$ for $t\in I$. Define $\wt\theta_t=\text{\rm id}$ on $E_1\cup E_2$. Then $\wt\theta_t$ is an isotopy of $\cbar$ rel $E_1\cup E_2$ connecting the identity with $\theta$.
\qed

\begin{lemma}\label{c2q} Let $F$ and $G$ be semi-rational maps. Let $(\phi, \psi)$ be a c-equivalence from $F$ to $G$. Then for any $\epsilon>0$, there exist a c-equivalence $(\phi_0, \psi_0)$ between $F$ and $G$ in the isotopy class of $\phi$ and an open set $U\supset\PPP'_F$ such that $\psi_0=\phi_0$ on $U$, $\phi_0$ is quasiconformal in $\cbar$ and $K(\phi_0|_U)<1+\epsilon$.
\end{lemma}

\beginp
Let $x\in\PPP'_F$ be a parabolic periodic point of $F$ with period $p\ge 1$. Assume $p=1$ for the simplicity. Let $\wt\VVV_x$ be an attracting flower of $(F, x)$ such that $\phi$ is holomorphic in $\wt\VVV_x$ and $\psi$ is isotopic to $\phi$ rel $\wt\VVV_x$. By Lemma \ref{p2l}, there exist an attracting flower $\VVV_x$ of $(F, x)$ with $\VVV_x\subset\wt\VVV_x$, a disk $D_x$ with $\VVV_x\Subset D_x$ and a $(1+\epsilon)$-quasiconformal map $\xi_x$ defined on a domain containing the closure of $D_x\cup F(D_x)$ such that $\xi_x=\phi$ on $\VVV_x$ and $\xi_x\circ F=G\circ\xi_x$ on $D_x$.

Let $\wt U_x\subset\cbar$ be a disk such that $\VVV_x\Subset\wt U_x$ and $\wt U_x\smm\overline{\VVV_x}$ is disjoint from $\PPP_F$. By choosing $D_x$ close enough to $\VVV_x$, one may assume that $D_x\cup\phi^{-1}\circ\xi_x(D_x)\Subset\wt U_x$ since $\phi^{-1}\circ\xi_x=\text{\rm id}$ on $\VVV_x$. Thus there exists a homeomorphism $\theta_x$ of $\cbar$ isotopic to the identity rel $(\cbar\smm\wt U_x)\cup\overline{\VVV_x}$ such that $\theta_x=\phi^{-1}\circ\xi_x$ on $D_x\smm\overline{\VVV_x}$.

Choose $\theta_x$ as above for all parabolic periodic points $x\in\PPP'_F$. Let $\theta$ be the composition of all these maps and let $\phi_1=\phi\circ\theta$. Then there exists a small fundamental set $\UUU_1$ of $F$ such that $\phi_1$ is isotopic to $\phi$ rel $\UUU_1\cup\PPP_F$. Let $\psi_1$ be the lift of $\phi_1$, i.e. $\phi_1\circ F=G\circ\psi_1$. Then $(\phi_1, \psi_1)$ is also a c-equivalence from $F$ to $G$ and there exists an open set $\wt U\supset\PPP'_F$ such that $\psi_1=\phi_1$ on $\wt U$ and $K(\phi_1|_{\wt U})<1+\epsilon$.

Pick a quasi-disk $U_x\Subset\wt U$ with $x\in U_x$ for each point $x\in\PPP'_F$ such that they have disjoint closures. Then their images under $\phi_1$ are also quasi-disks. Thus $\phi_1$ can be further modified to be a global quasiconformal map without changing it on any $U_x$ such that the modified map $\phi_0$ is isotopic to $\phi_1$ rel $\bigcup U_x\cup\PPP_F$. Let $\psi_0$ be the lift of $\phi_0$. Then they satisfy the conditions.
\qed

\vskip 0.24cm

Let $f$ and $g$ be rational maps. Let $(\phi, \psi)$ be a pair of quasiconformal maps of $\cbar$ such that $\psi$ is isotopic to $\phi$ rel $\PPP_f$ and $\phi\circ f=g\circ\psi$ on $\cbar$. Let
$$
K_b[\phi]=\inf\{K(\wt\phi|_{U}):\, \wt\phi\text{ is isotopic to $\phi$ rel $\PPP_f$ and $U\supset\PPP'_f$ is open}\}.
$$

\REFTHM{qc-conj} With the above assumption, there exists a quasiconformal conjugacy $h$ between $f$ and $g$ in the isotopy class of $\phi$ rel $\PPP_f$ such that $K(h)\le K(\phi)$. Moreover, if $f$ has no Thurston obstructions, then $h$ can be chosen such that $K(h)\le K_b[\phi]$.
\ENDTHM

The existence of the quasiconformal conjugacy $h$ is proved in \cite{Mc2}. The second part of the theorem is obtained in \cite{C}. Refer to the next section for the definition of Thurston obstruction.

\vskip 0.24cm

{\noindent\it Proof of Theorem \ref{unicity}}. Let $f$ and $g$ be geometrically finite rational maps with infinite post-critical sets. Suppose that $(\phi, \psi)$ is a c-equivalence between them. One may choose $\phi$ to be quasiconformal, and hence $K_b[\phi]=1$ by Lemma \ref{c2q}. Thus there exists a holomorphic conjugacy between $f$ and $g$ in the isotopy class of $\phi$, by Theorem \ref{qc-conj}, since $f$ has no Thurston obstructions (refer to Theorem \ref{McMullen} in the next section).
\qed

\section{Thurston obstructions and connecting arcs}

\subsection{Thurston obstructions}

By a {\bf marked branched covering} $(F,\PPP)$ we mean a branched covering $F$ of $\cbar$ with $\deg F\ge 2$ and a closed set $\PPP\subset\cbar$ such that $\PPP_F\subset\PPP$ and $F(\PPP)\subset\PPP$. A marked branched covering $(F,\PPP)$ will be written as $F$ if $\PPP=\PPP_F$.

A simple closed curve on $\cbar\smm\PPP$ is called {\bf essential} if it does not bound a disk in $\cbar\smm\PPP$, or {\bf peripheral}\/ if it encloses a single point of $\PPP$.

A {\bf multicurve}\/ $\Gamma$ of $(F, \PPP)$ is a finite nonempty collection of disjoint simple closed curves in $\cbar\smm\PPP$, each essential and non-peripheral, and no two isotopic rel $\PPP$. It is called {\bf stable}\/ if for any $\gamma\in\Gamma$, every  essential and non-peripheral component of $F^{-1}(\gamma)$ is isotopic rel $\PPP$ to a curve in $\Gamma$.

A multicurve determines a transition matrix $M(\Gamma)=(a_{\beta\gamma})$ by the formula
$$
a_{\beta\gamma}=\sum_{\delta}\frac 1{\deg(F:\,\delta\to\gamma)}
$$
where the sum is taken over all components $\delta$ of $F^{-1}(\gamma)$ which are isotopic to $\beta$ rel $\PPP$. Let $\lambda(\Gamma)\ge 0$ denote the spectral radius of $M(\Gamma)$. A stable multicurve $\Gamma$ on $\cbar\smm\PPP$ is called a {\bf Thurston obstruction} of $(F, \PPP)$ if $\lambda(\Gamma)\ge 1$.

Two marked semi-rational maps $(F, \PPP)$ and $(G, \QQQ)$ are called {\bf c-equivalent} if there exists a fundamental set $\UUU$ of $F$ and a pair of orientation-preserving homeomorphisms $(\phi, \psi)$ of $\cbar$ such that:

(a) $\phi\circ F=G\circ\psi$,

(b) $\phi$ is holomorphic in $\UUU$, and

(c) $\psi=\phi$ in $\UUU\cup\PPP$ and $\psi$ is isotopic to $\phi$ rel $\UUU\cup\PPP$.

This definition coincides with the definition in \cite{BCT} of combinatorial equivalence when $\PPP$ is finite. Refer to \cite{DH2, Mc1} for the definition of hyperbolic orbifold and \cite{BCT} for the following theorem.

\REFTHM{Thurston}{\bf (Marked Thurston Theorem)} Let $(F,\PPP)$ be a marked branched covering of $\cbar$ with hyperbolic orbifold and with $\#\PPP<\infty$. Then $(F, \PPP)$ is c-equivalent to a marked rational map $(f, \QQQ)$ if and only if $(F, \PPP)$ has no Thurston obstructions. Moreover, the marked rational map $(f, \QQQ)$ is unique up to holomorphic c-equivalence.
\ENDTHM

\REFTHM{McMullen}{\bf (McMullen)} Let $(f, \PPP)$ be a marked rational map and let $\Gamma$ be a multicurve on $\cbar\smm\PPP$. Then $\lambda(\Gamma)\le 1$. The equality $\lambda(\Gamma)=1$ holds only in the following cases:

$\bullet$ $f$ is post-critically finite and the signature of the orbifold of $f$ is $(2,2,2,2)$.

$\bullet$ $\PPP_f$ is an infinite set, and $\Gamma$ includes the essential curves in a finite system of annuli permuted by $f$. These annuli lie in Siegel discs or Herman rings for $f$, and each annulus is a connected component of $\cbar\smm\PPP_f$.
\ENDTHM

Refer to \cite{BCT, Mc1} for the above theorem. If $f$ is a geometrically finite rational map, then $f$ has no rotation domains. Hence if $\PPP_f'\ne\emptyset$, then $\la(\Ga)<1$ for any multicurve $\Ga$ on $\cbar\smm\PPP$.

\REFTHM{hyperbolic} {\bf (Sub-hyperbolic version)} Let $(F, \PPP)$ be a marked sub-hyperbolic semi-rational map with $\PPP'_F\ne\emptyset$ and $\#(\PPP\smm\PPP_F)<\infty$. Then $(F, \PPP)$ is c-equivalent to a marked rational map $(f, \QQQ)$ if and only if $(F, \PPP)$ has no Thurston obstructions. Moreover, the marked rational map $(f, \QQQ)$ is unique up to holomorphic c-equivalence.
\ENDTHM

This theorem was proved for the case $\PPP=\PPP_F$ in \cite{CT1,JZ}. One may easily check that the proof in \cite[\S3.3]{CT1} still works in this slightly stronger version. The theorem will be used in \S7 and \S9. We will prove the following theorem, which is stronger than Theorem \ref{existence}. Refer to \S4.2 for the definition of connecting arcs for marked semi-rational maps.

\begin{theorem}\label{existence-1}
Let $(G,\QQQ)$ be a marked semi-rational map with $\PPP'_G\ne\emptyset$ and $\#(\QQQ\smm\PPP_G)<\infty$. Then $(G, \QQQ)$ is c-equivalent to a marked rational map if and only if $(G, \QQQ)$ has neither Thurston obstructions nor connecting arcs.
\end{theorem}

The following lemma (refer to \cite{Mc1}) is useful for checking if there is a Thurston obstruction. A multicurve $\Gamma$ is called {\bf irreducible} if for each pair $(\gamma,\beta)\in\Gamma\times\Gamma$, there is an integer $n\ge 1$ such that $F^{-n}(\beta)$ has a component $\delta$ isotopic to $\gamma$ rel $\sP_F$ and $F^k(\delta)$ is isotopic rel $\sP_F$ to a curve in $\Gamma$ for $1\le k<n$.

\REFLEM{irreducible} For any multicurve $\Gamma$ with $\lambda(\Gamma)>0$, there is an irreducible multicurve $\Gamma_0\subset\Gamma$ such that $\lambda(\Gamma_0)=\lambda(\Gamma)$.
\ENDLEM

\subsection{Connecting arcs}

Let $(G, \QQQ)$ be a marked semi-rational map with $\#(\QQQ\smm\PPP_G)<\infty$ and with parabolic cycles in $\PPP'_G$. An open arc $\beta\subset\cbar\smm\QQQ$ which joins two points $z_0, z_1\in\PPP'_G$ is a connecting arc if:

$\bullet$ either $z_0\neq z_1$, or $z_0=z_1$ and both components of $\cbar\smm\overline{\beta}$ contain points of $\QQQ$,

$\bullet$ $\beta$ is disjoint from a fundamental set of $G$, and

$\bullet$ $\beta$ is isotopic rel $\QQQ$ to a component of $G^{-p}(\beta)$ for some integer $p>0$.

\vskip 0.24cm

{\noindent\bf Example}. Let $G$ be the formal mating of the quadratic polynomial $P(z)=z^2+\frac 14$ with itself (refer to \cite{T1} for a detailed definition of mating). It will be a semi-rational map with two points in $\PPP'_G$ if we preserve its complex structure near $\PPP'_G$. Consider two external rays, each with angle zero. They form an invariant arc and hence the arc is a connecting arc of $G$. It is easy to check that $G$ has no Thurston obstructions since it is combinatorially equivalent to a Blaschke product.

\vskip 0.24cm

A connecting arc is invariant under c-equivalence by Lemma \ref{equivalent}. The following lemma gives a stronger version of the definition.

\REFLEM{connecting0} Let $\beta\subset\cbar\smm\QQQ$ be a connecting arc. Then there exist a connecting arc $\alpha$ isotopic to $\beta$ rel $\QQQ$ and an integer $p\ge 1$ such that $G^{-p}(\alpha)$ has a component $\wt\alpha$ isotopic to $\alpha$ rel $\QQQ$ and $\wt\alpha$ coincides with $\alpha$ in a neighborhood of its endpoints.
\ENDLEM

\beginp
Denote by $\beta: (0,1)\to\beta$ a parametrization, which is a homeomorphism. Let $z_0, z_1\in\PPP'_G$ be the two endpoints of $\beta$, with $\beta(t)\to z_0$ as $t\to 0$. Then there exists an integer $p\ge 1$ such that $G^{-p}(\beta)$ has a component $\wt\beta$ isotopic to $\beta$ rel $\QQQ$ and $G^{p}(z_i)=z_i$ for $i=0,1$. Denote by $\wt\beta: (0,1)\to\wt\beta$ the parametrization such that $G^p\circ\wt\beta(t)=\beta(t)$.

Since $\beta$ is disjoint from a fundamental set, both $z_0$ and $z_1$ are parabolic periodic points of $G$. By Proposition \ref{regular} condition (3), there exist constants $0<t_0<t_1<1$ and regular repelling flowers $\VVV'_0, \VVV'_1\subset\cbar\smm\QQQ$ at the parabolic fixed points $(G^{p}, z_0)$ and $(G^{p}, z_1)$, respectively, such that $\beta(0,t_0)\subset\VVV'_0$ and $\beta(t_1,1)\subset\VVV'_1$. Obviously, $\VVV'_0$ is disjoint from $\VVV'_1$ if $z_0\neq z_1$. If $z_0=z_1$ we may require that $\VVV'_0=\VVV'_1$.

By the definition of semi-rational maps, there are critical orbits converging to the point $z_i$ between any two adjacent repelling petals of $\VVV'_i$. Thus there exist constants $0<s_0<s_1<1$ such that either

(1) both $\wt\beta(0, s_0]$ and $\beta(0,t_0]$ are contained in the same petal of $\VVV'_0$, and both $\wt\beta[s_1, 1)$ and $\beta[t_1, 1)$ are contained in the same petal of $\VVV'_1$, or

(2) $z_0=z_1$, both $\wt\beta(0, s_0]$ and $\beta[t_1,1)$ are contained in the same petal of $\VVV'_0=\VVV'_1$, and both $\wt\beta[s_1, 1)$ and $\beta(0, t_1]$ are contained in the same petal of $\VVV'_0$.

If condition (2) holds, we replace $\wt\beta$ by the component of $G^{-2p}(\beta)$ that is isotopic to $\beta$ rel $\PPP_G$. Then condition (1) holds. Thus we may assume that condition (1) holds. Consequently, the rotation numbers of $G^{p}$ at $z_0$ and $z_1$ are both equal to $1$.

Denote by $g_i$ ($i=0,1$) the inverse map of $G^{p}$ restricted to the repelling flower $\VVV'_i$.  Then there exists an integer $m\ge 1$ such that $\beta(t_0, t_1)$ is disjoint from $\overline{g_i^m(\VVV'_i)}$. Thus there exists an arc $\gamma_i\subset\VVV'_i\smm \overline{g_i^m(\VVV'_i)}$ connecting $\beta(t_i)$ with a point $w_i\in\partial g_i^m(\VVV'_i)$ such that $\gamma_i$ is disjoint from $\beta(t_0, t_1)$. Obviously, $\gamma_0$ is disjoint from $\gamma_1$ if $z_0\neq z_1$. Otherwise we may assume that they are disjoint by a suitable choice of $w_i$.

Since the rotation number of $G^{p}$ at $z_0$ and at $z_1$ is equal to $1$, there exists an arc $\alpha_i\subset g_i^m(\VVV'_i)$ connecting $z_i$ with $w_i$ such that $g_i(\alpha_i)\subset\alpha_i$ and $\alpha_0$ is disjoint from $\alpha_1$. Set
$$
\alpha=\alpha_0\cup\overline{\gamma_0}\cup\beta(t_0,t_1)\cup\overline{\gamma_1}\cup\alpha_1.
$$
Since $\VVV'_i$ is disjoint from $\QQQ$, $\alpha$ is isotopic to $\beta$ rel $\QQQ$. Moreover, $G^{-p}(\alpha)$ has a component $\wt\alpha$ isotopic to $\alpha$ rel $\QQQ$, and both $g_0(\alpha_0)$ and $g_1(\alpha_1)$ are contained in $\wt\alpha$. Now the lemma follows from the fact that $g_i(\alpha_i)\subset\alpha_i$.
\qed

\REFTHM{connecting}{\bf (No connecting arcs)} Any marked geometrically finite rational map has no connecting arcs.
\ENDTHM

\beginp Let $(g, \QQQ)$ be a marked geometrically finite rational map. Assume that $\beta_0\subset\cbar\smm\QQQ$ is a connecting arc, i.e., $\beta_0$ joins two parabolic periodic points $z_0, z_1\in\PPP'_g$, $\beta_0$ is disjoint from an attracting flower of $z_0$ and one of $z_1$, and $\beta_0$ is isotopic rel $\QQQ$ to a component of $g^{-p}(\beta_0)$ for some integer $p\ge 1$.

There exist repelling flowers $\VVV'_0$, $\VVV'_1$ at $z_0$ and $z_1$, respectively, such that $\beta_0$ is cut into three arcs $\beta_{0,0}$, $\beta_{0,1}$ and $\beta_{0,2}$ where $\beta_{0,i}\subset\VVV'_i$ for $i=0,1$ and the closure of $\beta_{0,2}$ is disjoint from $\QQQ$. Let $\beta_n$ be the component of $g^{-np}(\beta)$ isotopic rel $\QQQ$ to $\beta_0$ for $n\ge 1$. Then $\beta_n$ is also cut into three arcs $\beta_{n,0}$, $\beta_{n,1}$ and $\beta_{n,2}$ such that $g^{np}(\beta_{n,j})=\beta_{0,j}$ for $j=0,1,2$. Thus $\beta_{n,i}\subset\VVV'_i$ and hence $\text{diam}_s \beta_{n,i} \to 0$ as $n\to\infty$ for $i=0, 1$. By Lemma \ref{shrinking}, $\text{diam}_s \beta_{n,2}\to 0$ as $m\to\infty$. So $\text{diam}_s \beta_n\to 0$ as $n\to\infty$. This shows that $z_0=z_1$ and one component of $\cbar\smm(\beta\cup\{z_0\})$ is disjoint from $\QQQ$. This is a contradiction.
\qed

\subsection{Thurston's algorithm}

This part is not needed until \S9; one may skip it on a first reading. Let $(F, \PPP)$ be a marked semi-rational map with $\PPP'_F\neq\emptyset$ and $\#(\PPP\smm\PPP_F)<\infty$.

\vskip 0.24cm

{\bf Thurston sequences}.
Let $x_i$ ($i=1,2,3$) be three distinct points in $\PPP$. Then there exists a unique homeomorphism $\theta_1$ of $\cbar$ normalized by $\theta_1(x_i)=x_i$ such that $f_1:=F\circ\theta_1^{-1}$ is a rational map by the Uniformization Theorem. There is also a unique normalized homeomorphism $\theta_2$ of $\cbar$ such that $f_2:=\theta_1\circ F\circ\theta_2^{-1}$ is a rational map. Continuing this process inductively, we produce a sequence of normalized homeomorphisms $\{\theta_n\}$ of $\cbar$ and a sequence of rational maps $\{f_n\}$ such that $f_n\circ\theta_{n}=\theta_{n-1}\circ F$. We call $\{f_n\}$ a {\bf Thurston sequence}\/ of $(F, \PPP)$.

\vskip 0.24cm

{\bf Lift of the c-equivalence}. Denote by $\PPP^s_F$ the set of super-attracting periodic points of $F$. Then $\PPP^s_F\subset\PPP_F$ and $F(\PPP^s_F)=\PPP^s_F$. Assume that $F$ is holomorphic in a neighborhood of $\PPP^s_F$ and c-equivalent to a marked rational map. Then there exists a normalized c-equivalence $(\phi_0,\phi_1)$ from $(F, \PPP)$ to a marked rational map $(f, \QQQ)$ on a fundamental set $\UUU$ of $F$ such that $\PPP^s_F\subset\UUU$. Refer to \cite{Sh} for the construction of $\phi_0$ near $\PPP^s_F\smm\PPP'_F$.

Let $\phi_2$ be the lift of $\phi_1$. Then $\phi_2$ is isotopic to $\phi_1$ rel $F^{-1}(\PPP\cup\UUU)$. Inductively, we obtain a sequence of homeomorphisms $\{\phi_n\}$ of $\cbar$ such that $\phi_{n+1}$ is isotopic to $\phi_n$ rel $F^{-n}(\PPP\cup\UUU)$ and $f\circ\phi_{n+1}=\phi_n\circ F$. See the diagram below.
$$
\xymatrix{
\ar@{.>}[d] & \ar@{.>}[d]  & \ar@{.>}[d] \\
{\rule[-1ex]{0ex}{4ex}{}} \cbar\ar[d]_{f} & {\rule[-1ex]{0ex}{4ex}{}} \cbar\ar[d]_{F} \ar[l]_{\phi_2}
\ar[r]^{\theta_2} & {\rule[-1ex]{0ex}{4ex}{}} \cbar\  \ar[d]^{f_2}  \\
{\rule[-1ex]{0ex}{4ex}{}} \cbar\ar[d]_{f} & {\rule[-1ex]{0ex}{4ex}{}} \cbar\ar[d]_{F} \ar[l]_{\phi_1}
\ar[r]^{\theta_1} & {\rule[-1ex]{0ex}{4ex}{}} \cbar\  \ar[d]^{f_1}  \\
{\rule[-1ex]{0ex}{4ex}{}} \cbar & {\rule[-1ex]{0ex}{4ex}{}} \cbar \ar[l]_{\phi_0}
\ar[r]^{\text{id}} & {\rule[-1ex]{0ex}{4ex}{}} \cbar}
$$

Let $\zeta_0=\phi_0^{-1}$. Then $\zeta_0$ is holomorphic in $\phi_0(\UUU)$. Let $\zeta_n=\theta_n\circ\phi_n^{-1}$ for $n\ge 1$. Then $f_{n}\circ\zeta_{n}=\zeta_{n-1}\circ f$. Consequently, $\zeta_{n}$ is holomorphic in $f^{-n}(\phi_0(\UUU))$.

\REFTHM{algorithm}
The sequence $\{f_n\}$ converges uniformly to the rational map $f$ and $\{\zeta_n\}$ converges uniformly to the identity as $n\to\infty$.
\ENDTHM

To prove this theorem, we need the lemma below, which is more general than needed here but which will be used in \S9.

\vskip 0.24cm

{\bf Combinatorial quotient maps}. Let $(f, \QQQ)$ be a marked geometrically finite rational map with $\PPP'_f\neq\emptyset$ and $\#(\QQQ\smm\PPP_f)<\infty$. Let $\UUU$ be a fundamental set of $f$ with $\PPP^s_f\subset\UUU$. Let $\{h_n\}$ ($n\ge 0$) be a sequence of quotient maps of $\cbar$ such that $f\circ h_{n+1}=h_n\circ f$ and $h_n$ is isotopic to the identity rel $f^{-n}(\overline{\UUU}\cup\QQQ)$.

\REFLEM{algorithm1}
The sequence $\{h_n\}$ converges uniformly to the identity as $n\to\infty$.
\ENDLEM

\beginp Let $h_{t,0}$, $t\in I=[0,1]$, be an isotopy of quotient maps connecting $h_{0, 0}=\text{\rm id}$ with $h_{1, 0}=h_0$ such that $h_{t, 0}$ is a quotient map of $\cbar$ for all $t\in I$ and $h_{t,0}^{-1}(w)=w$ for $w\in\overline{\UUU}\cup\QQQ$ and $t\in I$. Let $h_{t,1}$ be the lift of $h_{t,0}$. Then $h_{t,0}\circ f=f\circ h_{t,1}$, $h_{0,1}=\text{\rm id}$, $h_{1,1}=h_1$ and $h_{t,1}^{-1}(w)=w$ for $w\in f^{-1}(\overline{\UUU}\cup\QQQ)$ and $t\in I$. Inductively, let $h_{t,n}$ be the lift of $h_{t,n-1}$. Then
$$
\begin{cases}
h_{0,n}=\text{\rm id}, \quad h_{1,n}=h_n, \\
h_{t,n}^{-1}(w)=w\text{ for $w\in f^{-n}(\overline{\UUU}\cup\QQQ)$ and $t\in I$, and} \\
h_{t,n-1}\circ f=f\circ h_{t,n}.
\end{cases}
$$

Let $\beta_n(z)=\{h_{t,n}(z):\, t\in I\}$. Then $f: \beta_{n+1}(z)\to\beta_{n}(f(z))$ is injective, and
$$
\text{\rm dist}_s(h_n(z), z)\le \text{\rm diam}_s\beta_n(z).
$$
We want to prove that $\text{\rm diam}_s\beta_n(z)\to 0$ as $n\to\infty$ uniformly for $z\in\cbar$.

For any disk $U\subset\cbar$ and any integer $n\ge 0$, we denote by $C_n(U)$ the maximum of the diameters of the components of $f^{-n}(U)$. Then $C_n(U)\to 0$ as $n\to\infty$ if the closure of $U$ is disjoint from $\PPP_f$, by Lemma \ref{shrinking}.

Let $w\in\cbar\smm(\UUU\cup\QQQ)$ be a point. Then $\beta_0(w)$ is disjoint from $\QQQ$. If it is simple, then there exists a disk $U\supset\beta_0(w)$ such that $\overline{U}\cap\QQQ=\emptyset$. Thus there exists an open set $D_w\ni w$ such that $\beta_0(z)\subset U$ for all points $z\in D_w$. Therefore $\text{\rm diam}_s\beta_n(z)\le C_n(U)\to 0$ as $n\to\infty$ for $z\in f^{-n}(D_w)$.

In general, let $d:=\text{\rm dist}_s(\beta_0(w), P)>0$. We can cut the path $\beta_0(w)$ into $k$ sub-paths such that each of them has diameter
less than $d$. Thus each sub-path is contained in a disk $U_i$ ($1\le i\le k$) whose closure is disjoint from $\QQQ$. Denote their union by $U$. Similarly, there exists an open set $D_w\ni w$ such that $\beta_0(z)\subset U$ for all points $z\in D_w$. Therefore $\text{\rm diam}_s\beta_n(z)\le \sum_{i=1}^kC_n(U_i)\to 0$ as $n\to\infty$ for $z\in f^{-n}(D_w)$.

Assume that $w\in\QQQ\smm\UUU$ is not periodic. Let $U\ni w$ be a disk whose closure is disjoint from $\QQQ\smm\{w\}$. Then $C_n(U)\to 0$ as $n\to\infty$. Since $\beta_0(w)=w$, there exists an open set $D_w\ni w$ such that $\beta_0(z)\subset U$ for $z\in D_w$. Thus $\text{\rm diam}_s\beta_n(z)\le C_n(U)\to 0$ as $n\to\infty$ for $z\in f^{-n}(D_w)$.

Assume $w\in\QQQ\smm\UUU$ is a repelling periodic point with period $p\ge 1$. Let $U\ni w$ be a disk whose closure is disjoint from $\QQQ\smm\{w\}$ such that the component of $f^{-p}(U)$ containing the point $w$ is compactly contained in $U$. Then $C_n(U)\to 0$ as $n\to\infty$. As above, there exists an open set $D_w\ni w$ such that $\text{\rm diam}_s\beta_n(z)\le c_n(U)\to 0$ as $n\to\infty$ for $z\in f^{-n}(D_w)$.

Suppose now $w\in\QQQ\smm\UUU$ is a parabolic periodic point with period $p\ge 1$. Let $\VVV_w$ be a repelling flower of $(F^p,w)$ such that its closure is disjoint from $\QQQ\smm\{w\}$ and $\VVV_w\cup\UUU$ is a neighborhood of $w$. Then $C_n(\VVV_w)\to 0$ as $n\to\infty$. Since $\beta_0(w)=w$, there exists an open set $D_w\ni w$ such that $\beta_0(z)\subset\VVV_w\cup\UUU$ for $z\in D_w$. In particular, $\beta_0(z)\subset\VVV_w$ if $z\in D_w\smm\UUU$. Therefore $\text{\rm diam}_s\beta_n(z)\le C_n(\VVV_w)\to 0$ as $n\to\infty$ for $z\in f^{-n}(D_w\smm\UUU)$.

The union of these open sets $D_w$ forms an open cover of $\cbar\smm\UUU$. Hence there is a finite cover. Thus $\text{\rm diam}_s\beta_n(z)\to 0$ as $n\to\infty$ uniformly for $z\in\cbar\smm f^{-n}(\UUU)$. On the other hand, $\beta_n(z)=z$ for $z\in f^{-n}(\UUU)$. Therefore $\text{\rm diam}_s\beta_n(z)\to 0$ as $n\to\infty$ uniformly for $z\in\cbar$. This completes the proof.
\qed

\vskip 0.24cm

Let $f$ be a geometrically finite rational map with $\PPP'_f\neq\emptyset$. Let $\UUU$ be a fundamental set of $f$ with $\PPP^s_f\subset\UUU$. Assume that $j_0$ is a quasiconformal map of $\cbar$ which is holomorphic in $\UUU$ and normalized by fixing three points in $\PPP_f$. Then there exists a unique normalized quasiconformal map $j_1$ of $\cbar$ such that $f_1:=j_0\circ f\circ j_1^{-1}$ is a rational map. Inductively, there exists a sequence of normalized quasiconformal maps $\{j_n\}$ ($n\ge 1$) such that
$$
f_n:=j_{n-1}\circ f\circ j_{n}^{-1}
$$
is a rational map.

\begin{lemma}\label{algorithm2}
The sequence $\{j_n\}$ converges uniformly to the identity and $\{f_n\}$ converges uniformly to $f$ as $n\to\infty$.
\end{lemma}

\beginp
Since $f_n=j_{n-1}\circ f\circ j_{n}^{-1}$ is a rational map, the sequence $\{j_n\}$ is uniformly quasiconformal. So it has a subsequence which converges uniformly to a quasiconformal map $j_{\infty}$. Since $j_n$ is holomorphic in $f^{-n}(\UUU)$, the map $j_{\infty}$ is holomorphic in $\FFF_f$ and hence is holomorphic on the whole sphere since $\JJJ_f$ has zero Lebesgue measure \cite{U}. Thus it is the identity since it fixes three points in $\PPP_f$. It follows that the whole sequence $\{j_n\}$ uniformly converges to the identity as $n\to\infty$. Consequently, $\{f_n\}$ converges uniformly to $f$.
\qed

\vskip 0.24cm

{\noindent\it Proof of Theorem \ref{algorithm}}. By Lemma \ref{c2q}, there exist a fundamental set $\UUU_1\subset\UUU$ of $F$ with $\PPP^s_F\subset\UUU_1$ and a quasiconformal map $\psi_0$ of $\cbar$ normalized by fixing $x_i$ ($i=1,2,3$), such that $\psi_0$ is isotopic to $\phi_0$ rel $\UUU_1\cup P$. Let $\psi_{n+1}$ be the lift of $\psi_n$. Then $\phi_0\circ\psi_0^{-1}$ is isotopic to the identity rel $\phi_0(\UUU)\cup P$ and
$$
(\phi_n\circ\psi_n^{-1})\circ f=f\circ(\phi_{n+1}\circ\psi^{-1}_{n+1}).
$$
Thus $\{\phi_n\circ\psi_n^{-1}\}$ converges uniformly to the identity as $n\to\infty$, by Lemma \ref{algorithm1}.

Set $j_n=\zeta_n\circ\phi_n\circ\psi_n^{-1}$ for $n\ge 0$. Then $j_0=\psi_0^{-1}$ is quasiconformal in $\cbar$ and holomorphic in $\psi_0(\UUU_1)$. Moreover,
$$
f_{n+1}\circ j_{n+1}=j_{n}\circ f.
$$
By Lemma \ref{algorithm2}, $\{j_n\}$ converges uniformly to the identity and $\{f_n\}$ converges uniformly to $f$ as $n\to\infty$. Thus $\{\zeta_n\}$ also converges uniformly to the identity as $n\to\infty$.
\qed

\section{Basic properties of pinching}

\subsection{Definition of pinching}

{\bf Pinching model}. For any $r>1$ and $t\ge 0$, define a quasiconformal map
$$
w_{t,r}=w_t:\, \AA(r)\to\AA(r^{1+t})
$$
by $\arg w_t(z)=\arg z$ and $\log|w_t(z)|=\varrho(\log|z|)$, where
$$
\varrho: (-\log r, \log r)\to (-(1+t)\log r,\ (1+t)\log r)
$$
is defined by $\varrho(-x)=-\varrho(x)$ and
$$
\varrho(x)=\begin{cases}
e^{2t}x & \text{if } 0\le x\le \dfrac{\log r}{2e^{2t}}, \\
\dfrac{1}{2}\left(\log\frac{2x}{\log r}+1+2t\right)\log r &  \text{if } \dfrac{\log r}{2e^{2t}}<x<\frac{1}{2}\log r, \\
x+t\log r & \text{if } \frac{1}{2}\log r\le x<\log r.
\end{cases}
$$
The family $\{w_t\}$ $(t\ge 0)$ is called a {\bf pinching model}.

Let $\nu_t$ be the Beltrami differential of $w_t$. Let $r(t)=r^{1/(2e^{2t})}$ and $r'=r(0)=\sqrt{r}$. The following proposition is easy to check.

\REFPROP{model} The pinching model $w_t(z)$ satisfies the following properties:

(1) $w_t(z)$ is conformal on $\AA(r)\smm\overline{\AA(r')}$.

(2) $\nu_t(z)=\nu_{t_0}(z)$  on $\AA(r)\smm\overline{\AA(r(t_0))}$ for $t\ge t_0\ge 0$.

(3) For any $t_0\ge 0$, let $E_{t_0}$ be a component of $\AA(r)\smm\overline{\AA(r(t_0))}$; then
$$
\md w_t(E_{t_0})=\dfrac{2t_0+1}{4}\md \AA(r)
$$
for $t\ge t_0$ and hence $\md w_t(E_{t_0})\to\infty$ as $t\ge t_0\to\infty$.

(4) The map $w_{t}(rz)/r^{1+t}$ restricted to $\AA(1/r,1)$ converges uniformly to a homeomorphism $w:\, \AA(1/r,1)\to\DD^*$ as $t\to\infty$, where $w$ is defined as:
$$
\begin{cases}
w(z)=z & \text{ if }{\sqrt{1/r}}\le |z|<1, \\
\arg w(z)=\arg z & \text{ if }1/r<|z|<{\sqrt{1/r}}, \\
\log|w(z)|=-\dfrac{1}{2}\left(1+\log\dfrac{\log r}{2\log(r|z|)}\right)\log r
& \text{ if }1/r<|z|<{\sqrt{1/r}}.
\end{cases}
$$
\ENDPROP

{\bf Remark}. The pinching model $w_{t}$ has usually been defined as $w_t(z)=|z|^{t}z$ (refer to \cite{Mas}). We choose the above technical definition  to ensure the convergence of the quasiconformal conjugacy path $\{\phi_t\}$ defined below.

\vskip 0.24cm

{\bf Multi-annuli}. Let $f$ be a rational map with non-empty quotient space $\sR_f$. A {\bf multi-annulus} $\sA\subset\sR_f$ is a finite disjoint union of annuli whose boundaries are pairwise disjoint simple closed curves in $\sR_f$ such that each component of $\pi_f^{-1}(e(\sA))$ is an arc, where $e(\sA)$ denotes the union of the equators of the annuli in $\sA$.

A multi-annulus $\sA\subset\sR_f$ is called {\bf non-separating}\ if for any choice of finitely many components of $\pi_f^{-1}(\sA)$,  the union $T$ of their closures does not separate the Julia set $\JJJ_f$, i.e. $\cbar\smm T$ has exactly one component intersecting $\JJJ_f$ (refer to \cite{T2} Example $3'$ for a non-separating case, and Examples $6'$ and 8 for two separating cases).

Let $\sA\subset\sR_f$ be a multi-annulus. Then each component of $\pi_f^{-1}(e(\sA))$ is an eventually periodic arc. The multi-annulus $\sA$ is called {\bf starlike}\ if for each component $\beta$ of $\pi_f^{-1}(e(\sA))$, $r(\beta)$ is eventually repelling, $a(\beta)$ is eventually attracting, and $a(\beta_1)\neq a(\beta_2)$ for any two distinct components $\beta_1, \beta_2$ of $\pi_f^{-1}(e(\sA))$. Obviously, a starlike multi-annulus is non-separating.

\vskip 0.24cm

{\bf Pinching paths}. Let $f$ be a rational map with $\sR_f\neq\emptyset$. Let $\sA=\bigcup A_i\subset\sR_f$ be a non-separating multi-annulus. Let $\chi_i$ be a conformal map from $A_i$ onto $\AA(r_i)$ and let $\mu_{i,t}$ be the Beltrami differential of $w_{t,r_i}\circ\chi_i$, where $w_{t,r_i}$ is the pinching model on $\AA(r_i)$. Set $\mu_t=\mu_{i,t}$ on each $A_i\subset\sA$. Let $\tilde\mu_t$ be the pullback of $\mu_t$, i.e.,
$$
\tilde\mu_t(z)=\begin{cases}
\mu_t(\pi(z))\dfrac{\overline{\pi_f'(z)}}{\pi_f'(z)}\qquad
& \text{for } z\in\pi_f^{-1}(\sA), \\
0 & \text{otherwise}. \end{cases}
$$
Then there exists a quasiconformal map $\phi_t:\,\cbar\to\cbar$ whose Beltrami differential is $\tilde\mu_t$. Set $f_t=\phi_t\circ f\circ\phi_t^{-1}$. Then $f_t$ is a rational map. The quasiconformal map $\phi_t$ has a natural projection
$$
\Phi_t: \sR_f\to\sR_{f_t}.
$$

We call the path $f_t=\phi_t\circ f\circ\phi_t^{-1}$ $(t\ge 0)$ the {\bf pinching path} starting from $f$ supported on $\sA$, or a {\bf simple pinching path} if $\sA$ is starlike.

\vskip 0.24cm

Note that the family $\{f_t\}$ is defined only up to holomorphic conjugation and hence represents a family in $\fm_d$, the complex orbifold of holomorphic conjugate classes of rational maps with degree $d=\deg f$. It is convenient to consider $\{f_t\}$ as a family in the space of rational maps when we study its convergence in $\fm_d$. For this purpose, we need to make a normalization for the map $\phi_t$.

One favorite choice of a normalization of $\phi_t$ is fixing three points in $\PPP_f$. There always exists a component $U_0$ of $\cbar\smm\overline{\pi_f^{-1}(\sA)}$ such that both $U_0\cap\PPP_f$ and $\partial U_0\cap\JJJ_f$ are infinite sets. Throughout this paper we always make a normalization for the map $\phi_t$ by fixing three distinct points in $U_0\cap\PPP_f$. Then both $\{\phi_t\}$ and $\{f_t\}$ are continuous families.

Such a choice of $U_0$ is necessary. As we will see later, some components of $\cbar\smm\overline{\pi_f^{-1}(\sA)}$ may touch $\JJJ_f$ at only finitely many points, and images of the components under $\phi_t$ will shrink to single points as $t\to\infty$.

\vskip 0.24cm

Let $A_i(t)=\chi_{i}^{-1}(\AA(r_i(t)))$ for $t\ge 0$ and $A'_i=\chi_{i}^{-1}(\AA(r_i'))=A_i(0)$. Then $A_i(t_1)\subset A_i(t_2)$ if $t_1\ge t_2$. Denote by $\sA'$, $\sA(t)$ the union of $A'_i$ and $A_i(t)$ for all components $A_i$ of $\sA$, respectively. The following proposition is a direct consequence of Proposition \ref{model}.

\REFPROP{A(t)} Let $A_i$ be a component of $\sA$. Then the following conditions hold:

(1) $\mu_t(z)=0$ on $A_i\smm A'_i$.

(2) $\mu_t(z)=\mu_{t_0}(z)$ on $A_i\smm A_i(t_0)$ for $t\ge t_0\ge 0$.

(3) Let $t_0\ge 0$ and let $E$ be a component of $A_i\smm\overline{A_i(t_0)}$. Then
$$
\md \Phi_t(E) =\dfrac{2t_0+1}{4}\md A_i
$$
for $t\ge t_0$ and hence $\md\Phi_t(E)$ tends to infinity as $t_0\to\infty$.
\ENDPROP

\subsection{Bands and skeletons}

We now need to analyze in more detail the structure of the lifts of a multi-annulus. Let $f$ be a geometrically finite rational map and let $\sA\subset\sR_f$ be a non-separating  multi-annulus. A component $B$ of $\pi_f^{-1}(\sA)$ is called a {\bf band}. It is of {\bf level $0$} if it is periodic, or {\bf level $n$} with $n\ge 1$ if $f^n(B)$ is periodic but $f^{n-1}(B)$ is not periodic.

Any band $B$ is bounded by two eventually periodic arcs with a common attracting end and a common repelling end. We denote them by $a(B)$ and $r(B)$, respectively. A band $B$ is periodic if and only if both of its endpoints are periodic. Consequently, bands of different levels have disjoint closures.

The {\bf core arc} of a band $B$ is the lift of the equator of $\pi_f(B)$ to $B$. Its two endpoints are exactly the endpoints of $B$. A {\bf skeleton of level $n$}\/ is a component of the union of core arcs of all the level $n$ bands together with their endpoints.

By the non-separating condition, each skeleton has exactly one complementary component intersecting $\JJJ_f$. The {\bf fill-in}\/ of a skeleton $S$, denote by $\wh S$, is the union of $S$ together with all its complementary components disjoint from $\JJJ_f$. Thus $\wh S$ is a full continuum, i.e. $\cbar\smm\wh S$ is connected. A {\bf filled-in skeleton} is the fill-in of a skeleton. Each component of $f^{-1}(\wh S)$ is also a filled-in skeleton.

There is an integer $n_0\ge 1$ such that for any skeleton $S$ of level $n\ge n_0$, $\wh S\cap\PPP_f=\emptyset$. Thus one need only check finitely many levels of skeletons to see whether $\sA$ is non-separating. The following proposition is easy to verify.

\begin{proposition}\label{infty} Let $\wh S$ be a periodic filled-in skeleton and $x\in\wh S$ be a periodic point. Let $k\ge 1$ be the number of components of $\wh S\smm\{x\}$. Let $D_x\ni x$ be a sufficiently small disk, such that $D_x\smm\wh S$ has $k$ components $U_i$ whose closures contain the point $x$. Then $U_i$ contains infinitely many points of $\PPP_f$ if there exists a periodic band $B$ such that $a(B)=x$ and $U_i\cap B\neq\emptyset$.
\end{proposition}

In general, the map $f$ need not be injective on a filled-in skeleton. However, $f$ is injective on periodic filled-in skeletons.

\REFPROP{injective}
Let $\wh S$ be a periodic filled-in skeleton. Then $f$ is injective in a neighborhood of $\wh S$.
\ENDPROP

\beginp We need only consider the case that $\wh S\neq S$. Let $U$ be a component of $\wh S\smm S$. Then $U$ is a component of $\cbar\smm S$ by the definition of $\wh S$. We claim that $U$ is bounded by either two periodic arcs $\alpha, \beta\subset S$ with $a(\alpha)=a(\beta)$ and $r(\alpha)=r(\beta)$, or one periodic arc $\gamma\subset S$ with $a(\gamma)=r(\gamma)$. Otherwise, if $U$ is bounded by distinct periodic arcs $\gamma_1,\cdots, \gamma_n$ with $n\ge 3$, let $p\ge 1$ be an integer such that $f^p$ fixes all these arcs; then $f^p(U)=U\subset\FFF_f$ by the non-separating condition. Thus all the $\gamma_i$ have the same attracting ends. This is a contradiction. Therefore $f$ is injective on $U$. This implies that $f$ is injective in a neighborhood of $\wh S$.
\qed

\vskip 0.24cm

By a {\bf band-tree} of level $n$ we mean a connected component of the closure of the union of all the level $n$ bands and filled-in skeletons.

\begin{figure}[htbp]
\begin{center}
\includegraphics[scale=0.5]{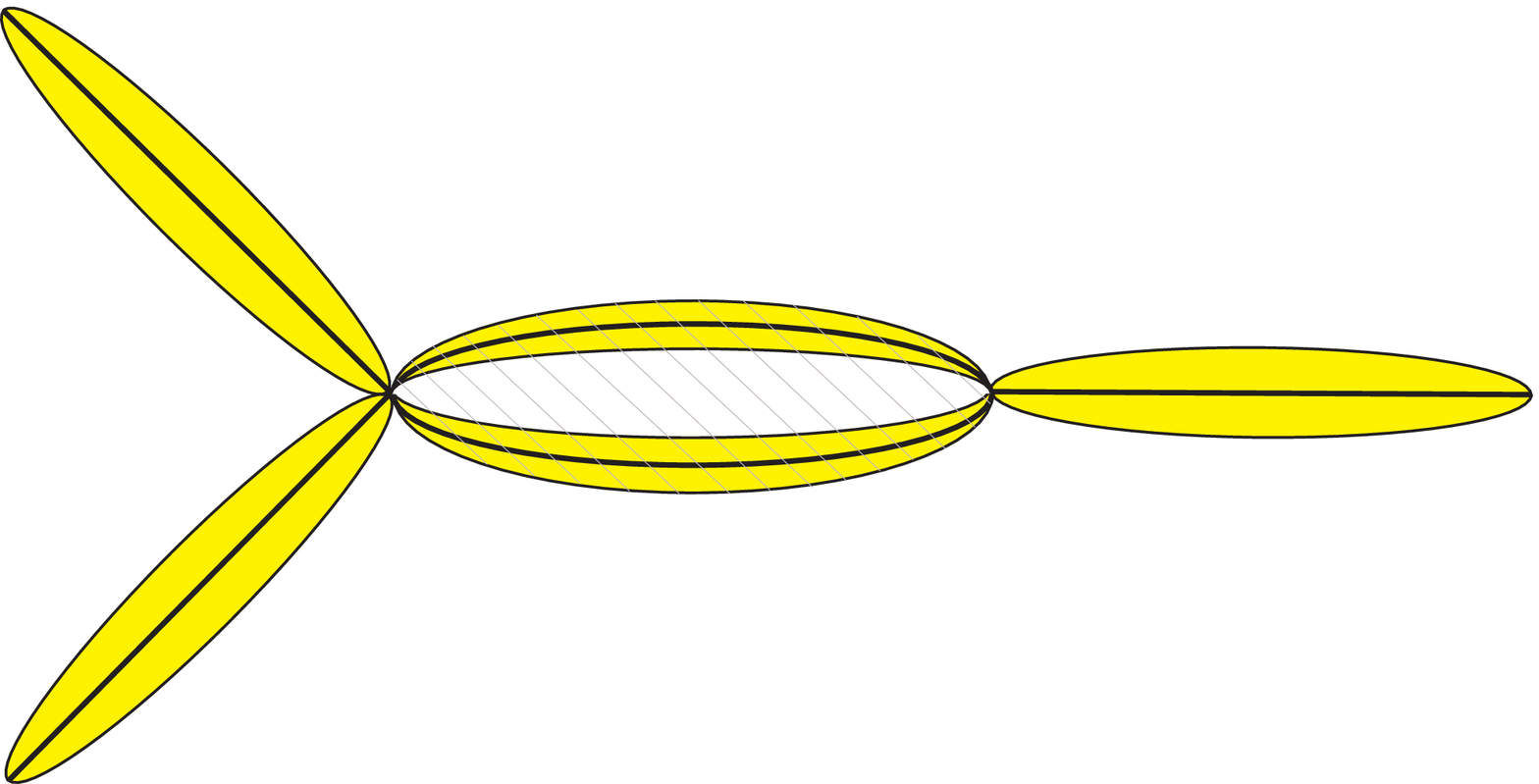}
\end{center}
\begin{center}{\sf Figure 4. A band-tree.}
\end{center}
\end{figure}

Pick pairwise disjoint disks $N(T)\supset T$ for all periodic band-trees $T$ such that $N(T)\smm T$ is disjoint from the critical values of $f$ and $\partial N(T)\cap\PPP_f=\emptyset$. Then each component of $f^{-1}(N(T))$ is also a disk containing exactly one component of $f^{-1}(T)$.

For each level $1$ band-tree $T_1$ with $f(T_1)=T$, denote by $N(T_1)$ the component of $f^{-1}(N(T))$ that contains $T_1$. Then $N(T_1)$ is disjoint from all periodic band-trees. Since $T_1$ is disjoint from $\PPP'_f$, taking $N(T)$ small enough, we may assume that $N(T_1)\smm T_1$ is disjoint from $\PPP_f$.

For each $n>1$ and each component $T_n$ of $f^{-n+1}(T_1)$, let $N(T_n)$ be the component of $f^{-n+1}(N(T_1))$ that contains $T_n$; then $N(T_n)$ is a disk disjoint from band-trees of level $k\le n$ except for $T_n$, and $N(T_n)\smm T_n$ is disjoint from $\PPP_f$.

Note that there is an integer $n_0\ge 1$ such that each band-tree of level $n\ge n_0$ is disjoint from $\PPP_f$. By Lemma \ref{shrinking}, $\text{diam}_s N(T_n)\to 0$ as $n\to\infty$ uniformly for all $n$-level band-trees $T_n$. We have proved:

\REFPROP{Koebe}
There exist a constant $M>0$ and a disk $N(T)\supset T$ for each band-tree $T$ satisfying the following conditions:

(a) $\partial N(T)\cap\PPP_f=\emptyset$.

(b) $f: N(T)\smm T\to N(f(T))\smm f(T)$ is a covering if $T$ is not periodic.

(c) $\md(N(T)\smm T)>M$.

(d) $N(T_n)$ is disjoint from all band-trees of level $k\le n$ except for $T_n$.

(e) $\text{diam}_s N(T_n) \to 0$ as $n\to\infty$ uniformly for all $n$-level band-trees $T_n$.
\ENDPROP

We will call such a choice of the neighborhood $N(T)$ a {\bf Koebe space} of the band-tree $T$, and the collection $\{N(T)\}$ a {\bf Koebe space system}. For a band $B$ or a skeleton $S$, we also use $N(B)$ or $N(S)$ to denote the Koebe space $N(T)$ if $B\subset T$ or $S\subset T$.

For any point $z_0\in\cbar$, its {\bf $\omega$-limit set} $\omega(z_0)$ is defined to be the set of points $z\in\cbar$ such that there exists a sequence of positive integers $\{n_k\}\to\infty$ as $k\to\infty$ such that $\{f^{n_k}(z_0)\}\to z$ as $k\to\infty$.

\begin{proposition}\label{Koebe1}
(1) For any compact set $E\subset\cbar$ with $E\cap\overline{\pi_f^{-1}(\sA)}=\emptyset$, there is a Koebe space system $\{N(T)\}$ such that $N(T)\cap E=\emptyset$ for every band-tree $T$.

(2) Let $z_0\in\cbar$ be a point such that $\omega(z_0)$ is disjoint from all periodic band-trees. Then there is a Koebe space system $\{N(T)\}$ such that $z_0\notin N(T)$ for every band-tree $T$.
\end{proposition}

\beginp
(1) This is a direct consequence of Proposition \ref{Koebe} (e).

(2) Choose the Koebe space system $\{N(T)\}$ such that for every periodic band-tree $T_0$, $N(T_0)$ is disjoint from the closure of the orbit of $z_0$. Then $z_0\notin N(T)$ for every band-tree $T$.

\qed

\subsection{Nested neighborhoods of a skeleton}

Denote by $\BBB_n$ the union of all bands of level $k\le n$ for $n\ge 0$. Let $\phi_{t,n}$ be the normalized quasiconformal map of $\cbar$ whose Beltrami differential $\mu(\phi_{t,n})$ is the truncation of $\mu_t$ up to the $n$-level bands, i.e.
$$
\mu(\phi_{t,n})=\left\{
\begin{array}{ll}
\mu(\phi_t) & \text{on }\BBB_n, \\
0 & \text{elsewhere}.
\end{array}\right.
$$

\REFLEM{piece}
For any fixed $t\ge 0$, $\{\phi_{t,n}\}$  converges uniformly to $\phi_t$ as $n\to\infty$.
\ENDLEM

\beginp
For any fixed $t\ge 0$, the sequence $\{\phi_{t,n}\}$ is uniformly quasiconformal and hence a normal family. Let $\psi_t$ be one of its limits; then the Beltrami differential of $\psi_t\circ\phi_t^{-1}$ vanishes everywhere. So $\psi_t\circ\phi_t^{-1}$ is conformal on $\cbar$ and thus $\psi_t=\phi_t$ by the normalization.
\qed

\vskip 0.24cm

The main objective in this sub-section is to prove the following lemma.

\REFLEM{nested-annuli}
Fix $n\ge 0$. Let $S$ be a skeleton of level $j\le n$. Then there exist a constant $M>0$, a sequence of positive numbers $\{t_k\}\to\infty$ as $k\to\infty$, and a sequence of nested disks $\{U_k\}$ in $\cbar$ such that:

(a) $\overline{U_{k+1}}\subset U_k$,

(b) $\cap_{k\ge 0}U_k=\wh S$, and

(c) $\md \phi_{t,n}(U_k\smm\overline{U_{k+1}})>M$ for $t\ge t_k$.
\ENDLEM

\beginp
Let $N(T)$ be a Koebe space system such that $N(T)$ for all the band-trees $T$ of level $k\le n$ are pairwise disjoint. Begin by assuming that $j=0$, i.e., $S$ is periodic. Then there is an integer $p\ge 1$ such that each periodic arc on $S$ is fixed under $f^p$. Then each periodic point $x\in S$ is fixed under $f^p$. Let $g=f^p|_{N(S)}$.

\vskip 0.24cm

{\bf Nested sub-bands and the choice of $\{t_k\}$}. For each component $A_i\subset\sA$, let
$$
A_{i,0}=A_i\Supset A_{i,1}\Supset\cdots\Supset A_{i,k}\Supset\cdots
$$
be a sequence of nested annuli such that $\partial A_{i,k}$ consists of horizontal circles in $A_i$, and the equator $e(A_{i,k})=e(A_i)$ satisfies. By Proposition \ref{A(t)}, there exists a sequence of constants $\{t_k\ge 0\}$ such that $A_i(t_k)\subset A_{i,k}$. Thus for each component $E$ of $A_{i,k}\smm\overline{A_{i, k+1}}$, the modulus of $\Phi_t(E)$ is a constant for $t\ge t_k$.  Therefore we may further assume that
$$
\md \Phi_t(E)=3^k\quad\text{for }t\ge t_k.
$$

Denote by $\sA_k$ the union of $A_{i,k}$ for all components $A_i$ of $\sA$. Then $\cap_{k\ge 0}\sA_k=e(\sA)$. Let $B_0$ be the union of all bands intersecting  $S$. Let $B_k=B_0\cap\pi^{-1}_f(\sA_k)$. Then the $\{B_k\}$ $(k\ge 0)$ are nested sub-bands with $\cap_{k\ge 0}\overline{B_k}=S$.

\vskip 0.24cm

{\bf Nested calyxes}. Let $x\in S$ be a parabolic fixed point of $g$. Pick a calyx $\WWW_0(x)\subset N(S)$ of $g$ at the point $x$ such that $\overline{\WWW_0(x)}\cap\overline{B_0}= \{x\}$. Noticing that $\pi(\WWW_0(x))$ is the disjoint union of finitely many once-punctured disks, there is a holomorphic map $\xi_x: \pi(\WWW_0(x))\to\DD^*$ such that $\xi_x$ is a conformal map on each component.

Pick a constant $r\in (0,1)$. Set $r_0=1, r_1=r$ and $r_k=r^{1+3+\cdots +3^{k-1}}$ for $k\ge 2$. Let
$$
\WWW_k(x)=\pi^{-1}\circ\xi_x^{-1}(\DD^*(r_k))\cap\WWW_0(x)\quad\text{ for }k\ge 0.
$$
Then the $\{\WWW_k(x)\}$ are nested calyxes at the point $x$. Let $\WWW_k=\cup\WWW_k(x)$, the union over all parabolic periodic points $x\in S$.

\vskip 0.24cm

{\bf Nested attracting/repelling flowers}. Let $x\in S$ be a parabolic fixed point of $g$. By Proposition \ref{regular}, there exist a regular attracting flower $\VVV^+_0(x)$ and a regular repelling flower $\VVV^-_0(x)$ of $g$ at the point $x$ in $N(S)$, such that:

$\bullet$ for each band $B\subset B_0$, $\VVV^{\pm}_0(x)$ is disjoint from $B$ if $x\notin\overline{B}$; and $\pi(\partial \VVV^{\pm}_0(x)\cap B)$ is a vertical arc in $\pi(B)$ if $a(B)=x$ (or $r(B)=x$),

$\bullet$ for each component $\beta$ of $\partial \VVV^{\pm}_0(x)\cap\WWW_0(x)$, $\xi_x\circ\pi(\beta)$ is a straight line in $\DD^*$, and

$\bullet$ The sets $\VVV^{\pm}_0(x)$ for all parabolic periodic points $x\in S$, are pairwise disjoint .

Define $\VVV^{\pm}_{k+1}(x)=g^{\pm 3^{k}}(\VVV^{\pm}_{k}(x))$ for $k\ge 0$ inductively. Then $\{\VVV^{\pm}_k(x)\}$ $(k\ge 0)$ are nested attracting and repelling flowers at the point $x$. Let $\VVV^{\pm}_{k}=\cup\VVV^{\pm}_k(x)$, the union over all parabolic fixed points $x\in S$.

\vskip 0.24cm

{\bf Nested disks at hyperbolic points}. Let $x\in S$ be an attracting or repelling fixed point of $g$. Pick a disk $D_x\subset N(S)$ with $x\in D_x$ satisfying the following conditions:

$\bullet$ $\overline{g(D_x)}\subset D_x$ (or $\overline{g^{-1}(D_x)}\subset D_x$ if $x$ is repelling).

$\bullet$ $D_x\cap\VVV^{\pm}_0=\emptyset$.

$\bullet$ The sets $D_x$ for all periodic attracting and repelling points $x\in S$, are pairwise disjoint.

Note that the quotient space $\TT_x=D_x/\langle g\rangle$ is a torus. Denote by $\pi_x$ the natural projection from $D_x$ to $\TT_x$. Then $\pi_x(D_x\cap B_0)$ are mutually isotopic annuli. Pick a Jordan curve $\gamma\subset\TT_x$ such that the intersection of $\gamma$ with each annulus is a vertical line.
Then $\gamma$ has a lift in $D_x$ that bounds a disk $\VVV_0^h(x)\subset D_x$.

Define $\VVV_{k+1}^h(x)=g^{3^k}(\VVV_k^h(x))$ (or $g^{-3^k}(\VVV_k^h(x))$ if $x$ is repelling) for $k\ge 0$ inductively. Then $\{\VVV_{k}^h(x)\}$ $(k\ge 0)$ forms nested disks with $\cap_{k\ge 0}\VVV_k^h(x)=\{x\}$. Let $\VVV^h_k=\cup\VVV_k^h(x)$, the union over all attracting and repelling points $x\in S$.

\vskip 0.24cm

{\bf Nested disk neighborhoods of $\wh S$}. Let $\VVV_k=\VVV^h_k\cap\VVV^{+}_k\cup\VVV^-_k$ and let
$$
U_k=B_k\cup\VVV_k\cup\WWW_k\cup \wh S.
$$
Then $\{U_k\}$ $(k\ge 0)$ are nested disks with $\overline{U_{k+1}}\subset U_k$ and $\cap_{k\ge 0} U_k=\wh S$.

\begin{proposition}\label{quadrilateral}
Let $Q$ be a component of $(B_{k-i}\smm\overline{B_{k+1}})\cap(\VVV_k\smm\overline{\VVV_{k+1}})$ where $k\ge 0$ and $i=0,1$. We obtain a topological quadrilateral by setting the horizontal sides to be
$$
\alpha_1=(B_{k-i}\smm\overline{B_{k+1}})\cap\partial\VVV_k\quad\text{and}\quad \alpha_2=(B_{k-i}\smm\overline{B_{k+1}})\cap\partial\VVV_{k+1}.
$$
Then when $t\ge t_{k+1}$,
$$
\md\phi_{t,n}(Q)=\left\{
\begin{array}{ll}
1 & \text{ if}\quad i=0, \\
3/4 & \text{ if}\quad i=1.
\end{array}\right.
$$
\end{proposition}

\beginp
By the choice of $\partial\VVV_k$, $\pi_f(\alpha_1)=\pi_f(\alpha_2)$ is a vertical line in the annulus $\pi_f(B_0)$. By the choice of $B_k$, $\pi_f(\partial B_k)$ are horizontal circles in $\pi_f(B_0)$. Thus $Q\cap g^s(\partial\VVV_k)$ is a horizontal line in $Q$ for $s\in\ZZ$ whenever it is not empty, and $\partial B_s\cap Q$ is a vertical line in $Q$ for $k<s<l$ if it exists. Recall that for each component $E$ of $B_k\smm\overline{B_{k+1}}$, $\md \pi_{f_t}\circ\phi_t(E)=3^k$ for $t\ge t_{k+1}$, and $\VVV_{k+1}^h(x)=g^{3^k}(\VVV_k^h(x))$ (or $g^{-3^k}(\VVV_k^h(x))$ if $x$ is repelling) for $k\ge 0$. Thus when $t\ge t_{k+1}$,
$$
\begin{cases}
\md\phi_{t,n}(Q)=\dfrac{3^k}{3^k}=1 & \text{if $i=0$}, \\
\md\phi_{t,n}(Q)=\dfrac{3^k}{3^{k-1}+3^k}=\dfrac{3}{4} & \text{if $i=1$}.
\end{cases}
$$
\qed

{\bf Modulus control}. Let $N_k=U_k\smm\overline{U_{k+1}}$. In the following we always assume $t\ge t_k$. Then $\md \phi_{t,n}(N_k)$ is independent of $t$. We want to show that there exists a constant $M>0$ such that
$$
\md \phi_{t,n}(N_k)\ge M.
$$
The intersections of $N_k$ with $B_k$, with $\WWW_k$ and with $\VVV_k$ are topological quadrilaterals. Their moduli are bounded from below by the following discussion.

\begin{figure}[htbp]
\begin{center}
\includegraphics[scale=0.9]{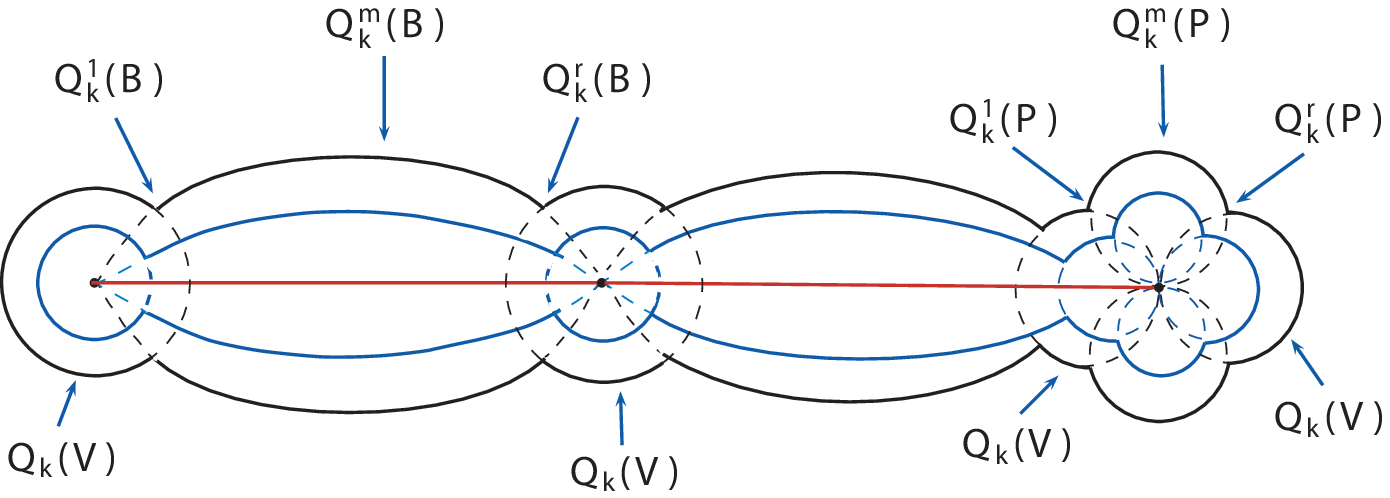}
\end{center}
\begin{center}{\sf Figure 5. Nested disks.}
\end{center}
\end{figure}

{\it Case 1}. For each component $B$ of $B_0\smm\wh S$, there is a unique component of $B_k\cap N_k$ contained in $B$ (denote it by $Q_k(B)$) which becomes a quadrilateral by setting its horizontal sides to be $\partial Q_k(B)\cap\partial B_k$ and $\partial Q_k(B)\cap\partial B_{k+1}$. Let $M_k(B)=\md \phi_{t,n}(Q_k(B)) $.

The set $Q_k(B)$ is cut into three sub-quadrilaterals by $\partial\VVV_k$. Denote the left, right and middle pieces by $Q_k^l(B)$, $Q_k^r(B)$ and $Q_k^m(B)$, respectively. Denote by $M_{k,l}(B)$, $M_{k,r}(B)$ and $M_{k,m}(B)$ the moduli of their $\phi_{t,n}$-images, respectively. Note that $\phi_{t,n}(\partial \VVV_k\cap Q_k(B))$ are still vertical arcs in $\phi_{t,n}(Q_k(B))$. By the choices of $B_k$ and $\VVV_k$, we have:
$$
M_{k,m}(B)=3 M_{k-1}(B).
$$
By Proposition \ref{quadrilateral},
$$
M_{k,l}(B)=M_{k, r}(B)=1\quad\text{for }k\ge 0.
$$
Thus
$$
\dfrac{1}{M_{k}(B)}=\dfrac{1}{3M_{k-1}(B)}+\dfrac{1}{M_{0,l}(B)}+\dfrac{1}{M_{0,r}(B)}=\dfrac{1}{3}\frac{1}{M_{k-1}(B)}+2.
$$
It follows that
\begin{align*}
\dfrac{1}{M_{k}(B)} & =\dfrac{1}{3^k}\frac{1}{ M_0(B)}+2\left(1+\frac{1}{3}+\cdots+\dfrac{1}{3^{k-1}}\right) \\
& \le \dfrac{1}{M_0(B)}+3=: \dfrac1{M_1}.
\end{align*}

{\it Case 2}. For each component $W$ of $\WWW_0\smm\wh S$, there is a unique component of $\WWW_k\cap N_k$ contained in $W$ (denote it by $Q_k(W)$) which becomes a quadrilateral by setting its horizontal sides to be $\partial Q_k(W)\cap\partial\WWW_k$ and $\partial Q_k(W)\cap\partial\WWW_{k+1}$. Similarly to Case 1, we may show that there exists a constant $M_2>0$ such that
$$
\md \phi_{t,n}(Q_k(W))\ge M_2\quad\text{for }k\ge 0.
$$

{\it Case 3}. For each component $V$ of $\VVV_0\smm\wh S$, there is a unique component of $\VVV_k\cap N_k$ contained in $V$ (denote it by $Q_k(V)$) which becomes a quadrilateral by setting its horizontal sides to be $\partial Q_k(V)\cap\partial\VVV_k$ and $\partial Q_k(V)\cap\partial \VVV_{k+1}$. Let $M_k(V)=\md \phi_{t,n}(Q_k(V))$.

For $k\ge 1$, there are exactly two components of $B_{k-1}$ intersecting $Q_k(V)$. Denote the intersections by $Q_k^\ell(V)$ and $Q_k^r(V)$ for $k\ge 0$. These are quadrilaterals whose vertical sides lie on $\partial B_{k+1}$ and $\partial B_{k-1}$. By Lemma \ref{quadrilateral}, we have
$$
\md\phi_{t,n}(Q_{k}^\ell(V))=\md\phi_{t,n}(Q_{k}^r(V))=3/4\quad \text{for }k\ge 1.
$$

Set $Q_{k}^m(V)=Q_k(V)\smm B_{k}$ for $k\ge 0$. By the choice of $\VVV_k$, we get three pairwise disjoint sub-quadrilaterals of $Q_{k}^m(V)$ whose moduli equal to that of $Q_{k-1}(V)$, and whose vertical sides as subsets of the vertical sides of $Q_k^m(V)$. Then we can apply the standard Gr\"{o}tzsch inequality to get
$$
\md\phi_{t,n}(Q_{k}^m(V))\ge 3 \md\phi_{t,n}(Q_{k-1}(V))=3 M_{k-1}(V)\text{ for }k\ge 1.
$$
Let
$$
Q_k^{\ell m}(V)=Q_k^{\ell}\cap Q_k^m(V)\text{ and }
Q_k^{rm}(V)=Q_k^{r}\cap Q_k^m(V).
$$
Also by Proposition \ref{quadrilateral},
$$
\md\phi_{t,n}(Q_{k}^{\ell m}(V))=\md\phi_{t,n}(Q_{k}^{rm}(V))=1.
$$
As above, $\phi_{t,n}(\partial Q_k^m(V)\cap Q_k^r(V))$ is a vertical line in $\phi_{t,n}(Q_k^r(V))$; and $\phi_{t,n}(\partial Q_k^m(V)\cap Q_k^\ell(V))$ is a vertical line in $\phi_{t,n}(Q_k^\ell(V))$. By Lemma \ref{mod4}, we have
\begin{align*}
\frac{1}{M_{k}(V)} & \le\frac{1}{3}\frac{1}{M_{k-1}(V)}+\frac{4}{3}+\frac{4}{3} \\
& =\frac{1}{3^k}\frac{1}{M_0(V)}+\left(1+\frac{1}{3}+\cdots+\frac{1}{3^{k-1}}\right)\frac{8}{3} \le \frac{1}{M_0(V)}+4=:\frac{1}{M_3}.
\end{align*}

The numbers of components of $N_k\cap B_k$, $N_k\cap\WWW_k$ and $N_k\cap\VVV_k$ are independent of $k\ge 0$; denote them by $n_1, n_2$ and $n_3$, respectively. Applying Lemma \ref{mod3}, we conclude that
$$
\dfrac 1{\md N_k}\le\frac{n_1}{M_1}+\frac{n_2}{M_2}+\frac{n_3}{M_3}=:\frac{1}{M}.
$$

\vskip 0.24cm

Now suppose that $S$ is a skeleton of level $j$ with $0<j\le n$. Let $U_k$ be the domain chosen as above for the periodic skeleton $f^j(S)$. Let $U_k(S)$ be the component of $f^{-j}(U_k)$ containing $S$; then
$$
f^j: U_0(S)\smm\overline{U_k(S)}\to U_0\smm\overline{U_k}
$$
is a covering of degree $d\ge 1$. Thus
$$
\md\phi_{t,n}(U_0(S)\smm\overline{U_k(S)}) =\frac 1d \md\phi_{t,n}(U_0\smm\overline{U_k}).
$$
This completes the proof.
\qed

\section{Simple pinching}

In this section, we will prove Theorem \ref{s-pinching}. Let $f$ be a geometrically finite rational map and let $\sA\subset\sR_f$ be a starlike multi-annulus. Let $f_t=\phi_t\circ f\circ\phi_t^{-1}$ $(t\ge 0)$ be the simple pinching path supported on $\sA$. Recall that the map $\phi_t$ is normalized by fixing three distinct points $z_1, z_2$ and $z_3$ in $U_0\cap\PPP_f$, where $U_0$ is a component of $\cbar\smm\overline{\pi_f^{-1}(\sA)}$ such that both $U_0\cap\PPP_f$ and $\partial U_0\cap\JJJ_f$ are infinite sets. By making a holomorphic conjugacy, we may assume that $z_3=\infty$, for the sake of simplicity.

\subsection{The angular space system}

Recall that for each component $A$ of $\sA$, $A'\Subset A$ is an annulus essentially contained in $A$ whose boundary consists of horizontal curves in $A$. Denote by $A''$ the annulus bounded by the equators of the two annular components of $A\smm\overline{A'}$. Then $A'\Subset A''\Subset A$, and $\partial A''$ consists of smooth curves. Denote by $\sA''$ the union of $A''$ for each component $A$ of $\sA$. For each band $B$, we denote by $B''$ the component of $\pi^{-1}(\sA'')$ contained in $B$.

Let $\{N(T)\}$ be a Koebe space system of the band-trees. For each periodic band $B$, pick a smooth disk $U_B\ni a(B)$ such that:

$\bullet$ $\partial U_B\cap\PPP_f=\emptyset$,

$\bullet$ $\overline{U_B}\subset N(T)$ if $B\in T$,

$\bullet$ $\overline{f(U_B)}\subset U_{f(B)}$ and

$\bullet$ $\partial U_B\cap\partial B''$ contains exactly two points.

Let $D(B)=B''\cup U_B$. Then $\{D(B)\}$ are pairwise disjoint disks, due to the assumption that the multi-annulus $\sA$ is starlike. For each level $n$ band $B_n$, $f^n(B_n)=B$ is a periodic band. Let $D(B_n)$ be the component of $f^{-n}(D(B))$ containing $a(B_n)$. Then $D(B)\Subset N(B)$ and $\{D(B)\}$ are pairwise disjoint disks for all bands $B$.

For each band $B$, denote by $B'\subset B$ the band of $\sA'$. Then $B'\Subset D(B)\cup r(B)$. We call $D(B)$ the {\bf angular space}\ of $B'$, and $\{D(B)\}$ an {\bf angular space system}.

\begin{figure}[htbp]
\begin{center}
\includegraphics[scale=0.5]{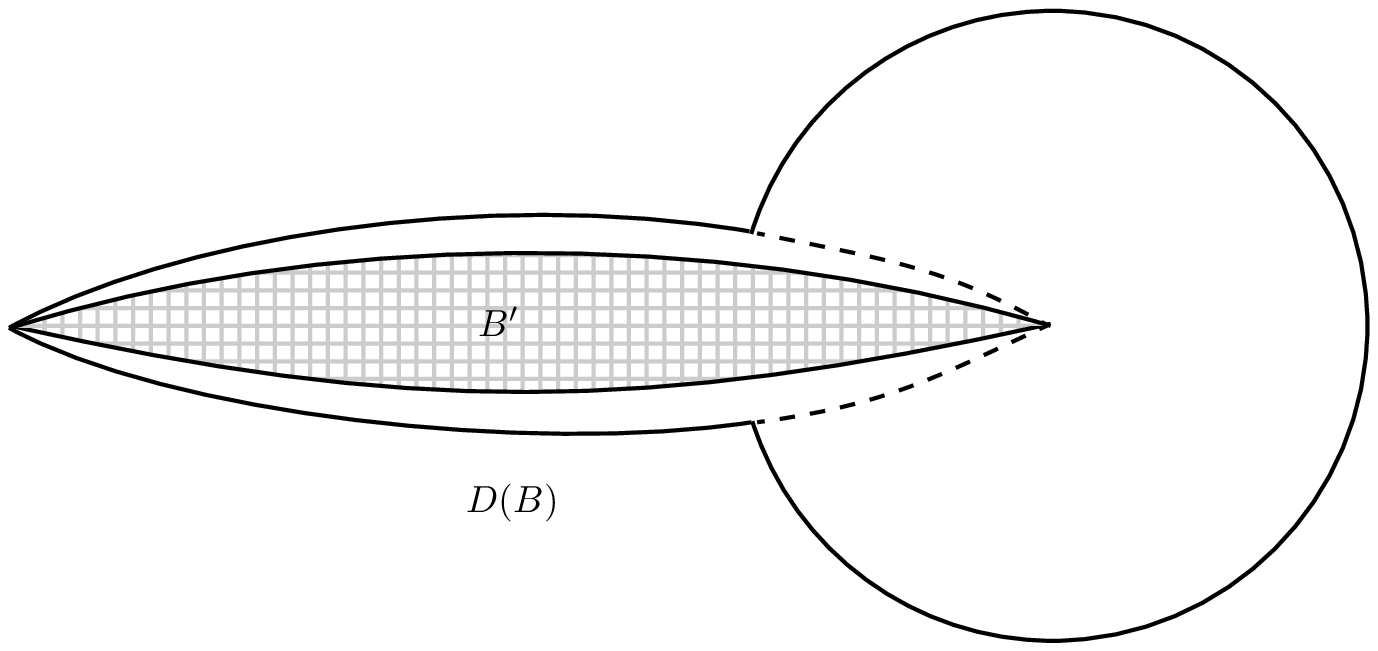}
\end{center}
\begin{center}{\sf Figure 6. The angular space $D(B)$ of $B'$.}
\end{center}
\end{figure}

\REFLEM{angle}
There is a constant $K_0<\infty$ such that for any band $B$, any point $z\in\partial D(B)\smm r(B)$ and any holomorphic injection $\zeta$ from $N(B)$ into $\CC$, there exists an arc $\delta\subset D(B)-\overline{B'}$ which joins $z$ to $r(B)$, such that
$$
L(\zeta(\delta))<K_0\cdot \text{\rm dist}(\zeta(z), \zeta(B')),
$$
where $L(\cdot)$ denotes the Euclidean length and $\text{\rm dist}(\cdot,\cdot)$ denotes the Euclidean distance.
\ENDLEM

\beginp We only need to prove the lemma for bands $B$ with $r(B)\in f^{-1}(\PPP_f)$ and for the holomorphic injection $\zeta=\text{id}$, by Koebe distortion Theorem.

Let $B$ be a periodic band with period $p\ge 1$. Since $r(B)$ is repelling, there exists a disk $V\subset N(B)$ with $r(B)\in V$ such that $\partial V$ intersects $\partial B''$ at exactly two points on $\partial D(B)$, $f^p$ is injective on $V$ and $V\Subset f^p(V)$.

Let $\lambda$ be the multiplier of $f^p$ at the fixed point $r(B)$. Then $|\lambda|>1$. By the Linearization Theorem, there exists a conformal map $\psi:\,f^p(V)\to\psi(f^p(V))\subset\CC$ such that $\psi(r(B))=0$ and $\psi\circ f^p(z)=\lambda\cdot\psi(z)$ on $V$.

Let $\beta=\partial D(B)\smm\{r(B)\}$. Denote by $\beta_1, \beta_2$ the two arcs of $(\partial B''\smm\{r(B)\})\cap V$. Then $f^p(\beta_i)\supset\beta_i$. Let $\gamma_i=\psi(\beta_i)$; then one endpoint of $\gamma_i$ is the origin and the other endpoint of $\gamma_i$, denote it by $w_i$, is on $\psi(\partial V)$.

For any two distinct points $w, w'\in\overline{\gamma_i}$, denote by $\gamma_i(w, w')$ the arc on $\gamma_i$ with endpoints $w$ and $w'$.
Since $\beta_i$ is smooth,  the Euclidean length $L(\gamma_i(w_i/\lambda, w_i))$ is bounded. Moreover,
$$
L\left(\gamma_i\left(\frac{w_i}{\lambda^2}, \frac{w_i}{\lambda}\right)\right)=\frac{1}{|\lambda|}L\left(\gamma_i\left(\frac{w_i}{\lambda}, w_i\right)\right).
$$
Thus $L_0=L(\gamma_i(0,w_i))$ is bounded and for any point $w\in\gamma_i(0, w_i)$,
$$
L(\gamma_i(0,w))=|\lambda|\cdot L(\gamma_i(0,w/\lambda)).
$$

There exists a constant $r_0>0$ such that $\DD(w,r_0)\subset f^p(V\smm B')$ for all points $w\in\gamma_i(w_i/\lambda, w_i)$. For each $k\ge 1$ and any point $w\in\gamma_i(w_i/\lambda^{k+1}, w_i/\lambda^k)$, since $\lambda\circ\psi(V\cap B')=\psi(f^p(V)\cap B')$, we have $\DD(w, r_0/|\lambda|^k)\subset f^p(V\smm B')$.

By the Koebe Distortion Theorem, there exists a constant $K_1<\infty$ such that for any point $z\in\beta_i$,
$$
L(\beta_i(r(B), z))<K_1\,\text{\rm dist}_e(z, B').
$$

For each point $z\in\beta\smm V$, $\text{\rm dist}(z, B')$ is bounded from below, while $L(\beta(r(B), z))$ is bounded from above since $\partial D_{a(B)}$ is smooth. Thus there exists a constant $K_2$ such that
$$
L(\beta(r(B), z))<K_2\,\text{\rm dist}_e(z, B').
$$
Combining the two inequalities above, we prove the lemma for the periodic band $B$.

Now let $B_k$ be a band of level $k\ge 1$ with $r(B_k)\in f^{-1}(\PPP_f)$. Let $V$ be the disk defined as above for the periodic band $B=f^k(B_k)$ and let $p\ge 1$ be the period of $B$. Let $U, W$ be the components of $f^{-k}(V)$ and $f^{-k}(f^p(V))$ that contain $B_k$, respectively. Then there exists a conformal map $\xi: U\to W$ such that $f^k\circ\xi=f^k\circ f^p$ on $U$ and $\xi(B_k\cap U)=B_k\cap W$. This shows that $r(B_k)$ is a repelling fixed point of $\xi$. By the same argument as the above, we may verify the lemma for the band $B_k$.
\qed

\subsection{Modulus control}

Applying Lemma \ref{mod2} and Lemma \ref{angle}, we will control the modulus of $\phi_t(A)$ for certain annuli $A\subset\cbar$ in the following lemmas.

\REFLEM{mod-phi-1}
Let $A\subset\cbar$ be an annulus which contains $\overline{\pi_f^{-1}(\sA)}$. Then there exists a constant $K>1$ such that
$$
\frac{1}{K}\md A \le\md \phi_{t}(A) \le K\md A
$$
for $t\ge 0$.
\ENDLEM

\beginp
Since $\overline{\pi_f^{-1}(\sA)}\subset A$, there exists a Koebe space system $\{N(T)\}$ such that $N(T)\subset A$ for every band-tree $T$, by Proposition \ref{Koebe1}. Let $D(B)$ be an angular space system such that $\overline{D(B)}\subset N(B)$.

Let $\rho$ be an extremal metric on the annulus $A$. It can be chosen to be $\rho(z)=|(\log\chi)'(z)|$, where $\chi: A\to\AA(r)$ is a conformal map for some $r>1$. For each $n\ge 0$, denote by $\BBB'_n$ the union of bands of the multi-annulus $\sA'$ with level $k\le n$ . Define
$$
\rho_n(z)=\begin{cases}
0 & \text{on the closure of }\BBB'_n, \\
\rho(z) & \text{otherwise}.
\end{cases}
$$

Let $B$ be a band of level $k\le n$ and let $B'\subset B$ be a band of $\sA'$. Let $\alpha$ be an arc in $D(B)$ which joins two boundary points $z_1, z_2\in\partial D(B)$. If $\alpha\cap\overline{B'}=\emptyset$, then $L(\rho_n, \alpha)=L(\rho,\alpha)$.

Now we assume that $\alpha$ intersects $\overline{B'}$. Consider the map $\log\chi(z)$; it has a univalent branch $\zeta(z)$ on the disk $N(B)$. By Lemma \ref{angle}, there exist two arcs $\delta_1,\delta_2\subset D(B)-\overline{B'}$ which join $z_1, z_2$ with $r(B)$, respectively, such that
$$
L(\zeta(\delta_i))<K_0\cdot \text{\rm dist}(\zeta(z_i), \zeta(B'))\text{ for }i=1,2,
$$
where $K_0<\infty$ is the constant from Lemma \ref{angle}. Set $\delta=\delta_1\cup\{r(B)\}\cup\delta_2$. Then $\delta\subset (D(B)\smm\overline{B'})\cup\{r(B)\}$ and
$$
L(\rho, \delta)=L(\zeta(\delta_1))+L(\zeta(\delta_2))<K_0\cdot L(\rho_n, \alpha).
$$
It follows that
$$
\width(\rho_n, A)>\frac{1}{K_0}\width(\rho, A)\quad\text{and}\quad
\height(\rho_n, A)>\frac{1}{K_0}\height(\rho, A).
$$
Let $\rho_{t,n}$ be the pushforward of $\rho_n$ by $\phi_{t,n}$, i.e.,
$$
\rho_{t,n}(w)=\begin{cases}
0, & \text{on the closure of }\phi_{t,n}(\BBB'_n), \\
|(\phi_{t,n}^{-1})'(w)|\rho_n(\phi_{t,n}^{-1}(w)) & \text{otherwise}.
\end{cases}
$$
Then $\rho_{t,n}$ is a conformal metric on $\phi_{t,n}(A)$ since $\phi_{t,n}$ is holomorphic in $\cbar\smm\overline{\BBB'_n}$. Therefore
$$
\begin{cases}
\area(\rho_{t,n}, \phi_{t,n}(A))=\area(\rho_n, A)\le\area(\rho,A), \\
\width(\rho_{t,n}, \phi_{t,n}(A))=\width(\rho_n, A)>\dfrac{1}{K_0}\width(\rho, A)\quad\text{and} \\
\height(\rho_{t,n}, \phi_{t,n}(A))=\height(\rho_n, A)>\dfrac{1}{K_0}\height(\rho, A).
\end{cases}
$$
By Lemma \ref{mod2}, we have
$$
\dfrac1{K_0^2} \md A <\md \phi_{t,n}(A) <K_0^2\md A .
$$
Letting $n\to\infty$ and applying Lemma \ref{piece} and Theorem \ref{continuity} completes the proof.
\qed

\vskip 0.24cm

From Lemma \ref{mod-phi-1} and the normalization condition, we have:

\REFCOR{non-degenerate} For each domain $D$ compactly contained in $\cbar\smm \overline{\pi_f^{-1}(\sA)}$, there exists a constant $C>0$ such that $\text{diam}_s \phi_t(D)>C$ for $t\ge 0$.
\ENDCOR

\REFLEM{mod-phi-2} Let $S$ be a skeleton. Then for any disk $U\supset S$ and any constant $M>0$, there exist a disk $D\supset S$ with $D\Subset U$ and a constant $t_0>0$, such that
$$
\md \phi_{t}(U\smm\overline{D})\ge M \quad\text{ for }t\ge t_0.
$$
\ENDLEM

\beginp
We only need to prove the lemma for periodic skeletons. Let $S$ be a periodic skeleton. Let $\{N(T)\}$ be a Koebe space system and let $\{D(B)\}$ be an angular space system with $\overline{D(B)}\subset N(B)$. Since $N(B)$ is disjoint from periodic bands for each band $B$ with level $n\ge 1$, by the Koebe Distortion Theorem, there exists a constant $M_1<\infty$ such that for any band $B$ with level $n\ge 1$,
$$
\text{diam}\,\phi_{t,0}(D(B))\le M_1\,\text{dist}\,(\phi_{t,0}(D(B)), \partial\phi_{t,0}(N(B))).
$$

By Lemma \ref{nested-annuli}, for any constant $M_2<\infty$, there exist a constant $t_0\ge 0$ and disks $V_1, V_2$ with $S\subset V_2\Subset V_1\subset U\cap N(S)$, such that
$$
\md \phi_{t,0}(V_1\smm\overline{V_2})\ge M_2\text{ for }t\ge t_0.
$$

Given any $t\ge t_0$, by Lemma \ref{mod1}, there exists a round annulus $A_t=\AA(w_0, r_1, r_2)$ with center $w_0\in\phi_{t,0}(S)$ and with $A_t$ contained essentially in $\phi_{t,0}(V_1\smm\overline{V_2})$ such that
$$
\md A_t=\frac{\log(r_2/r_1)}{2\pi} \ge M_2-\dfrac{5\log 2}{2\pi}.
$$

Denote by $C_2$ and $C_1$ the outer and the inner boundary components of $A_t$, respectively. Since $w_0\notin\phi_{t,0}(N(B))$, if $\phi_{t,0}(D(B_1))\cap C_2\neq\emptyset$ for some band $B_1$ with level $n\ge 1$, then
$$
\text{\rm dist}(w_0, \phi_{t,0}(D(B_1)))\ge r_2/(1+M_1).
$$
Similarily,  if $\phi_{t,0}(D(B_2))\cap C_1\neq\emptyset$ for some band $B_2$ with level $n\ge 1$, then
$$
\text{\rm dist}(w_0, \phi_{t,0}(D(B_2)))\le r_1(1+M_1).
$$
Given any $n\ge 1$, the map $\phi_{t,n}\circ\phi_{t,0}^{-1}$ is holomorphic except on the closure of the union of all bands of level $1\le k\le n$.  Define a conformal metric on $A_t$ by
$$
\rho_n(w)=\begin{cases}
0 & \text{on the closure of }\phi_{t,0}(\BBB'_n), \\
\rho(w)=\dfrac{|dw|}{|w-w_0|} & \text{otherwise}.
\end{cases}
$$
Then
$$
\area(\rho_n, A_t)\le\area(\rho, A_t)=2\pi\log(r_2/r_1)=4\pi^2\md A_t.
$$
For any arc $\alpha$ in $A_t$ joining its two boundary components, as in the proof of Lemma \ref{mod-phi-1},
$$
L(\rho_n, \alpha)\ge\frac{1}{K_0}\left(\log\frac{r_2}{r_1}-2\log(1+M_1)\right),
$$
where $K_0$ is the constant in Lemma  \ref{angle}. Thus
$$
\md(\phi_{t,n}\circ\phi_{t,0}^{-1})(A_{t})\ge\frac{\height(\rho_n, A_t)^2}{\area(\rho_n, A_t)}
\ge \frac{\md A_{t}}{K_0^2}\left(1-\frac{\log(1+M_1)}{\pi\md A_t}\right)^2.
$$
Combining  these inequalities with the fact that $\phi_{t,0}^{-1}(A_{t})$ is contained in $V_1\smm\overline{V_2}$ essentially, we get
$$
\md \phi_{t,n}(V_1\smm\overline{V_2})\ge
\frac{M_2-(5\log 2)/(2\pi)}{K_0^2}\left(1-\frac{2\log(1+M_1)}{2\pi M_2-5\log 2}\right)^2\quad\text{for }t\ge t_0.
$$
Note that both constants $M_1$ and $K_0$ are independent of the choice of $n$. Let $n\to\infty$. We get the lemma by Theorem \ref{continuity} and Lemma \ref{piece}.
\qed

\REFLEM{mod-phi-3}
Let $z_0\in\cbar$ be a point which is not contained in any skeleton. Then for any disk $U\subset\cbar$ with $z_0\in U$ and any constant $M>0$, there exist a constant $t_0\ge 0$ and a disk $D\Subset U$ with $z_0\in D$, such that
$$
\md \phi_{t}(U\smm\overline{D})\ge M\quad\text{for }t\ge t_0.
$$
\ENDLEM

\beginp
If $z_0\in\FFF_f$, then $\{\phi_t\}$ is uniformly quasiconformal in a neighborhood of $z_0$ for $t\ge 0$. The lemma is trivial in this case. Now we assume that $z_0\in\JJJ_f$.

{\it Case 1}. Assume that $\omega(z_0)$ is disjoint from the closures of all periodic bands. By Proposition \ref{Koebe1}, there exists a Koebe space system $\{N(T)\}$ such that $z_0\notin N(T)$ for every band-tree $T$. Let $\{D(B)\}$ be an angular space system with $\overline{D(B)}\subset N(B)$. Then there exists a constant $M_1<\infty$ such that for each band $B$ and any holomorphic injection $\zeta: N(B)\to\CC$,
$$
\text{diam }\zeta(D(B))\le M_1\text{ dist\,}(\zeta(D(B)),\partial\zeta(N(B))).
$$
For any $M>0$, let $A$ be a round annulus with modulus $\md A=M_2$ and contained essentially in $U\smm\{z_0\}$. Using the argument in the proof of Lemma \ref{mod-phi-2}, we have a constant $K_0<\infty$ such that
$$
\md \phi_{t,n}(A)\ge\frac{M_2}{K_0^2}\left(1-\frac{\log(1+M_1)}{\pi M_2}\right)^2\quad\text{ for }t\ge 0.
$$
Note that both constants $M_1$ and $K_0$ are independent of the choice of $n$. Let $n\to\infty$. We get the lemma by Theorem \ref{continuity} and Lemma \ref{piece}.

{\it Case 2}. Assume that $r(S)\in\omega(z_0)$ for some periodic skeleton $S$. Let $\{N(T)\}$ be a Koebe space system. By Lemma \ref{mod-phi-2}, there exist a constant $M>0$, two disks $U_0$ and $U_1$ in $N(S)$ with $S\subset U_1\Subset U_0$ and a constant $t_0\ge 0$, such that
$$
\md \phi_t(U_0\smm\overline{U_{1}}) \ge M\quad \text{ for }t\ge t_0.
$$
Without loss of generality, we may assume that $S$ is fixed, i.e., $f(S)=S$. Then $f^{-1}(U_1)$ has exactly one component $U'$ such that $S\subset U'$. Since $r(S)$ is repelling, we may assume that $(U'\cap\JJJ_f)\subset U_1$ (recall the construction of $U_n$ in Lemma \ref{nested-annuli}).

Now let $f^{n_k}(z_0)\in U_1$ be the first return of the orbit $\{f^n(z_0)\}$. Then $f^{n_k-1}(z_0)\notin U_1$. Let $V_k, V'_k$ be the components of $f^{-n_k}(U_0)$ and $f^{-n_k}(U_1)$ that contain $z_0$ along its orbit, respectively. Then $V_k\smm\overline{V'_k}$ is an annulus around $z_0$, which shrinks to the point $z_0$ by Lemma \ref{shrinking}. Note that
$$
\md \phi_t(V_k\smm\overline{V'_k}) \ge\frac{M}{K}\quad \text{ for }t\ge t_0,
$$
where $K=\max_{k\ge 1}\deg_{z_k} f^k$ over all point $z_k$ with $f^k(z_k)=r(S)$ but $f^{k-1}(z_k)\neq r(S)$. Therefore, we may choose as many of them as possible such that they are pairwise disjoint, which forms the desired annulus.
\qed

\subsection{Proof of Theorem \ref{s-pinching}}

\REFLEM{normal-phi}
The family $\{\phi_{t}\}$ ($t\ge 0$) is equicontinuous. Let $\varphi$ be a limit of the family as $t\to\infty$. Then $\varphi(S)$ is a single point for each skeleton $S$. Conversely, for any point $w\in\cbar$, $\varphi^{-1}(w)$ is either a single point or a skeleton.
\ENDLEM

\beginp
We begin by proving that $\{\phi_{t}\}$ is equicontinuous. Pick a disk $D_0\subset\CC$ such that $\overline{\pi_f^{-1}(\sA)}\subset D_0$. From Lemma \ref{mod1} and Corollary \ref{non-degenerate}, we only need to prove that for any $M>0$ and point $z_0\in D_0$, there exists a disk $D_{z_0}\ni z_0$ such that $\md \phi_t(D_0\smm\overline{D_{z_0}})>M$ for all $t\ge 0$.

Assume that $z_0$ is contained in an $m$-level skeleton $S$. By Lemma \ref{mod-phi-2}, there exist a disk $D\supset S$ with $\overline{D}\subset D_0$ and a constant $t_0>0$, such that
$$
\md \phi_{t}(D_0\smm\overline{D})\ge M
$$
for $t\ge t_0$. Since $\{\phi_t, t\le t_0\}$ is uniformly quasiconformal, there exists a disk $D_{z_0}\ni z_0$ with $D_{z_0}\subset D$, such that $\md \phi_{t}(D_0\smm\overline{D_{z_0}})\ge M$ for $t\le t_0$. Combining these facts we get $\md \phi_{t}(D_0\smm\overline{D_{z_0}})\ge M$ for all $t\ge 0$.

Now we assume that $z_0$ is not contained in any skeleton. By Lemma \ref{mod-phi-3}, there exists a disk $D_{z_0}\ni z_0$ with $\overline{D_{z_0}}\subset D_0$ such that
$$
\md \phi_{t}(D_0\smm\overline{D_{z_0}}) \ge M
$$
for all $t\ge 0$.

Now we have proved that $\{\phi_{t}\}$ is equicontinuous. Let $\varphi$ be a limit of the family as $t\to\infty$. From Lemma \ref{quotient}, $\varphi$ is a quotient map of $\cbar$. By Lemma \ref{mod-phi-2}, for each skeleton $S$, $\varphi(S)$ is a single point. Conversely, for each point $w\in\cbar$, if $\varphi^{-1}(w)$ contains at least two points, we claim that $\varphi^{-1}(w)$ is a skeleton. Otherwise, let $z_1, z_2$ be two distinct points in $\varphi^{-1}(w)$ which are not contained in one skeleton. By Lemma \ref{mod-phi-2} and \ref{mod-phi-3}, there exist constants $t_0>0$ and $M>5\log 2/(2\pi)$, and an annulus $A\subset\cbar$, such that the two components $D_1, D_2$ of $\cbar\smm\overline{A}$ are disks which contain the points $z_1$ and $z_2$, respectively, and
$$
\md \phi_t(A) \ge M\quad\text{ for }t\ge t_0.
$$
Note that both $D_1$ and $D_2$ intersect $W$. Thus neither $\varphi(D_1)$ nor $\varphi(D_2)$ is a single point. Therefore, there is a positive distance between them and hence $\varphi(z_1)\neq\varphi(z_2)$. This is a contradiction.
\qed

\vskip 0.24cm

Let $\varphi$ be a limit of the family $\{\phi_t\}$ as $t\to\infty$. Let $B$ be a periodic band with period $p\ge 1$ and let $E$ be a component of $B\smm S$, where $S$ is the periodic skeleton with $S\subset\overline{B}$. Let $\chi: \pi_f(E)\to\AA(1/r,1)$ be a conformal map such that $|\chi\circ\pi_f(z)|\to 1$ as $z\to\partial B\smm S$. Let $W=\varphi(E)$. Then $g=\varphi\circ f^p\circ\varphi^{-1}: W\to W$ is a well-defined conformal map.
From Proposition \ref{model}, we have:

\begin{lemma}\label{model-varphi}
There exists a universal covering $\pi: W\to\DD^*$ with $\pi(w_1)=\pi(w_2)$ if and only if $w_1=g^{k}(w_2)$ for some integer $k\in\ZZ$ such that the following diagram commutes,
\[
\xymatrix{
E\ \ar[d]_{\chi\circ\pi_f}\ar[r]^{\varphi} & \ W\ \ar[d]_{\pi} \\
\AA(1/r,1)\ \ar[r]^{w} & \ \DD^* }
\]
where $w$ is the map defined in Proposition \ref{model} (4).
\end{lemma}

\vskip 0.24cm

{\noindent\it Proof of Theorem \ref{s-pinching}}. From Lemmas \ref{normal-phi} and \ref{piece}, $\{\phi_{t}\}$ $(t\ge 0)$ is equicontinuous. Let $\{t_n\}$ be a sequence in $[0, \infty)$ with $\{t_n\}\to\infty$ as $n\to\infty$ such that $\{\phi_{t_n}\}$ converges uniformly to a quotient map $\varphi$.  Let $f_{t_n}=\phi_{t_n}\circ f\circ\phi_{t_n}^{-1}$. By Corollary \ref{non-degenerate}, $\varphi$ is injective on
$$
W:=\cbar\smm\overline{\pi_f^{-1}(e(\sA))}.
$$
So $\{f_{t_n}\}$ converges uniformly to a map $g$ on any compact set in $W$ and $\deg(g|_{\varphi(W)})=\deg f$. By Lemma \ref{rational-s}, $g$ is a rational map and $\{f_{t_n}\}$ converges uniformly to $g$ on $\cbar$.

Each skeleton $S$ intersects $\JJJ_f$ at exactly one point $r(S)$. So $\varphi$ is injective on $\JJJ_f$ and hence it is a homeomorphism from $\JJJ_f$ to $\varphi(\JJJ_f)$.

Since $\varphi(\JJJ_f)$ is a completely invariant perfect set and is contained in the closure of the periodic point set, we have $\varphi(\JJJ_f)=\JJJ_g$. Since $\JJJ_g\cap\PPP_g=\varphi(\JJJ_f\cap\PPP_f)$ and $f$ is geometrically finite, $g$ is a geometrically finite rational map.

If there exists another sequence $\{t'_n\}$ in $[0,\infty)$ with $\{t'_n\}\to\infty$ as $n\to\infty$ such that $\{\phi_{t'_n}\}$ converges uniformly to $\wt\varphi$, then $\{\phi_{t'_n}\circ f\circ\phi_{t'_n}^{-1}\}$ converges uniformly to another rational map $\wt g$. Set $\psi=\wt\varphi\circ\varphi^{-1}$. It is a well-defined homeomorphism and $\psi\circ g=\wt g\circ\psi$. Moreover, $\psi$ is holomorphic on the Fatou set $\FFF_g$. Thus $\psi$ is a conformal map of $\cbar$ by Theorem \ref{unicity}. By the normalization condition, we have $\psi=\text{id}$ and hence $\wt\varphi=\varphi$. Thus $\{\phi_t\}$ converges uniformly to $\varphi$ and $\{f_t\}$ converges uniformly to $g$ as $t\to\infty$.
\qed

\section{Simple plumbing}

In this section we will prove the following theorem. Theorem \ref{s-plumbing} and the sufficiency part of Theorem \ref{existence} are direct consequences of this theorem.

\REFTHM{s-plumbing-1}
Let $(G, \QQQ)$ be a marked semi-rational map with parabolic cycles in $\PPP'_G$ and $\#(\QQQ\smm\PPP_G)<\infty$. Suppose that $(G, \QQQ)$ has neither Thurston obstructions nor connecting arcs. Then there exist a marked rational map $(g, \QQQ_1)$, a sub-hyperbolic rational map $f$ and a simple pinching path $\{f_t\}$ ($t\ge 0$) starting from $f$ such that $\{f_t\}$ converges to $g$ uniformly on $\cbar$ as $t\to\infty$ and $(G, \QQQ)$ is c-equivalent to $(g, \QQQ_1)$.
\ENDTHM

\subsection{Simple plumbing surgery}

{\bf Step 1. The cut-glue surgery}.
Let $(G, \QQQ)$ be a marked semi-rational map with parabolic cycles in $\PPP'_G$ and $\#(\QQQ\smm\PPP_G)<\infty$. Denote by $Y\subset\PPP'_G$ the set of parabolic cycles. Pick a calyx for each cycle in $Y$ such that the closure of their union $\WWW$ is disjoint from $\QQQ\smm Y$. The quotient space $\WWW/\langle G\rangle$ is a disjoint union of once-punctured disks. Thus there is a natural holomorphic projection $\pi: \WWW\to\DD^*$ such that for each sepal $W$ of $\WWW$ with period $p\ge 1$, $\pi:\ W\to\DD^*$ is a universal covering and $\pi(z_1)=\pi(z_2)$ if and only if $z_1=G^{kp}(z_2)$ for some integer $k\in\ZZ$.

Given any $0<r<1$, let $\WWW(r)=\pi^{-1}(\DD^*(r))$ and $\RRR(r)=\WWW\smm\overline{\WWW(r)}$. Then $G(\RRR(r))=\RRR(r)$. Thus there is a conformal map $\tau:\ \RRR(r^2)\to\RRR(r^2)$ such that:

(i) $\tau(z)$ and $z$ are contained in the same attracting petal but in different sepals.

(ii) $\tau^2=\text{id}$ and $G\circ\tau=\tau\circ G$.

The map $\tau$ is unique up to composition with some iterates of $G$.

Define an equivalence relation on $\cbar\smm\overline{\WWW(r^2)}$ by $z_1\sim z_2$ if $\tau(z_1)=z_2$. Then the quotient space is a punctured sphere with finitely many punctures. Thus there exists a finite set $X\subset\cbar$ and a holomorphic surjective map
$$
p:\, \cbar\smm\overline{\WWW(r^2)}\to\cbar\smm X
$$
such that $p(z_1)=p(z_2)$ if and only if $z_1=\tau(z_2)$.

Let $\SSS=p(\partial\WWW(r)\smm Y)\cup X$. It is a finite disjoint union of trees whose vertex set is $X$. Let $\BBB=p(\RRR(r^2))$. It is a finite disjoint union of disks, and $\SSS\cup X\subset\overline{\BBB}$.

\begin{figure}[htbp]
\begin{center}
\vspace{-2cm}
\includegraphics[scale=0.8]{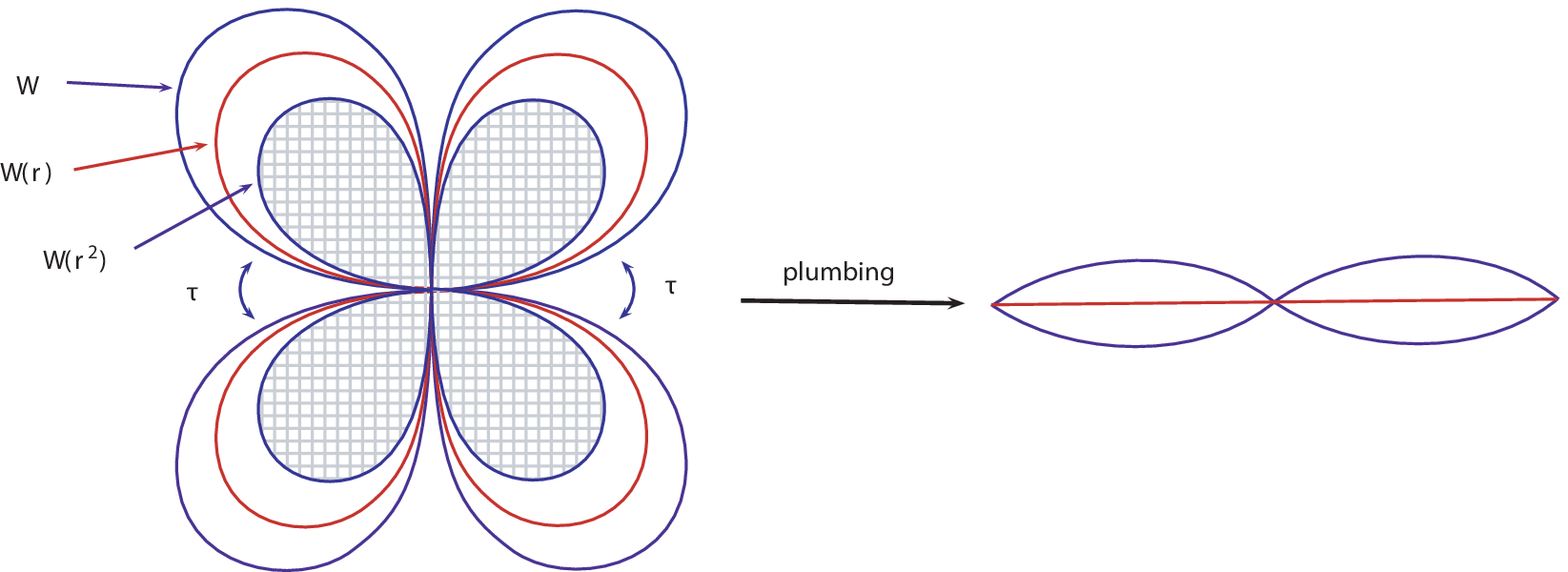}
\end{center}
\begin{center}{\sf Figure 7. Simple plumbing.}
\end{center}
\end{figure}

\vskip 0.24cm

{\bf Step 2. The induced map after surgery}.
The semi-rational map can be pushed forward to the quotient space $\cbar\smm\overline{\WWW(r^2)}$ except on $\WWW_1:=G^{-1}(\WWW)\smm\WWW$ since $G\circ\tau=\tau\circ G$ on $\RRR(r^2)$, i.e. there is continuous map $F_0:\,\cbar\smm p(\overline{\WWW_1})\to\cbar$ such that
$$
F_0\circ p=p\circ G\text{ on }\cbar\smm{G^{-1}(\overline\WWW)}.
$$
Since $G$ is holomorphic in a neighborhood of $\overline{\WWW}$, $F_0$ is holomorphic in a neighborhood of $\overline{\BBB}$. Obviously, $F_0(X)=X$, $F_0(\SSS)=\SSS$ and $F_0(\BBB)=\BBB$.

For each component $B$ of $\BBB$, $\overline{B}\cap X$ contains exactly two points. Let $p\ge 1$ be the period of $B$. The $\{F_0^{kp}(z)\}$ converges to a point in $\overline{B}\cap X$ for any point $z\in B$.  Denote it by $a(B)$. Since each attracting petal of $G$ at a point in $Y$ contains infinitely many points of $\PPP_G$, $a(B)$ is an accumulation point of $p(\PPP_G\smm Y)$.

Denote by $r(B)$ the other point in $\overline{B}\cap X$. By condition (i) in the definition of the map $\tau$, $r(B)$ is not an accumulation point of $p(\PPP_G\smm Y)$. Moreover, the closures of any two distinct components $B_1, B_2$ of $\BBB$, are either disjoint or touch each other at $r(B_1)=r(B_2)$. Consequently the set $X$ can be decomposed into $X=X_a\sqcup X_r$ such that $a(B)\in X_a$ and $r(B)\in X_r$ for each component $B$ of $\BBB$. Both $X_a$ and $X_r$ are fixed by $F_0$.

Since $G$ is holomorphic in a neighborhood of $\overline{\WWW}$, $F_0$ is holomorphic in a neighborhood of $\overline{\BBB}$. In the following we want to prove that cycles in $X_a$ are attracting and cycles in $X_r$ are repelling.

Let $y\in Y$ be a point with period $p\ge 1$. Pick an attracting flower $\VVV$ of $G$ at the point $y$ such that each component of  $\partial\RRR(r^2)\smm\{y\}$ intersects $\partial\VVV$ at exactly two points and $\tau(\partial\VVV\cap\RRR(r^2))=\partial\VVV\cap\RRR(r^2)$. Then $D=p(\VVV\smm\overline{\WWW(r^2)})$ is a disjoint union of once-punctured disks whose punctures are contained in $X_a$, and $\overline{F_0^k(D)}\subset D$. This shows that each cycle in $X_a$ is attracting. On the other hand, pick a repelling flower $\VVV'$ at $y$ such that each component of $\partial\RRR(r^2)\smm\{y\}$ intersects $\partial\VVV'$ at exactly two points and $\tau(\partial\VVV'\cap\RRR(r^2))=\partial\VVV'\cap\RRR(r^2)$. Then $D'=p(\VVV'\smm\overline{\WWW(r^2)})$ is a once-punctured disk with puncture $x\in X_r$, and $\overline{F_0^{-k}(D')}\subset D'$, where the inverse branch is taken along $F_0^{k}(x)=x$. This shows that each cycle in $X_r$ is repelling.

\vskip 0.24cm

{\bf Step 3. Extension of the inverse map of the projection}. Note that restricted to $\cbar\smm\overline{\WWW}$ the holomorphic projection $p$ is injective and $p(\cbar\smm\overline{\WWW})=\cbar\smm\overline{\BBB}$. We want to extend its inverse map to be a quotient map of $\cbar$ as follows.

Note that $\pi(\WWW)=\DD^*$ and $\pi(\RRR(r))=\AA(r,1)$. Let $w:\, \AA(r,1)\to\DD^*$ be the homeomorphism defined in Proposition \ref{model} (4), i.e.
$$
\begin{cases}
w(z)=z & \text{ if } \sqrt{r}\le |z|<1, \\
\arg w(z)=\arg z  & \text{ if }r<|z|<{\sqrt r}, \\
\log|w(z)|=\dfrac{1}{2}\left(1+\log\dfrac{\log(1/r)}{2\log(|z|/r)}\right)\log{r} & \text{ if }r<|z|<{\sqrt r}.
\end{cases}
$$
Then there exists a unique homeomorphism $\wt w:\, \RRR(r)\to\WWW$ such that $\pi\circ\wt w=w\circ\pi$ and $\wt w=\text{\rm id}$ on $\RRR(\sqrt r)$. It follows that $G\circ\wt w=\wt w\circ G$ on $\RRR(r)$. Define
$$
q=\begin{cases}
p^{-1}: & \cbar\smm\overline{\BBB}\to\cbar\smm\overline{\WWW}, \\
{\wt w}\circ p^{-1}: & \BBB\smm\SSS\to\RRR(r)\to\WWW.
\end{cases}
$$
Then $q: \cbar\smm\SSS\to\cbar\smm Y$ is a homeomorphism. It is quasiconformal on any compact subset of $\cbar\smm\SSS$ and holomorphic in $\cbar\smm\overline{\BBB'}$, where $\BBB'\subset\BBB$ is the interior of $q^{-1}(\overline{\WWW(\sqrt r)})$. It can be extended to a quotient map of $\cbar$ by setting $q(\SSS)=Y$. See the following commutative diagram:
\[
\xymatrix{
\BBB\smm\SSS\ \ar[d]_{\pi\circ p^{-1}}\ar[r]^{q} & \ \WWW\ \ar[d]_{\pi} \\
\AA(r,1)\ \ar[r]^{w} & \ \DD^* }
\]
From $G\circ\wt w=\wt w\circ G$ on $\RRR(r)$, we obtain that
$$
G\circ q=q\circ F_0\quad\text{on }\cbar\smm q^{-1}(\overline{\WWW_1}).
$$

\vskip 0.24cm

{\bf Step 4. Construction of a marked sub-hyperbolic semi-rational map}. Each component $E$ of $\overline{\WWW}$ is a full continuum which contains exactly one point of $\QQQ$. Pick a disk $U(E)\supset E$ such that $U(E)\smm E$ contains no critical values of $G$, $\partial U(E)$ is disjoint from $\QQQ$ and all these domains $U(E)$ have disjoint closures. Denote by $U_0$ their union. Then $U(E)$ contains at most one critical value of $G$. Thus each component of $G^{-1}(U_0)$ is a disk containing exactly one component of $G^{-1}(\overline{\WWW})$.

Let $U_1$ be the union of all components of $G^{-1}(U_0)$ which contain a component of $\overline{\WWW_1}$. Since $\overline{\WWW_1}$ is disjoint from $\PPP'_G$ once $U_0$ is close enough to $\WWW$, we may assume that $U_1\smm\overline{\WWW_1}$ is disjoint from $\QQQ$.

\vskip 0.24cm

Define a branched covering $F$ of $\cbar$ satisfying the following conditions:

(a) $F(z)=F_0(z)$ on $\cbar\smm q^{-1}(U_1)$ and hence $G\circ q=q\circ F$ on $\cbar\smm q^{-1}(U_1)$.

(b) $F: q^{-1}(U_1\smm Y_1)\to q^{-1}(U_0)\smm X_r$ is a covering, where $Y_1=G^{-1}(Y)\smm Y$.

Then $F$ is a geometrically finite branched covering of $\cbar$ with
$$
q^{-1}(\PPP_G\smm Y)\cup X_a\subset\PPP_F\subset q^{-1}(\PPP_G\smm Y)\cup X,
$$
and $\PPP'_F=q^{-1}(\PPP'_G\smm Y)\cup X_a$. In particular, $q(\PPP_F)=\PPP_G$. Set
$$
\PPP=q^{-1}(\QQQ\smm Y)\cup X.
$$
Then $q(\PPP)=\QQQ$, $(\PPP_F\cup F(\PPP))\subset\PPP$ and $\#(\PPP\smm\PPP_F)<\infty$.

Note that $F$ is holomorphic in a neighborhood of $q^{-1}(\PPP'_G)$. Since each cycle in $\PPP'_G\smm Y$ is either attracting or super-attracting for $G$, each cycle in $q^{-1}(\PPP'_G\smm Y)$ is either attracting or super-attracting for $F$. On the other hand, $F_0$ is holomorphic in a neighborhood of $\overline{\BBB}$ and each cycle in $X_a$ is attracting. Thus $(F, \PPP)$ is a marked sub-hyperbolic semi-rational map.

\vskip 0.24cm

{\bf Step 5. Lift of the quotient map}. For each component $D$ of $U_1$, $G: D\to G(D)$ is a covering with at most one critical value. On the other hand, $F: q^{-1}(D)\to F(q^{-1}(D))$ is also a branched covering with at most one critical value, and they have the same degree. Thus the quotient map
$$
q: F(q^{-1}(D))=q^{-1}(G(D))\to G(D)
$$
can be lifted to a quotient map $\wt q: q^{-1}(D)\to D$ that coincides with $q$ on the boundary, i.e. $G\circ\wt q=q\circ F$ on $q^{-1}(D)$.

Define $\wt q=q$ on $\cbar\smm q^{-1}(U_1)$. Then $\wt q$ is a quotient map of $\cbar$ and
$$
G\circ\wt q=q\circ F.
$$
Since each component of $q^{-1}(U_1)$ is a disk containing at most one point of $\PPP$, $\wt q$ is isotopic to $q$ rel $(\cbar\smm q^{-1}(U_1))\cup\PPP$.

\REFLEM{G2F}
If $\UUU$ is a fundamental set of $F$, $q(\UUU)$ contains a fundamental set of $G$.
\ENDLEM

\beginp
We only need to show that $q(\UUU)$ contains an attracting flower at each cycle in $Y$. Pick a disk $D_x$ around each point $x\in X_a$ in $\UUU\smm q^{-1}(U_1)$ such that $\partial D_x\cap\SSS$ is a single point, $F$ is univalent and holomorphic in $D_x$ and $\overline{F(D)}\subset D$, where $D$ denotes the union of the $D_x$. Then $\VVV=q(D)\smm Y$ is a disjoint union of disks, $G$ is holomorphic in $\VVV$ and $\overline{G(\VVV)}\subset\VVV\cup Y$.

Let $\{G^n(w)\}$ $(n\ge 0)$ be an orbit converging to a cycle in $Y$ but with $G^n(w)\notin Y$ for all $n\ge 0$. Let $z_n=q^{-1}(G^n(w))$. Then $F(z_n)=z_{n+1}$ once $n$ is large enough and $\{z_n\}$ converges to an attracting cycle in $X_a$ as $n\to\infty$. Thus $z_n\in D$ and hence $G^n(w)\in q(D)$ once $n$ is large enough. Therefore $\VVV$ contains an attracting flower at each cycle in $Y$. \qed

\REFLEM{no-obs1}
The marked semi-rational map $(F, \PPP)$ has no Thurston obstructions.
\ENDLEM

\beginp Assume by contradiction that $\Gamma$ is an irreducible multicurve of $(F, \PPP)$ with $\lambda(\Gamma)\ge 1$. Recall that $\SSS=q^{-1}(Y)$ is a star. Thus we may further assume that for each $\gamma\in\Gamma$, $\#(\gamma\cap\SSS)$ is minimal in its isotopy class. Since $F: \SSS\to\SSS$ is bijective, $k=\#(\gamma\cap\SSS)<\infty$ is a constant for $\gamma\in\Gamma$.

If $k=0$, then for each $\gamma\in\Gamma$, $q(\gamma)$ is essential and non-peripheral in $\cbar\smm Q$ since $q: \PPP\smm X\to\QQQ\smm Y$ is injective and $q(X)=Y\subset\PPP'_G$. Thus $\Gamma_1=\{q(\gamma):\, \gamma\in\Gamma\}$ is a multicurve of $(G, \QQQ)$.  Noticing that $G\circ\wt q=q\circ F$ and that $\wt q$ is isotopic to $q$ rel $\PPP$, we have $\lambda(\Gamma)=\lambda(\Gamma_1)<1$. This is a contradiction.

Now we assume that $k>0$. Then there exists at most one component of $F^{-1}(\gamma)$ isotopic to a curve in $\Gamma$ rel $\PPP$ for each $\gamma\in\Gamma$ since $F: \SSS\to\SSS$ is bijective. Thus for each $\gamma\in\Gamma$, there is exactly one curve $\beta\in\Gamma$ such that $F^{-1}(\beta)$ has a component isotopic to $\gamma$ rel $\PPP$ since $\Gamma$ is irreducible. Therefore each entry of the transition matrix $M(\Gamma)$ is less than or equal to $1$. Because  $\lambda(\Gamma)\ge 1$, there is a curve $\gamma\in\Gamma$ such that $\gamma$ is isotopic to a component $\delta$ of $F^{-p}(\gamma)$ rel $\PPP$ for some integer $p\ge 1$ and $F^p$ is injective on $\delta$.

Let $\UUU$ be a fundamental set of $F$ that is disjoint from every curve in $\Gamma$. Since $q(\UUU)\smm Y$ contains a fundamental set of $G$ and $q$ is injective on $\cbar\smm\SSS$, $q(\gamma)$ is disjoint from a fundamental set of $G$.

Suppose $\gamma$ intersects at least two components of $\SSS$. Let $\beta$ be a component of $q(\gamma)\smm Y$ such that $\beta$ joins two distinct points of $Y$. Then $\beta$ is isotopic to a component of $G^{-kp}(\beta)$ rel $\QQQ$ for some integer $k>0$ since $\#(\gamma\cap\SSS)$ is minimal in its isotopy class and $\gamma$ is isotopic to a component of $F^{-p}(\gamma)$ rel $\PPP$. Thus $\beta$ is a connecting arc of $(G, \QQQ)$. This is a contradiction.

Suppose that $\gamma$ intersects exactly one component $S$ of $\SSS$. Let $U_1, U_2$ be the two components of $\cbar\smm\gamma$. Since $\gamma$ is non-peripheral, each $U_i$ contains at least two points of $\PPP$,  and either one of them is not contained in $S$ or one of them is contained in $\PPP'_F$ since $\SSS$ contains exactly one isolated point of $\PPP$. In the latter case, $U_i$ contains infinitely many points of $\PPP$. Consequently, each $U_i$ contains at least one point of $\PPP\smm S$. Thus $U_i\smm S$ has a component $U'_i$ which contains at least one point of $\PPP\smm S$.

Let $y=q(S)$. Then there exists a component $\beta$ of $q(\gamma)\smm\{y\}$ such that $\beta\cup\{y\}$ separates $q(U'_1)$ from $q(U'_2)$. In other words, each component of $\cbar\smm(\beta\cup\{y\})$ contains at least one point of $\QQQ$ since $q(\PPP)=\QQQ$. As above, $\beta$ is isotopic to a component of $G^{-kp}(\beta)$ rel $\QQQ$ for some integer $k>0$. Thus $\beta$ is a connecting arc of $(G, \QQQ)$. This is a contradiction. Thus $(F, \PPP)$ has no Thurston obstructions.
\qed

\subsection{Proof of Theorems \ref{existence} and \ref{s-plumbing}}

{\noindent\it Proof of Theorem \ref{s-plumbing-1}}.
Let $(G, \QQQ)$ be a marked semi-rational map with parabolic cycles in $\PPP'_G$ and $\#(\QQQ\smm\PPP_G)<\infty$. Suppose that $(G, \QQQ)$ has neither Thurston obstructions nor connecting arcs. Let $(F, \PPP)$ be the marked sub-hyperbolic semi-rational map constructed in step 4. By Lemma \ref{no-obs1} and Theorem \ref{hyperbolic}, there exist a marked rational map $(f, \PPP_1)$ and a c-equivalence $(\phi_0,\phi_1)$ from $(F, \PPP)$ to $(f, \PPP_1)$ on a fundamental set $\UUU$ of $F$.

Recall that $F(\BBB)=\BBB$ and $F$ is conformal in a neighborhood of $\overline{\BBB}$. Pick a disk $D_x\subset\UUU$ around each point $x\in X_a$ such that $\partial D_x\cap\BBB$ is an arc and $\overline{F(D)}\subset D$, where $D$ denotes the union of the $D_x$. Then $\sA=\pi_f\circ\phi_0(D\cap\BBB)$ is a multi-annulus in $\sR_f$. From step 1 and condition (b) in step 4, we know that $\sA$ is starlike.

\vskip 0.24cm

\begin{proposition}\label{coincide}
For each component $B$ of $\BBB$, let $B'$ be the component of $\pi_f^{-1}(\sA)$ such that $a(B')=\phi_0(a(B))$. Then $r(B')=\phi_0(r(B))$.
\end{proposition}

\beginp
Let $p\ge 1$ be the period of $B$. Let $\wt B=\phi_0(B)$. Then $\wt B$ coincides with $B'$ in a neighborhood of $a(B')$. Since $\deg_{a(B')}f^{p}=1$, there is a unique component $\wt B_k$ of $f^{-kp}(\wt B)$ whose closure contains the point $a(B')$ for $k\ge 1$. Consequently, $\wt B_k$ coincides with $B'$ in a neighborhood of $a(B')$.

Note that $\phi_0(r(B))\in\PPP_1$ is a repelling periodic point of $f$. Cut $\wt B$ into three disks $\wt B=U\cup V\cup W$ such that $U\subset B'$, $\overline{V}\cap\PPP_1=\emptyset$ and $W$ is contained in a linearizable disk of the repelling periodic point $\phi_0(r(B))$. Then each $\wt B_k$ is also cut into three domains $\wt B_k=U_k\cup V_k\cup W_k$ such that $f^{kp}(U_k)=U$, $f^{kp}(V_k)=V$ and $f^{kp}(W_k)=W$. Since $f^p(B')=B'$, we have $U_k\subset B'$ for $k\ge 1$. On the other hand, $\text{\rm diam}_s V_k\to 0$ as $k\to\infty$ by Lemma \ref{shrinking} and $\text{\rm diam}_s W_k\to 0$ as $k\to\infty$ since $W$ is contained in a linearizable disk. Thus $r(B')=\phi_0(r(B))$.
\qed

\vskip 0.24cm

By this proposition, we may assume by modification that $\phi_0$ is holomorphic in $\BBB\cup\UUU$ and $\phi_1$ is isotopic to $\phi_0$ rel $\BBB\cup\UUU\cup\PPP$. Then $\phi_0(\BBB)$ is the union of periodic bands of $\sA$, $\phi_0(\SSS)$ is the union of periodic skeletons of $\sA$ and $\phi_1(F^{-1}(\SSS))$ is the union of skeletons of level at most $1$ by the equation $\phi_0\circ F=f\circ\phi_1$.

Let $f_t=\phi_t\circ f\circ\phi_t^{-1}$ $(t\ge 0)$ be the simple pinching path of $f$ supported on $\sA$. By Theorem \ref{s-pinching}, $\{f_t\}$ converges uniformly to a rational map $g$ and $\{\phi_t\}$ converges uniformly to a quotient map $\varphi$ of $\cbar$ as $t\to\infty$, and $g\circ\varphi=\varphi\circ f$. Let $\QQQ_1=\varphi(\PPP_1)$. Then $(g, \QQQ_1)$ is a marked rational map.

For any point $w\in\cbar$, $q^{-1}(w)$ is either a single point or a component of $\SSS$, and $\wt q^{-1}(w)$ is either a single point or a component of $F^{-1}(\SSS)$. Let
$$
\zeta_0=\varphi\circ\phi_0\circ q^{-1}\quad\text{and}\quad \zeta_1=\varphi\circ\phi_1\circ\wt q^{-1}.
$$
These are well-defined quotient maps of $\cbar$, and $\zeta_0\circ G=g\circ\zeta_1$.

From the definitions of $\wt q$ and $q$, there exists a fundamental set $\UUU_1\subset\UUU$ of $F$ such that $\wt q\circ q^{-1}$ is isotopic to the identity rel $q(\overline{\UUU_1}\cup\overline{\BBB}\cup\PPP)$. By Lemma \ref{G2F}, there exists a fundamental set $\UUU_G$ of $G$ such that $\UUU_G\subset q(\UUU_1)$. Note that $q^{-1}(\overline{\UUU_G}\cup\QQQ)\subset\overline{\UUU_1}\cup\overline{\BBB}\cup\PPP$. Thus $\zeta_1$ is isotopic to
$$
\zeta_1\circ(\wt q\circ q^{-1})=\varphi\circ\phi_1\circ q^{-1}
$$
rel $\overline{\UUU_G}\cup\QQQ$. Since $\phi_1$ is isotopic to $\phi_0$ rel $\overline{\UUU_1}\cup\overline{\BBB}\cup\PPP$, $\zeta_0$ is isotopic to $\varphi\circ\phi_1\circ q^{-1}$ rel $\overline{\UUU_G}\cup\QQQ$. Therefore $\zeta_1$ is isotopic to $\zeta_0$ rel $\overline{\UUU_G}\cup\QQQ$.

By Lemma \ref{equivalent}, $\phi_0(\UUU_1)$ is a fundamental set of $f$. Recall that $\phi_0(\BBB)$ is the union of all periodic bands of $\sA$. Let $\BBB'\subset\BBB$ be the domain defined in step 3. Then $\phi_0(\BBB')$ is the union of all periodic bands of $\sA'\subset\sA$, and $q(\overline{\BBB'})=\overline{\WWW(\sqrt r)}$. Once $\UUU_1$ is small enough, we may assume that
$$
\phi_0(\UUU_1\smm\overline{\BBB'})\subset\FFF_f\smm\overline{\pi_f^{-1}(\sA')}.
$$
Then $\varphi$ is holomorphic in $\phi_0(\UUU_1\smm \overline{\BBB'})$. Thus $\zeta_0$ is holomorphic in $\UUU(G)\smm\overline{\WWW(\sqrt r)}$.

By the definition of $q$ in step 3, $q$ is quasiconformal on any domain compactly contained in $\cbar\smm\SSS$ and holomorphic in $\cbar\smm\overline{\BBB}$. Moreover, there exist holomorphic universal coverings $\pi_1: \BBB\smm\SSS\to \AA(r,1)$ and $\pi_2: \WWW=q(\BBB\smm\SSS)\to\DD^*$ such that

$\bullet$ $\pi_1(z_1)=\pi_1(z_2)$ if and only if $z_1=F^i(z_2)$ for some integer $i\in\ZZ$,

$\bullet$ $\pi_2(w_1)=\pi_2(w_2)$ if and only if $w_1=G^j(w_2)$ for some integer $j\in\ZZ$, and

$\bullet$ $|\pi_1(z)|\to 1$ as $z\to\partial\BBB\smm\SSS$ and the following diagram commutes,
\[
\xymatrix{
\BBB\smm\SSS\ \ar[d]_{\pi_1}\ar[r]^{q} & \ \WWW\ \ar[d]_{\pi_2} \\
\AA(r,1)\ \ar[r]^{w} & \ \DD^* }
\]
where $w$ is the map defined in Proposition \ref{model} (4). Note that the map $w$ commutes with any rotation of $\DD^*$. Since $\phi_0(\BBB)$ is the union of all periodic bands of $\sA$, by Lemma \ref{model-varphi}, $\zeta_0$ is holomorphic in $\WWW$ and hence in $\WWW\cup(\UUU_G\smm\overline{\WWW(\sqrt r)})=\UUU_G\smm Y=\UUU_G$. The last equality is because $\xi_0(\UUU_G)$ is a fundamental set of $g$ and hence is contained in the Fatou set of $g$.

Obviously, $\zeta_0:\, \QQQ\to\QQQ_1$ is bijective. Thus there exists a homeomorphism $\psi_0$ of $\cbar$ such that $\psi_0$ is isotopic to $\zeta_0$ rel $\overline{\UUU_G}\cup\QQQ$. Let $\psi_1$ be the lift of $\psi_0$. Then it a homeomorphism of $\cbar$ isotopic to $\psi_0$ rel $\overline{\UUU_G}\cup\QQQ$ by Theorem 1.12 in \cite{BD}, and
$$
g\circ\psi_1=\psi_0\circ G.
$$
Therefore $(G, \QQQ)$ is c-equivalent to $(g, \QQQ_1)$. See the following diagram:
\[\xymatrix{
\ \cbar\ \ar[d]_{g}&\ \cbar\ \ar[l]_{\varphi}\ar[d]_{f}& \ar[l]_{\phi_1}\ \cbar\ \ar[r]^{\wt q}\ar[d]_{F}& \ar@/_2pc/[lll]_{\zeta_1}\ \cbar\ \ar[d]^{G}\\
\ \cbar\ & \ \cbar\ \ar[l]_{\varphi} & \ar[l]_{\phi_0}\ \cbar\ \ar[r]^{q} & \ar@/^2pc/[lll]_{\zeta_0}\ \cbar\
}
\]
\qed

{\noindent\it Proof of Theorems \ref{existence-1} and \ref{existence}}.
Let $(G, \QQQ)$ be a semi-rational map with $\PPP'_G\neq\emptyset$ and $\#(\QQQ\smm\PPP_G)<\infty$. If $(G, \QQQ)$ is c-equivalent to a marked rational map $(g, \QQQ_1)$, then $(g, \QQQ_1)$ has neither Thurston obstructions nor connecting arcs, by Theorems \ref{McMullen} and \ref{connecting}. Neither does $(G, \QQQ)$. Conversely, if $(G, \QQQ)$ has neither Thurston obstructions nor connecting arcs, then $(G, \QQQ)$ is c-equivalent to a marked rational map, by Theorems \ref{s-plumbing-1} and \ref{hyperbolic}.
\qed

\vskip 0.24cm

{\noindent\it Proof of Theorem \ref{s-plumbing}}. This is a direct consequence of Theorem \ref{s-plumbing-1}.
\qed

\section{Distortion of univalent maps}

\subsection{Modulus difference distortion}

Let $V\subset\cbar$ be an open set and let $\phi: V\to\cbar$ be a univalent map. Define
$$
\sD_0(\phi, V)=\sup_{E_1, E_2\subset V}|\md A(E_1, E_2)-\md A(\phi(E_1), \phi(E_2))|,
$$
where $E_1, E_2$ are disjoint full continua in $V$ and $A(E_1, E_2):=\cbar\smm(E_1\cup E_2)$. Define
$$
\sD_1(\phi)(z, w)=\left|\log\frac{|\phi'(z)\phi'(w)||z-w|^2}{|\phi(z)-\phi(w)|^2}\right|
$$
for $(z,w)\in V\times V,\, z\neq w$, and define
$$
\sD_1(\phi, V)=\|\sD_1(\phi)(z,w)\|_{\infty}.
$$
Obviously, for $i=0,1$,
$$
\sD_i(\phi^{-1}, \phi(V))=\sD_i(\phi, V),
$$
and for any M\"{o}bius transformations $\beta$ and $\gamma$ of $\cbar$,
$$
\sD_i(\gamma\circ \phi\circ\beta, \beta^{-1}(V))=\sD_i(\phi, V).
$$

\REFTHM{d2d}
Suppose that $\sD_0(\phi, V)=\delta<\infty$. Then

(a) $\sD_1(\phi, V)\le 2\pi\delta$.

(b) Assume that $V$ contains $0, \infty$ and $\DD(1,r_0)$ for some $r_0>0$. If $\phi$ fixes $0,1$ and $\infty$, then there exists a constant $C(r_0)>0$ depending only on $r_0$ such that
$$
\text{\rm dist}_s(\phi(z), z)\le C(r_0)\delta.
$$
\ENDTHM

\beginp
(a) We want to prove that $\sD_1(\phi)(z,w)\le 2\pi\delta$ for $z,w\in V$ and $z\neq w$. Since $\sD_1(\phi, V)$ is invariant under M\"{o}bius transformations, we may assume that $z=0, w=1$, $\phi(0)=0$ and $\phi(1)=1$. Then
$\sD_1(\phi)(z,w)=|\log|\phi'(0)\phi'(1)||$.

Let $M_0(r), m_0(r)$ be the supremum and infimum of $|\phi(z)|$ on the disk $\DD(r)$, and let $M_1(r), m_1(r)$ be the supremum and infimum of $|\phi(z)-1|$ on the disk $\DD(1,r)$. Then
$$
A(\DD(M_0(r)), \DD(1,M_1(r)))\subset A(\phi(\DD(r)),\phi(\DD(1,r)))\subset A(\DD(m_0(r)), \DD(1,m_1(r))).
$$
By a direct computation, we have
$$
\md A(\DD(r_1),\DD(1,r_2))=\frac 1{2\pi}\log (\kappa(r_1, r_2)+\sqrt{\kappa(r_1, r_2)^2-1}),
$$
where $\kappa(r_1, r_2)=(1-r_1^2-r_2^2)/(2r_1r_2)$. Since both $M_0(r)/r$ and $m_0(r)/r$ converge to $|\phi'(0)|$, and both $M_1(r)/r$ and $m_1(r)/r$ converge to $|\phi'(1)|$ as $r\to 0$, we get
$$
|\md A(\DD(r), \DD(1,r))-\md A(\phi(\DD(r)),\phi(\DD(1,r)))|=\dfrac{|\log|\phi'(0)\phi'(1)||}{2\pi}+O(r).
$$
Letting $r\to 0$, we deduce that $|\log|\phi'(0)\phi'(1)||\le 2\pi\delta$.

(b) Set $\lambda_0=|\phi'(0)|$, $\lambda_1=|\phi'(1)|$ and $\lambda_2= \displaystyle\lim_{z\to\infty}|z/\phi(z)|$. Then $|\log\lambda_i\lambda_j|\le 2\pi\delta$ for $i\neq j$. It follows that $|\log\lambda_i|\le 3\pi\delta$ for $i=0,1,2$.

For any point $z\in V\smm\{\infty\}$, we have $\sD_1(\phi)(z,w)\le 2\pi\delta$ for $w\in V\smm\{z\}$. Letting $w\to\infty$, we get $|\log|\lambda_2\phi'(z)||\le 2\pi\delta$. Thus
$$
|\log|\phi'(z)||\le 5\pi\delta\quad\text{for }z\in V\smm\{\infty\}.\eqno{(8.1)}
$$
Therefore
$$
\left|\log\frac{|\phi(z_1)-\phi(z_2)|^2}{|z_1-z_2|^2}\right|\le 10\pi\delta\eqno{(8.2)}
$$
for any points $z_1, z_2\in V\smm\{\infty\}$ with $z_1\neq z_2$. Applying (8.2) for the pairs $(z,0)$ and $(z,1)$ with $|z|\le 2$, we obtain
$$
||\phi(z)|^2-|z|^2|\le 4(e^{10\pi\delta}-1),\text{ and}
$$
$$
||\phi(z)-1|^2-|z-1|^2|\le 9(e^{10\pi\delta}-1).
$$
It follows that
$$
|\,\text{Re}\, (\phi(z)-z)|\le 7(e^{10\pi\delta}-1), \text{ and} \eqno{(8.3)}
$$
$$
|\,\text{Im}\,(\phi(z)-z)|\cdot |\text{Im}\,(\phi(z)+z)|\le 18(e^{20\pi\delta}-1).\eqno{(8.4)}
$$

Let $s=\max\{1, r_0/7\}$. Without loss of generality, we assume that $e^{10\pi\delta}-1<s/21$. Then
$$
|\phi'(w)|\le e^{5\pi\delta}\le\sqrt{1+s/21}\le 2\quad\text{for }w\in V\smm\{\infty\}.
$$
For each point $w\in\DD(1,s)$, we have $\DD(w, 6r_0/7)\subset\DD(1, r_0)$. Applying Cauchy's integral formula on $\partial\DD(w, 6r_0/7)$, we get
$|\phi''(w)|\le 1/(3s)$. Thus
$$
|\,\text{Re}\,(\phi'(w)-\phi'(1))|\le |\phi'(w)-\phi'(1)|\le\left|\int_1^w\phi''(\zeta) d\zeta\right|\le 1/3.\eqno{(8.5)}
$$
Applying (8.3) for $z=1+s$, we get
$$
\left|\int_1^{1+s}\text{Re}\,(\phi'(\zeta)-1) d\zeta\right|=
|\,\text{Re}\,\phi(1+s)-(1+s)|\le 7(e^{10\pi\delta}-1)<\frac{s}{3}.
$$
By (8.5), we have
\begin{align*}
|\,\text{Re}\,\phi'(1)-1| & = \frac{1}{s}\left|\int_1^{1+s}\text{Re}\,(\phi'(1)-1)d\zeta\right| \\
& \le\frac{1}{s}\left|\int_1^{1+s}\text{Re}\,(\phi'(\zeta)-1)d\zeta\right|+\frac{1}{s}\left|\int_1^{1+s}\text{Re}\,(\phi'(\zeta)-\phi'(1)) d\zeta\right| \\
& \le
\frac{2}{3}.
\end{align*}
Thus
$$
|\,\text{Re}\,\phi'(w)-1|\le |\,\text{Re}\,\phi'(1)-1|+|\,\text{Re}\,(\phi'(w)-\phi'(1))|\le 1
$$
for $w\in\DD(1,s)$. It follows that $\text{Re}\,\phi'(w)\ge 0$ for $w\in\DD(1,s)$. Therefore
$$
\text{Im}\, \phi(1+is)=\int_{0}^s \text{Re}\,\phi'(\zeta)d|\zeta|\ge 0.
$$
Applying (8.4) for $z=1+is$, we obtain
$$
|\text{Im}\, (\phi(1+is)-s)|\le 18(e^{20\pi\delta}-1)/s.\eqno{(8.6)}
$$
For each point $z\in\DD\cap V$, let $w_1=z-(1+is)$ and $w_2=\phi(z)-\phi(1+is)$. Applying (8.2) and (8.3) for the points $z$ and $1+is$, we get
$$
||w_1^2|-|w_2^2||\le 9(e^{10\pi\delta}-1)\quad\text{and}
$$
$$
|\,\text{Re}\, (w_1-w_2)|\le 14(e^{10\pi\delta}-1).
$$
It follows that
$$
|\,\text{Im}\,(w_2-w_1)|\cdot |\text{Im}\,(w_2+w_1)|\le 51(e^{20\pi\delta}-1).\eqno{(8.7)}
$$
Since
$$
\text{Im}\,(\phi(z)+z)-\text{Im}\,(w_1+w_2)=\text{Im}\, \phi(i+is)+s\ge s,
$$
we have either $|\,\text{Im}\,(w_1+w_2)|\ge s/2$ or $|\,\text{Im}\,(\phi(z)+z)|\ge s/2$. In the former case,
$$
|\,\text{Im}\,(w_2-w_1)|\le 102(e^{20\pi\delta}-1)/s
$$
by (8.7). Combining with (8.6), we obtain
$$
|\,\text{Im}\,(\phi(z)-z)|\le |\,\text{Im}\,(w_2-w_1)|+|\text{Im}\, (\phi(1+is)-s)|\le  120(e^{20\pi\delta}-1)/s.
$$
Applying (8.4) in the latter case,
$$
|\,\text{Im}\,(\phi(z)-z)|\le 36(e^{20\pi\delta}-1)/s.
$$
In summary, we have
$$
|\,\text{Im}\,(\phi(z)-z)|\le 120(e^{20\pi\delta}-1)/s.
$$
Combining with (8.3), we get
$$
|\,\phi(z)-z|\le 127(e^{20\pi\delta}-1)/s\quad\text{for }z\in\DD\cap V.
$$
By considering $1/\phi(1/z)$ and comparing spherical distance with Euclidean distance, we can get the constant $C(r_0)$.
\qed

\REFCOR{identity}
Let $V\subset\cbar$ be an open set and let $a_i$ ($i=0,1,2$) be three distinct points in $V$. Let $\{\phi_n\}$ be a sequence of univalent maps from $V$ into $\cbar$ such that $\phi_n(a_i)=a_i$ for $i=0,1,2$. If $\sD_0(\phi_n, V)\to 0$ as $n\to\infty$, then $\{\phi_n\}$ converges uniformly to the identity on $V$ as $n\to\infty$.
\ENDCOR

{\noindent\bf Remark}. Denote by $S_{\phi}$ the Schwarzian derivative of $\phi$ and by $\lambda(z)|dz|$ the Poincar\'{e} metric on $V$ if $\cbar\smm V$ contains at least three points. It is proved in \cite{CT2} that $\|S_{\phi}(z)\lambda^{-2}(z)\|_{\infty}\le C\delta(|\log\delta|^2+1)$ for a universal constant $C>0$, where $\delta=\sD_0(\phi, V)$.

\vskip 0.24cm

In order to estimate $\sD_0(\phi,V)$ we need the following quantity.

\REFLEM{area}
Let $W\subset\cbar$ be an open set with $\#(\cbar\smm W)\ge 2$. Suppose that $E$ is a measurable set with $\overline{E}\subset W$. Then
$$
\area_p(E,W):=\sup_{h}\area(\rho_*,h(E))<\infty,
$$
where the supremum is taken over all univalent maps $h: W\to\CC^*$ and $\rho_*(z)=1/|z|$ is the density of a planar metric on $\CC^*$.
\ENDLEM

\beginp
We may assume that $W\subset\CC$. Denote by $d>0$ the Euclidean distance between $E$ and $\partial W$, and by $M$ the Euclidean area of $E$. Let $h: W\to\CC^*$ be a univalent map. Applying Koebe's $1/4$-Theorem for the disk $\DD(z, d)$ with $z\in E$, we have $|h'(z)/h(z)|\le 4/d$. Thus
$$
\area(\rho_*,h(E))=\iint_{E}
\left|\dfrac{h'(z)}{h(z)}\right|^2 dxdy \le\frac{16M}{d^2}.
$$
\qed

\subsection{Nested disk systems}

Let $X\subset\cbar$ be a finite set. For each point $x\in X$, let $D_x\subset\cbar$ be a disk with $x\in D_x$. We will call $\{D_x\}_{x\in X}$ a {\bf nested disk system} if $D_x\neq D_y$ for distinct points $x,y\in X$, and $D_x\cap D_y\neq\emptyset$ for $x,y\in X$ implies that either $D_x\subset D_y$ or $D_y\subset D_x$.

For any $r\in (0,1)$, let $D_x(r)=\chi^{-1}(\DD(r))$, where $\chi$ is a conformal map from $D_x$ to $\DD$ such that $\chi(x)=0$. Let $s(r): (0, 1]\to (0, 1]$ be a (non-strictly) increasing function with $s(r)\ge r$ and $s(r)\to 0$ as $r\to 0$. A nested disk system $\{D_x\}_{x\in X}$ is called {\bf $s(r)$-nested} if for any two disks with $D_y\subset D_x$, $D_y\cap D_x(r)\neq\emptyset$ for some $r\in (0,1]$, we have $D_y\subset D_x(s(r))$.

Let $\{D_x\}_{x\in X}$ be a nested disk system. Let $W\subset\cbar$ be an open set with $\cbar\smm\overline{W}\neq\emptyset$ such that $\cup_{x\in X}D_x\Subset W$. For each point $x\in X$, denote by $V_x$ the union of all domains $D_y$ with $y\neq x$ and $D_y\subset D_x$, and by $W_x$ the component of $W$ containing the point $x$. Then $V_x\subset D_x\Subset W_x$. Let $\lambda\in (0,1)$ be a constant. We will call that the nested disk system $\{D_x\}_{x\in X}$ is {\bf $\lambda$-scattered} in $W$ if for any point $x\in X$ and any univalent map $h: W_x\to\CC^*$,
$$
\area(\rho_*,h(V_x))\le\lambda\cdot\area(\rho_*,h(D_x)).
$$

Let $\{D_x\}_{x\in X}$ be an $s(r)$-nested disk system $\lambda$-scattered in $W$. Obviously, for any subset $X_0\subset X$, $\{D_x\}_{x\in X_0}$ is also an $s(r)$-nested disk system $\lambda$-scattered in $W$. For any univalent map $\phi: W\to\cbar$, $\{\phi(D_x)\}_{x\in X}$ is an $s(r)$-nested disk system $\lambda$-scattered in $\phi(W)$.

\REFTHM{distortion}
Suppose that $\{D_x\}_{x\in X}$ is an $s(r)$-nested disk system $\lambda$-scattered in $W$. Let $D=\cup_{x\in X}D_x$ and $D(r)=\cup_{x\in X}D_x(r)$. Then there exist a constant $r_0\in (0,1)$ and an increasing function $C(r)$ on $(0,r_0)$ with $C(r)\to 0$ as $r\to 0$, which depend only on $\lambda$ and $s(r)$, such that for any $r\in (0, r_0)$ and any univalent map $\phi:\cbar\smm\overline{D(r)}\to\cbar$,
$$
\sD_0(\phi, \cbar\smm\overline{W})\le C(r)\cdot\area_p(D,W),
$$
and for any annulus $A\subset\cbar$ with $\partial A\cap W=\emptyset$ and $\md A<\infty$,
$$
\left|\frac{\md A'}{\md A}-1\right|\le C(r),
$$
where $A'$ is the annulus bounded by $\phi(\partial A)$.
\ENDTHM

\beginp
For each point $x\in X$, let $k(x)=\#\{y\in X:\, D_y\supset D_x\}$. Then $k(x)\ge 1$ and $D_x$ is disjoint from $D_y$ if $k(x)=k(y)$. Denote by $I_k$ and $I_k(r)$ the union of $D_x$ and $D_x(r)$, respectively, for all points $x\in X$ with $k(x)=k$. Let $n=\max\{k(x):\, x\in X\}$. Then
$$
D=I_1\supset I_2\supset\cdots\supset I_{n-1}\supset I_n.
$$

Let $r_1\in (0, 1)$ be a constant with $s(r_1)<\min\{(1-\sqrt{\lambda})^3, 1/64\}$. Given any $r\in (0,r_1)$, let $X'_1=\{x\in X:\, k(x)=k\}$ and let
$$
X'_k=\{x\in X:\, k(x)=k\text{ and }D_x\cap I_j(r)=\emptyset\text{ for all } j<k\}
$$
for $2\le k\le n$. Let $I'_k$ and $I'_k(r)$ be the union of $D_x$ and $D_x(r)$, respectively, over $x\in X'_k$. Then $I'_1=I_1$ and $I'_1(r)=I_1(r)$.

Let $A\subset\cbar$ be an annulus with $W\subset A$ and $\md A<\infty$. Let $\chi_A$ be a conformal map from $A$ to a round annulus in $\CC$ whose core curve is the unit circle. Then $\rho_0(z)=|(\log\chi_A)'(z)|/{2\pi}$ is an extremal metric on $A$, $\width(\rho_0, A)=1$ and
$$
\height(\rho_0, A)=\area(\rho_0, A)=\md A .
$$
Define $\rho_k(z)$ for $1\le k\le n$ inductively by
$$
\rho_k(z)=\begin{cases}
\rho_{k-1}(z)  & \text{ on } A\smm I'_k, \\
\dfrac{\rho_0(z)}{(1-\sqrt{s(r)})^{k}} & \text{ on } I'_k\smm I'_k(s(r)), \\
0 & \text{ on } I'_k(s(r)).
\end{cases}
$$

{\bf Claim 1}. We have $\rho_k(z)=0$ on $I_k(s(r))$ for $k\ge 1$.

{\it Proof}. Claim 1 is true for $k=1$ since $I'_1=I_1$. For $k\ge 2$, we assume by induction that $\rho_j(z)=0$ on $I_j(s(r))$ for $j<k$. Let $x\in X$ be a point with $k(x)=k$. By definition, $\rho_k(z)=0$ on $D_x(s(r))$ if $x\in X'_k$. Now we assume that $x\not\in X'_k$, i.e. there is an integer $1\le j<k$ such that $D_x\cap I_j(r)\neq\emptyset$. Let
$$
D_x\subset D_{x_{k-1}}\subset\cdots\subset D_{x_{j+1}}
$$
be the unique sequence with $k(x_{i})=i$ ($j<i<k$). Then $D_{x_i}\subset I_i-I'_i$. Thus
$$
\rho_k(z)=\rho_{k-1}(z)=\cdots=\rho_j(z)
$$
on $D_x$ by the definition of $\rho_i$. Since $D_x\cap I_j(r)\neq\emptyset$, we have $D_x\subset I_j(s(r))$ and hence $\rho_k(z)=\rho_j(z)=0$ on $D_x(s(r))\subset D_x$.
\qed

{\bf Claim 2}. We have $\rho_k(z)=0$ on $\cup_{i=1}^k I_i(r)$ for $k\ge 1$.

{\it Proof}. By definition, $\rho_1(z)=0$ on $I_1(s(r))$. Thus $\rho_1(z)=0$ on $I_1(r)$ since $r\le s(r)$. Assume that $\rho_{k-1}(z)=0$ on $\cup_{i=1}^{k-1}I_i(r)$ for $k\ge 2$. Then $\rho_k(z)=\rho_{k-1}(z)=0$ on $\cup_{i=1}^{k-1}I_i(r)$ since $\cup_{i=1}^{k-1}I_i(r)$ is disjoint from $I'_k$ by the definition of $X'_k$. Now Claim 2 follows from Claim 1. \qed

{\bf Claim 3}.
$$
\rho_k(z)
\begin{cases}
\le\frac{\rho_0(z)}{(1-\sqrt{s(r)})^k} & \text{on }I_k, \\
\le\frac{\rho_0(z)}{(1-\sqrt{s(r)})^{j}} & \text{on }I_j\smm I_{j+1},\, 1\le j<k,\text{ and} \\
=\rho_0(z) & \text{on } A\smm I_1.
\end{cases}
$$

{\it Proof}. By definition, $\rho_{j+1}(z)=\rho_j(z)$ on $A\smm I'_{j+1}$ and
$$
A\smm I'_{j+1}\supset A\smm I_{j+1}=(A\smm I_1)\cup(I_1\smm I_2)\cup\cdots\cup(I_j\smm I_{j+1}).
$$
Thus $\rho_k(z)=\rho_0(z)$ on $A\smm I_1$ and $\rho_k(z)=\rho_j(z)$ on $I_j\smm I_{j+1}$. Now the claim is proved by induction.

{\bf Claim 4}. For each $1\le k\le n$,
$$
\height(\rho_k, A)\ge\height(\rho_{k-1}, A)\quad\text{and}\quad\width(\rho_k, A)\ge\width(\rho_{k-1}, A).
$$

{\it Proof}. For each point $x\in X'_k$, the map $\log\circ\chi_A$ has a univalent branch $h_x$ on $D_x$ and
$$
\md h_x(D_x\smm\overline{D_x(s(r))})=\frac{-\log s(r)}{2\pi}.
$$
Since $s(r)\le s(r_1)<1/64$, by Lemma 1.1, there exists an annulus $B_x\subset D_x\smm D_x(s(r))$ which separates $D_x(s(r))$ from $\partial D_x$, such that $h_x(B_x)$ is bounded by concentric Euclidean circles and
$$
\md h_x(B_x)\ge\frac{-\log s(r)}{4\pi}.
$$
Denote by $C_x\subset D_x$ the disk bounded by the outer boundary of $B_x$.

Let $\gamma\subset A$ be a locally rectifiable simple closed curve or arc which separates (joins) the two components of $\cbar\smm A$. Assume that $\beta$ is a component of $\gamma\cap C_x$ such that $\beta\cap D_x(s(r))\neq\emptyset$. Then there is an arc $\beta'$ in $C_x$ which joins the two endpoints of $\beta$ such that $h_x(\beta')$ is a straight line segment, since $h_x(C_x)$ is a Euclidean disk.

Denote by $d_1>d_2$ the diameters of the two circles of $\partial h_x(B_x)$. Then $d_2\le d_1\sqrt{s(r)}$. By the definition of $\rho_k(z)$ and Claim 3,
$$
L(\rho_{k-1}, \beta')\le (1-\sqrt{s(r)})^{-k+1}d_1 \le
(1-\sqrt{s(r)})^{-k}(d_1-d_2)\le L(\rho_k, \beta).
$$

Modifying every component of $\gamma\cap C_x$ as above and making the modification for each point $x\in X'_k$, we get a simple closed curve (or an arc) $\gamma'$ which separates (joins) the two components of $\cbar\smm A$. Since $\rho_k(z)\ge\rho_{k-1}(z)$ on $A-I'_k(s(r))$, we have $L(\rho_k, \gamma)\ge L(\rho_{k-1}, \gamma')$. The claim is proved. \qed

\vskip 0.24cm

By Claim 4,
$$
\height(\rho_n, A)\ge\height(\rho_0, A)\quad\text{and}\quad\width(\rho_n, A)\ge\width(\rho_0, A)=1.
$$
Thus $\area(\rho_n, A)\ge\area(\rho_0,A)$ by Lemma \ref{mod2}. We know $\{D_x\}_{x\in X}$ is $\lambda$-scattered in $W$. Since $\chi_A: W\subset A\to\AA(R)\subset\CC^*$ is univalent and the metric $\rho_0$ is the pullback of the metric $\rho_*$ by $\chi_A$, we deduce that for $x\in X$ with $k(x)=k<n$,
$$
\area(\rho_0, I_{k+1}\cap D_x)\le\lambda\cdot\area(\rho_0, D_x).
$$
It follows that $\area(\rho_0,I_{k+1})\le\lambda\cdot\area(\rho_0,I_k)$ and
$$
\area(\rho_0,I_k)\le\lambda^{k-1}\area(\rho_0, I_1)\le\lambda^{k-1}\area_p(D,W).
$$
Therefore
\begin{align*}
0 & \le\area(\rho_n, A)-\area(\rho_0,A) \\
  & \le\ds\sum_{k=1}^{n-1}\iint_{I_k\smm I_{k+1}}(\rho_n^2-\rho_0^2)dxdy+\iint_{I_n}(\rho_n^2-\rho_0^2)dxdy   \\
  & \le\ds\sum_{k=1}^n\left(\frac{1}{(1-\sqrt{s(r)})^{2k}}-1\right)\area(\rho_0,I_k) \\
  & \le\ds\sum_{k=1}^{\infty}\left(\frac{1}{(1-\sqrt{s(r)})^{2k}}-1\right)\lambda^{k-1}\area(\rho_0,I_1).
\end{align*}
Since $s(r)$ is increasing and $s(r)\to 0$ as $r\to 0$, there exists a constant $r_0\in (0,r_1)$ such that $(1-\sqrt{s(r_0)})^2>\lambda$. Thus
\begin{align*}
0 & \le\area(\rho_n, A)-\area(\rho_0,A) \\
  & \le\frac{2\sqrt{s(r)}-s(r)}{(1-\lambda)[(1-\sqrt{s(r)})^2-\lambda]}\area(\rho_0,I_1):=C(r)\cdot\area(\rho_0,I_1).
\end{align*}

Let $A'\subset\cbar$ be the annulus bounded by $\phi(\partial A)$. Let $\wt\rho(w)|dw|$ be the pullback of the metric $\rho_n(z)|dz|$ under the map $\phi$ on $A\smm \overline{D(r)}$ and zero otherwise, i.e. $\wt\rho(\phi(z))=\rho_n(z)/|\phi'(z)|$ for $z\in A\smm\overline{D(r)}$. Then $\area(\wt\rho,A')=\area(\rho_n,A)$,
$$
\height(\wt\rho, A')=\height(\rho_n, A)\quad\text{and}\quad\width(\wt\rho, A')=\width(\rho_n, A).
$$
By Lemma \ref{mod2},
$$
|\md A'-\md A |\le |\area(\rho_n, A)-\area(\rho_0,A)|\le C(r)\cdot\area(\rho_0,I_1).\eqno{(8.8)}
$$

Because $\area(\rho_0,I_1)\le\area_p(D,W)$, we obtain
$$
\sD_0(\phi, \cbar\smm\overline{W})\le C(r)\cdot\area_p(D,W).
$$

Now assume that $A\subset\cbar$ is an annulus with $\partial A\cap W=\emptyset$. Then $W$ is split into two disjoint open sets, $W=W'\cup W''$, such that $W'\subset A$ and $W''\cap A=\emptyset$. Consider the nested disk sub-system $\{D_x\}$ where $x\in X\cap A$. It is $s(r)$-nested and $\lambda$-scattered in $W'$. Thus the inequalities (8.8) still hold with $I_1$ replaced by $I_1\cap A$. Obviously, $\area(\rho_0,I_1\cap A)\le\area(\rho_0,A)=\md A$. We obtain
$$
\left|\frac{\md A'}{\md A }-1\right|\le C(r).
$$
\qed

\subsection{Application to rational maps}

Let $g$ be a geometrically finite rational map with $\FFF_g\neq\emptyset$. Let $X_0$ be a finite set with $g(X_0)=X_0$. Let $X=\cup_{n\ge 0}g^{-n}(X_0)$. For each point $x\in X$, denote by $n(x)\ge 0$ the minimal integer such that $g^{n(x)}(x)\in X_0$.

\vskip 0.24cm

{\bf Pullback system of $X$}. Pick a disk $U_x\ni x$ for each point $x\in X_0$ such that $\overline{U_x}\smm\{x\}$ contains no critical values of $g$ and $U_x\cap U_y=\emptyset$ if $x\neq y$.

For each point $x\in X$ with $n(x)=1$, let $y=g(x)$ and $U_{x}$ be the component of $g^{-1}(U_{y})$ that contains the point $x$, then $U_{x}$ is also a disk and $\overline{U_x}\smm\{x\}$ is disjoint from $g^{-1}(X_0)$. Since $x\notin\PPP'_g$, we may require that $\overline{U_{x}}\smm\{x\}$ is disjoint from $\PPP_g$ if $U_y$ is chosen to be small enough.

For each point $x\in X$ with $n(x)=n$, let $y=g^{n}(x)$ and let $U_{x}$ be the component of $g^{-n}(U_{y})$ that contains the point $x$; then $U_{x}$ is also a disk and $\overline{U_{y}}\smm\{y\}$ is disjoint from $\PPP_g$. Obviously, $\overline{U_{x}}\smm\{x\}$ is disjoint from $g^{-n}(X_0)$. By Lemma \ref{shrinking}, $\max_{n(x)=n}\{\text{diam}_s U_x\}\to 0$ as $n\to\infty$. We have proved:

\REFPROP{disk} There exist disks $U_x\ni x$ for all points $x\in X$ satisfying the following conditions:

(a) $\overline{U_x}\smm\{x\}$ contains no critical values of $g$.

(b) $\overline{U_x}\smm\{x\}$ is disjoint from $\PPP_g$ if $n(x)\ge 1$.

(c) $\overline{U_x}\smm\{x\}$ is disjoint from $g^{-n}(X_0)$ if $n(x)=n\ge 0$.

(d) $U_x$ is a component of $g^{-1}(U_{y})$ if $g(x)=y$ and $n(x)\ge 1$.

(e) $\max_{n(x)=n}\{\text{diam}_s U_x\}\to 0$ as $n\to\infty$.

We will call $\{U_x\}$ {\bf a pullback system of $X$}.
\ENDPROP

Let $\{U_x\}_{x\in X}$ be a pullback system of $X$. Given any $r\in (0,1]$, for each point $x\in X_0$, let $U_x(r)=\chi^{-1}(\DD(r))$, where $\chi: U_x\to\DD$ is a conformal map with $\chi(x)=0$.

For each point $x\in X$ with $n(x)=1$, let $y=g(x)$ and let $U_{x}(r)$ be the component of $g^{-1}(U_{y}(r))$ that contains the point $x$. For each point $x\in X$ with $n(x)=n$, let $y=g^n(x)$ and let $U_x(r)$ be the component of $g^{-n}(U_{y}(r))$ that contains the point $x$. Then $\{U_x(r)\}$ satisfies conditions (a)--(d) and hence is also a pullback system of $X$.

\vskip 0.24cm

The notation $U_x(r)$ may cause misunderstanding. For each point $x\in X$, let
$$
d(x)=
\begin{cases}
\deg_x g^{n(x)} & \text{if } n(x)\ge 1, \\
1 & \text{if }x\in X_0
\end{cases}
$$
and
$$
d(X)=\sup\{d(x):\, x\in X\}<\infty.
$$
Let $x\in X$ with $n(x)\ge 1$. Let $y=g^{n(x)}(x)$. Then $g^{n(x)}: U_x\smm\{x\}\to U_y\smm\{y\}$ is a covering of degree $d(x)$. Let $\chi_x$ and $\chi_y$ be conformal maps from $U_x$ and $U_y$, respectively, to the unit disk $\DD$, such that $\chi_x(x)=\chi_y(y)=0$. Then $\chi_y\circ g^{n(x)}\circ\chi_x^{-1}$ is a covering of $\DD^*$ of degree $d(x)$. It follows that $\chi_x(U_x(r^{d(x)}))=\DD(r)$.

\vskip 0.24cm

{\bf Nested disk system}. Let $\{U_x\}_{x\in X}$ be a pullback system of $X$. Let $n\ge 0$ be given.

Set $D_x=U_x$ if $n(x)=n$. For each point $x\in X$ with $n(x)=n-1$, let
$$
b(x)=\{y\in X:\, n(y)=n\text{ and }D_y\cap\partial U_x\neq\emptyset\}.
$$
Set $D_x$ to be the component of $U_x\smm\cup_{y\in b(x)}\overline{D_y}$ that contains the point $x$. Then $D_x$ is a disk. Inductively, for each point $x\in X$ with $0\le n(x)=k<n-1$, let
$$
b(x)=\{y\in X:\, k<n(y)\le n\text{ and }D_y\cap\partial U_x\neq\emptyset\}.
$$
Set $D_x$ to be the component of $U_x\smm\cup_{y\in b(x)}\overline{D_y}$ that contains the point $x$. Then $D_x$ is a disk. Obviously, $\{D_x\}_{n(x)\le n}$ is a nested disk system. The following properties are easy to check.

\vskip 0.24cm

{\noindent\bf Properties of the nested disk system}. Let $x\in X$ be a point.

(1) $D_x\subset U_x$ .

(2) If $D_y\subset D_x$ and $y\neq x$, then $n(y)>n(x)$.

(3) Given $r_0\in (0,1)$, if any $U_y$ with $n(y)>n(x)$ and $U_y\cap\partial U_x\neq\emptyset$ is disjoint from $U_x(r_0)$, then $U_x(r_0)\subset D_x$.

Recall that for a nested disk system $\{D_x\}$ and for $r\in (0,1)$, $D_x(r)$ is defined by $D_x(r)=\chi^{-1}(\DD(r))$ for a conformal map $\chi: D_x\to\DD$.
Applying the Schwarz Lemma, we have:

(4) $D_x(r)\subset U_x(r^{d(x)})$ for $r\in (0,1)$.

(5) If $U_x(r_0)\subset D_x$, then $U_x(r_0\cdot r^{d(x)})\subset D_x(r)$ for $r\in (0,1)$.

\begin{proposition}\label{nested}
Let $\{U_x\}_{x\in X}$ be a pullback system of $X$. There exist constants $r_0, \lambda\in (0,1)$ and an increasing function $s(r): (0,1]\to (0,1]$ with $s(r)\ge r$ and $s(r)\to 0$ as $r\to 0$, such that for any integer $n\ge 0$, the nested disk system $\{D_x\}$ generated from the pullback system $\{U_x(r_0)\}$ at step $n$ is $s(r)$-nested and $\lambda$-scattered in $\cup_{n(x)\le n}U_x$. Moreover, there exists a constant $r'_0\in (0, r_0)$ such that for $x\in g^{-n}(X_0)$ and $r\in (0,1]$,
$$
U_x(r'_0\cdot r^{d(x)})\subset D_x(r).
$$
\end{proposition}

\beginp
Denote by $p(X)\ge 1$ the maximum of the periods of cycles in $X_0$. By Proposition \ref{disk} (e), there exists a constant $r_1\in (0,1)$ such that $U_x(r_1)\cap U_y(r_1)=\emptyset$ if $x\neq y$ and $n(x), n(y)\le p(X)$.

Let $x\in X_0$ be a parabolic or repelling periodic point with period $p\ge 1$. Let $V$ be a repelling flower of $(g^p, x)$ if $x$ is parabolic, or a linearizable disk at the point $x$ if $x$ is repelling, i.e. $g^p$ is injective in $V$ and $\overline{V}\subset g^p(V)$. Since $\FFF_g\neq\emptyset$, there exists a disk $E_0\subset\FFF_g\cap V$ such that $\overline{E_0}\cap(\overline{X}\cup\PPP_g)=\emptyset$.  Let $\{E_{k}\}$ ($k\ge 1$) be the backward orbit of $E_0$ in $V$ under $g^p$, i.e. $g^p(E_{k+1})=E_{k}$ for $k\ge 0$. Then $E_{k}$ is also a disk and $\{E_{k}\}$ converges to the point $x$ as $k\to\infty$. Thus there exists an integer $k_0\ge 0$ such that $E_k\subset U_x(r_1)$ for $k\ge k_0$. For simplicity, we write $E_k$ in place of $E_{k+k_0}$.

Again by Proposition \ref{disk} (e), we have $U_y(r_1)\cap E_0=\emptyset$ if $n(y)$ is large enough, since $E_0\cap\overline{X}=\emptyset$. Therefore there exists a constant $r(x)\in (0, r_1)$ such that $U_y(r(x))\cap E_0=\emptyset$ for all points $y\in X\smm\{x\}$. It follows that $U_y(r(x))\cap E_1=\emptyset$ if $n(y)\ge p$. Hence $U_y(r(x))\cap E_1=\emptyset$ for all points $y\in X\smm\{x\}$ since $E_1\subset U_x(r_1)$ and $U_y(r_1)\cap U_x(r_1)=\emptyset$ if $y\neq x$ and $n(y)\le p$. Inductively, for each $k\ge 1$,
$$
U_y(r(x))\cap E_{k}=\emptyset\quad\text{for all points }y\in X\smm\{x\}.
$$

Let $r_2=\min\{r(x)\}$, the minimum over all parabolic and repelling points $x\in X_0$. Then for each parabolic or repelling point $x\in X_0$, we have a sequence of disks $\{E_{x,k}\}$ in $U_x$ which converge to the point $x$ such that $E_{x,k}$ are disjoint from $U_y(r_2)$ for all points $y\in X\smm\{x\}$.

Let $n_0\ge 2$ be an integer such that $x\not\in g^{-1}(\PPP_g)$ if $n(x)\ge n_0$. Then there exists a constant $r_0\in (0,r_2)$ such that the $U_x(r_0)$ are pairwise disjoint for all points $x\in g^{-n_0}(X_0)$. In the sequel we will write $U'_x:=U_x(r_0)$ for simplicity.

For each point $x\in g^{-n_0}(X_0)$ and each point $y\in X$, if $\text{\rm dist}_s(x, U'_y)\to 0$ then $n(y)\to\infty$ and hence $\text{diam}_s U'_y\to 0$ uniformly. Thus there exist a constant $r_3\in (0, r_0)$ and an increasing function $s_1: (0,r_3]\to (0,1]$ with $s_1(r)\ge r$ and $s_1(r)\to 0$ as $r\to 0$, such that for each point $x\in g^{-n_0}(X_0)$ and each point $y\in X$, if $U'_y\cap U'_x(r)\neq\emptyset$ for some $r\le r_r$, then $n(y)>n_0$ and $U'_y\subset U'_x(s_1(r))$.

Let $x\in X$ be a point with $n(x)>n_0$. Set $k=n(x)-n_0$. Then $n(g^k(x))=n(x)-k=n_0$. For any point $y\in X$ with $n(y)>n(x)$, if $U'_y\cap U'_x(r)\neq\emptyset$ for some $r\le r_1$, then
$$
U'_{g^k(y)}\cap U'_{g^k(x)}(r)=g^k(U'_y\cap U'_x(r))\neq\emptyset.
$$
Thus $g^k(U'_y)=U'_{g^k(y)}\subset U'_{g^k(x)}(s_1(r))=g^k(U'_x(s_1(r)))$ and hence $U'_y\subset U'_x(s_1(r))$.

For any integer $n\ge 0$, let $\{D_x\}$ be the nested disk system generated from the pullback system $\{U'_x\}$ at step $n$. From Property (3), we know that $U'_x(r_3)\subset D_x$ if $n(x)\le n$. Set $r'_0=r_0\cdot r_3$. By Property (5),
$$
U_x(r'_0\cdot r^{d(x)})=U'_x(r_3\cdot r^{d(x)})\subset D_x(r)\quad\text{for }r\in (0,1].
$$

Because $s_1(r)\to 0$ as $r\to 0$, there exists a constant $r_4\in (0, r_3)$ such that $s_1(r_4)<r_3$. For any two points $x,y\in X$ with $n(y)>n(x)\ge n$, if $D_y\cap D_x\neq\emptyset$, then $D_y\subset D_x$ by Property (2). Assume further that $D_y\cap D_x(r)\neq\emptyset$ for some $r<r_4$; then $U'_y\cap U'_x(r)\neq\emptyset$ since $D_x(r)\subset U'_x(r^{d(x)})\subset U'_x(r)$ by Property (4). Thus
$$
D_y\subset U'_y\subset U'_x(s_1(r))\subset D_x\left(\left( \frac{s_1(r)}{r_3} \right)^{ \frac{1}{d(x)} }\right).
$$
Set $s(r)=(s_1(r)/r_3)^{\frac{1}{d(X)}}$ if $0<r<r_4$ and $s(r)=1$ if $r_4\le r\le 1$. Then $\{D_x\}$ is an $s(r)$-nested disk system.

Now we want to prove that the nested disk system $\{D_x\}$ is $\lambda$-scattered in $\cup_{n(x)\le n}U_x$. Let $x\in X$ be a point with $n(x)\le n_0$.
If it is eventually attracting or super-attracting, then there exists a constant $r(x)\in (0,1)$ such that $D_x(s(r(x)))\smm\{x\}$ contains no eventually periodic points. Set $E_x:=D_x(r(x))$. Then for any disk $D_y$ with $y\neq x$ and $D_y\subset D_x$, $E_x\cap D_y=\emptyset$. Otherwise $y\in D_y\subset D_x(s(r(x)))$. This ia a contradiction.

Assume that the point $x\in X$ with $n(x)\le n_0$ is eventually parabolic or repelling. Set $z=g^{n(x)}(x)$. From the discussion at the beginning of the proof, there exists a sequence of disks $\{E_{z,k}\}$ in $U_z$ which converges to the point $z$ such that the disks are disjoint from $U_y(r_0)$ for all points $y\in X\smm\{z\}$. In particular, there exists an integer $k\ge 1$ such that $E_{z,k}\subset g^{n(x)}(D_x)$. Let $E_x\subset D_x$ be a component of $g^{-n(x)}(E_{x,k})$. Then for any disk $D_y\subset D_x$ with $y\neq x$, $D_y\cap E_x=\emptyset$. Otherwise, $g^{n(x)}(D_y)$ intersects $g^{n(x)}(E_x)$. It follows that $g^{n(x)}(U_y(r_0))=U_{g^{n(x)}(y)}(r_0)$ intersects $E_{x,k}$, since $n(y)>n(x)$. This is a contradiction.

In summary, for each point $x\in X$ with $n(x)\le n_0$, there exists a disk $E_x\subset D_x$ such that for any disk $D_y\subset D_x$ with $y\neq x$, $D_y\cap E_x=\emptyset$.

For any univalent map $h:U_x\to\CC^*$, applying the Koebe Distortion Theorem for $h$, we get a constant $\lambda_x\in (0,1)$ such that
$$
\area(\rho_*,h(E_x))\ge(1-\lambda_x)\area(\rho_*,h(U'_x)).\eqno{(8.9)}
$$
It follows that
$$
\area(\rho_*,h(E_x))\ge(1-\lambda_x)\area(\rho_*,h(D_x))\quad\text{and}
$$
$$
\area(\rho_*,h(D_x\smm E_x))\le\lambda_x\cdot\area(\rho_*,h(D_x)).\eqno{(8.10)}
$$
Set $\lambda=\max\{\lambda_x\}$, the maximum over all points $x\in X$ with $n(x)\le n_0$. If $n\le n_0$, the lemma is proved. Otherwise, for each point $x\in X$ with $n_0<n(x)\le n$, let $k=n(x)-n_0$ and $z=g^k(x)$. Then $n(z)=n(x)-k=n_0$. Let $E_x\subset D_x$ be a component of $g^{-k}(E_z)$, where $E_z$ is defined above. Then for any disk $D_y\subset D_x$ with $y\neq x$, we have $D_y\cap E_x=\emptyset$.

For any univalent map $h: U_x\to\CC^*$, applying the Koebe Distortion Theorem again for $h\circ g^{-k}$ on $U_z$, we get the inequality (8.9), where $\lambda_x$ should be replaced by $\lambda$. The inequality (8.10) follows. This completes the proof.
\qed

\vskip 0.24cm

Combining Proposition \ref{nested} and Theorem \ref{distortion}, we have:

\REFTHM{distortion1}
Let $\{U_x\}_{x\in X}$ be a pullback system of $X$. Let $V$ be an open set compactly contained in $\cbar\smm\overline{X}$. Then there exist a constant $r_0\in (0,1)$ and an increasing function $C(r)$ on $(0, r_0)$ with $C(r)\to 0$ as $r\to 0$, such that for any constant $r\in (0,r_0)$, if
$$
\phi: \cbar\smm\bigcup_{n(x)\le n}\overline{U_x(r)}\to\cbar
$$
is univalent for some integer $n\ge 0$, then $\sD_0(\phi, V)\le C(r)$.
\ENDTHM

Recall that if $q$ is a quotient map of $\cbar$ and $A\subset\cbar$ is an annulus, then $q^{-1}(A)$ is also an annulus. The following lemma will be used in \S9.

\begin{lemma}\label{prepare}
Let $\{U_x\}$ be a pullback system of $X$. Then there exists a constant $r_0\in (0,1)$ such that if

(a) $A_1\subset A_0$ are annuli in $\cbar$ such that $A_1$ is contained essentially in $A_0$ and for any point $x\in X$, $U_x\cap\partial A_0\neq\emptyset$ implies that $U_x\cap A_1=\emptyset$, and

(b) $q$ is a quotient map of $\cbar$ such that
$$
q^{-1}:\, \cbar\smm\bigcup_{n(x)\le n}\overline{U_x(r_0)}\to\cbar
$$
is univalent for some integer $n\ge 0$, then $\md q^{-1}(A_0)\ge(\md A_1)/2$.
\end{lemma}

\beginp
By Proposition \ref{nested}, there exist constants $0<r_1<r'_1<1$ such that for any integer $n\ge 0$, we have $U_x(r_1)\subset D'_x$ if $n(x)\le n$, where $\{D'_x\}$ is the nested disk system generated from the pullback system $\{U_x(r'_1)\}$ at step $n$.

Applying Proposition \ref{nested} for the pullback system $\{U_x(r_1)\}$, there exist constants $\lambda\in (0, 1)$ and $0<r_2<r'_2<r_1$ and an increasing function $s(r): (0,1]\to (0,1]$ with $s(r)\ge r$ and $s(r)\to 0$ as $r\to 0$, such that for any integer $n\ge 0$, the nested disk system $\{D_x\}$ generated from the disk system $\{U_x(r'_2)\}$ at step $n$ is $s(r)$-nested and $\lambda$-scattered in $W:=\cup_{n(x)\le n}U_x(r_1)$ and $U_x(r_2\cdot r^{d(x)})\subset D_x(r)$ if $n(x)\le n$.

Fix any $n\ge 0$. Then $U_x(r_1)\subset D'_x\subset U_x$ if $n(x)\le n$. Since $\{D'_x\}$ is a nested disk system, by assumption, there exists an annulus $B_n\subset A_0$ such that $\partial B_n\cap(\cup_{n(x)\le n}D'_x)=\emptyset$ and $A_1\subset B_n$. In particular, $\partial B_n\cap W=\emptyset$ since $W\subset\cup_{n(x)\le n} D'_x$.

By Theorem \ref{distortion}, there exists a constant $r_3\in (0,r_2)$, which is independent of the choice of $n$, such that for any univalent map
$\phi:\, \cbar\smm\cup_{n(x)\le n}\overline{D_x(r_3)}\to\cbar$, we have
$$
\left|\dfrac{\md B'_n}{\md B_n}-1\right|\le\frac{1}{2},
$$
where $B'_n$ is the annulus bounded by $\phi(\partial B_n)$.

Set $r_0=r_2\circ r_3^{d(X)}$. Then $U_x(r_0)\subset D_x(r_3)$ if $n(x)\le n$. Let $q$ be a quotient map of $\cbar$ such that $q^{-1}$ is univalent in $\cbar\smm\cup_{n(x)\le n}\overline{U_x(r_0)}$; then $q^{-1}$ is univalent in $\cbar\smm\cup_{n(x)\le n}\overline{D_x(r_3)}$. Therefore
$$
\md q^{-1}(A_0)\ge \md q^{-1}(B_n)\ge(\md B_n)/2\ge(\md A_1)/2.
$$
\qed

\section{Hyperbolic-parabolic deformation}

In this section, we will prove Theorem \ref{pinching}. Let $f$ a geometrically finite rational map and let $\sA\subset\sR_f$ be a non-separating multi-annulus. Then there exists an integer $n_0\ge 1$ such that all the filled-in skeletons of level $n_0$ are disjoint from $\PPP_f\cup\Omega_f$. Thus all the filled-in skeletons of level $n$ are disjoint from $f^{-1}(\PPP_f)$ for $n>n_0$. We will prove the theorem under the assumption $n_0=1$ for simplicity. The proof in the case $n_0>1$ has no essential difficulty. The following notation will be used in this section: For $n\ge 0$,

\vskip 0.24cm

$\BBB_n$ is the union of all bands of level $k\le n$,

$\wh\SSS_n$ is the union of all filled-in skeletons of level $k\le n$,

$\TTT_n$ is the union of all band-trees of level $k\le n$, and

$\III_n$ is the interior of $\TTT_n$.

\subsection{Piecewise pinching}

Recall that $\phi_{t,n}$ is the normalized quasiconformal map of $\cbar$ with Beltrami differential $\mu(\phi_{t,n})=\mu(\phi_t)$ on $\BBB_n$ and zero otherwise.

\REFLEM{piece-pinching}
Let $n\ge 0$ be fixed. Then $\{\phi_{t,n}\}$ converges uniformly to a quotient map $\varphi_n$ as $t\to\infty$. For each filled-in skeleton $\wh S\subset\wh\SSS_n$, $\varphi_n(\wh S)$ is a single point. Conversely, for any point $w\in\cbar$, $\varphi_n^{-1}(w)$ is either a single point or a filled-in skeleton in $\wh\SSS_n$.
\ENDLEM

\beginp
We begin by proving that $\{\phi_{t,n}\}$ is equicontinuous for the given $n$. For each filled-in skeleton $\wh S_i\subset\wh\SSS_n$, by Lemma \ref{nested-annuli}, there exist a constant $M>0$, a sequence $\{t_k\ge 0\}$ $(k\ge 0)$ with $t_k\to\infty$ as $k\to\infty$, and a sequence of nested disks $\{U_k\}$ such that $\overline{U_{k+1}}\subset U_k$, $\cap_{k\ge 0}U_k=\wh S_i$ and $\md \phi_{t,n}(U_k\smm\overline{U_{k+1}})>M$ for $t\ge t_k$. Thus
$$
\md \phi_{t,n}(U_0\smm\overline{U_{k+1}})\ge kM\quad\text{for }t\ge t_k.
$$
For any $\varepsilon>0$, by the normalization condition and the above inequality, there is an integer $k>0$ such that
$$
\text{diam}_s\phi_{t,n}(U_{k})<\varepsilon\quad\text{for } t\ge t_k.\eqno{(9.1)}
$$

Let $D_i= U_{k+1}$. Since $\{\phi_{t,n}\}$ is uniformly quasiconformal for $t<t_k$, there exists a constant $\delta_i>0$ such that for any two points $z_1\in D_i$ and $z_2\in\cbar$ with $\text{\rm dist}_s(z_1, z_2)<\delta_i$,
$$
\text{\rm dist}_s(\phi_{t,n}(z_1), \phi_{t,n}(z_2))<\varepsilon\quad\text{for }t\ge 0.\eqno{(9.2)}
$$

Let $D$ be the union of all the domains $D_i$ taken above for all filled-in skeletons $\wh S_i\subset\wh\SSS_n$. Let $W=\cbar\smm\wh\SSS_n$. Then $W$ is connected and $\cbar\smm D\subset W$. Thus there exists a domain $V\Subset W$ such that $D\subset V$. The family $\{\phi_{t,n}\}$ is uniformly quasiconformal on $V$ and hence is equicontinuous by the normalization condition. Thus there exists a constant $\delta_0<0$ such that for any points $z_1\in\cbar\smm D$ and $z_2\in\cbar$ with $\text{\rm dist}_s(z_1, z_2)<\delta_0$, the inequality (9.2) holds.

Set $\delta=\min\{\delta_i\}$. Then (9.2) holds for any two points $z_1, z_2\in\cbar$ with $\text{\rm dist}_s(z_1, z_2)<\delta$. Thus $\{\phi_{t,n}\}$ is equicontinuous.

Let $\varphi_n$ be a limit of $\phi_{t,n}$ as $t\to\infty$. By (9.1), $\varphi_n(\wh S_i)$ is a single point for each filled-in skeleton $\wh S_i\subset\wh\SSS_n$. Since $\varphi_n$ is a quasiconformal map on any domain compactly contained in $W$, it is injective on $W$.

If $\wt \varphi_n$ is also a limit of $\{\phi_{t,n}\}$, then $\wt\varphi_n\circ\varphi_n^{-1}$ is a well-defined homeomorphism of $\cbar$, which is holomorphic except on a finite set. Thus it is a global conformal map. So $\wt\varphi_n=\varphi_n$ by the normalization. Therefore $\{\phi_{t,n}\}$ converges uniformly to $\varphi_n$ as $t\to\infty$.
\qed

\subsection{The candidate pinching limit and the proof of Theorem \ref{pinching}}

Let $n\ge 1$ be fixed. Then $\{\phi_{t,n-1}\circ f\circ\phi_{t,n}^{-1}\}$ are rational maps which converge uniformly to a rational map $g_n$ as $t\to\infty$, by Lemma \ref{piece-pinching} and Lemma \ref{rational-s}. Obviously,
$$
\varphi_{n-1}\circ f=g_n\circ\varphi_{n}.
$$

Pick a Koebe space system $\{N(T)\}$ for the band-trees. Denote by $N_n$ the union of the Koebe spaces of the level $n$ band-trees. Then $\overline{N}_1$ is disjoint from $\PPP_f\cup\Omega_f$ by the assumption. Note that $\varphi_{1}\circ\varphi_0^{-1}$ is a well-defined quotient map of $\cbar$, which is univalent in $\cbar\smm\varphi_0(N_1)$. Thus there exists a homeomorphism $\theta_1$ of $\cbar$ such that
$\theta_1=\varphi_{1}\circ\varphi_0^{-1}$ on $\cbar\smm\varphi_{0}(N_1)$.

Define $G=g_{1}\circ\theta_1$. Then $G$ is a branched covering of $\cbar$ which is holomorphic in $\cbar\smm\varphi_0(\overline{N_1})$. Let $\wt\varphi_0=\theta^{-1}\circ\varphi_{1}$. Then
$$
\begin{cases}
\varphi_0\circ f=G\circ\wt\varphi_0 & \text{on }\cbar, \\
\wt\varphi_0=\varphi_0 & \text{on }\cbar\smm N_1.
\end{cases}
$$
Thus $\wt\varphi_0$ is isotopic to $\varphi_0$ rel $\cbar\smm N_1$ since $N_1$ is a disjoint union of disks whose closures are disjoint from $\PPP_G$. Obviously, $\PPP_G=\varphi_0(\PPP_f)\subset\cbar\smm\varphi_0(\overline{N}_1)$. See the commutative diagram:
\[
\xymatrix{
\cbar\ar[d]_{f}\ar[r]^{\varphi_{1}} & \cbar\ar[d]_{g_{1}} & \cbar\ar[l]_{\theta_{1}} \ar[d]^{G} \\
\cbar\ar[r]^{\varphi_{0}} & \cbar & \cbar\ar[l]_{\text{\rm id}}}
\]

Recall that $\PPP^s_G\subset\PPP_G$ is the set of super-attracting periodic points of $G$.

\REFLEM{f2G}
The map $G$ is a semi-rational map which is holomorphic in a neighborhood of $\PPP^s_G$. If $\UUU$ is a fundamental set of $f$, then $\varphi_0(\UUU)$ contains a fundamental set of $G$.
\ENDLEM

\beginp
Note that $\PPP_G=\varphi_0(\PPP_f)$ is contained in $\cbar\smm\varphi_0(\overline{N_1})$ and $G$ is holomorphic in $\cbar\smm\varphi_0(\overline{N_1})$. Thus $G$ is geometrically finite and holomorphic in a neighborhood of $\PPP_G$. Let $\wt Y=\varphi_0(\wh\SSS_0)$. Then $\wt Y\subset\PPP'_G$ and $G(\wt Y)=\wt Y$. We only need to show that each point $y\in\wt Y$ is a parabolic periodic point of $G$ and that there exists an attracting flower at $y$ in $\varphi_0(\UUU)$ such that each of its petals contains points of $\PPP_G$.

Let $\wh S=\varphi_0^{-1}(y)$. This is a periodic filled-in skeleton. Denote by $X$ the set of attracting or parabolic periodic points in $\wh S$. This is a non-empty finite set. For simplicity, we may assume that each point in $X$ is a fixed point of $f$.

Let $\UUU\subset\cbar\smm N_1$ be a fundamental set of $f$. For each point $x\in X$, if $x$ is attracting, then $x\in\UUU$. Thus there exists a disk $D_x\subset\UUU$ with $x\in D_x$ such that $f$ is injective on $D_x$ and $\overline{f(D_x)}\subset D_x$. Moreover, we may require that $\partial D_x$ intersects each component of $\wh S\smm\{x\}$ at a single point or a closed arc. Then for each component $U$ of $D_x\smm\wh S$, we see that $V=\varphi_0(U)$ is a disk and $\overline{G(V)}\subset V\cup\{y\}$. By Proposition \ref{infty}, the domain $U$ contains infinitely many points of $\PPP_f$. Thus the disk $V$ contains infinitely many points of $\PPP_G$.

Now suppose that $x\in X$ is parabolic. Then there exists an attracting flower $\VVV_x$ of $f$ at $x$ such that $\VVV_x\subset\UUU$. We may also require that each component of $\partial\VVV_x\smm\{x\}$ is either disjoint from $\wh S$ or intersects each component of $\wh S\smm\{x\}$ at a single point or a closed arc. Then for each component $U$ of $\VVV_x\smm\wh S$,  $V=\varphi_0(U)$ is a disk and $\overline{G(V)}\subset V\cup\{y\}$. By Proposition \ref{infty}, the domain $U$ contains infinitely many points of $\PPP_f$. Thus the disk $V$ contains infinitely many points of $\PPP_G$.

Denote by $\VVV_y$ the union of $\varphi_0(U)$ for all components $U$ of $D_x\smm\wh S$ and $\VVV_x\smm\wh S$ for all point $x\in X$. Let $\{w_n=G^{n}(w)\}$ be an orbit in $\cbar\smm\varphi_0(N_1)$ converging to the point $y$ as $n\to\infty$. Then $\{z_n=\varphi_0^{-1}(w_n)\}$ converges to a point $x\in X$ and $f(z_n)=z_{n+1}$. Thus once $n$ is large enough, $z_n$ is contained in $D_x$ if $x$ is attracting, or $\VVV_x$ if $x$ is parabolic. So $w_n$ is contained in $\VVV_y$ once $n$ is large enough. Therefore the point $y$ is a parabolic fixed point of $G$ and $\VVV_y\subset\varphi_0(\UUU)$ is an attracting flower of $(G,y)$ by Lemma \ref{parabolic}. Since each component of $\VVV_y$ contains infinitely many points of $\PPP_G$, the branched covering $G$ is a semi-rational map and $\varphi_0(\UUU)$ contains a fundamental set of $G$. \qed

\REFLEM{no-obs-2}
The semi-rational map $G$ has neither Thurston obstructions nor connecting arcs.
\ENDLEM

\beginp
Let $\Gamma$ be a multicurve in $\cbar\smm\PPP_G$. Since $\varphi_{0}(N_1)$ is a disjoint union of disks whose closures are disjoint from $\PPP_G$, we may assume that $\gamma$ is disjoint from $\varphi_{0}(N_1)$ for each $\gamma\in\Gamma$. Let $\Gamma'=\{\varphi_0^{-1}(\gamma):\, \gamma\in\Gamma\}$. Then $\Gamma'$ is a multicurve of $f$ and $\lambda(\Gamma)=\lambda(\Gamma')$ since $\varphi_{0}\circ f=G\circ\varphi_{0}$ on $\cbar\smm N_1$ and $\PPP_G=\varphi_0(\PPP_f)$. Thus $\lambda(\Gamma)<1$ since $f$ has no Thurston obstructions. Therefore $G$ has no Thurston obstructions.

Assume that $\beta_0\subset\cbar\smm\PPP_G$ is a connecting arc which connects two points $y^{\pm}$ in $\PPP'_G$, i.e., either $y^-\neq y^+$ or $y^-=y^+$ and both components of $\cbar\smm(\beta_0\cup\{y^+\})$ contain points of $\PPP_G$; $\beta_0$ is disjoint from a fundamental set of $G$ and $G^{-p}(\beta_0)$ has a component $\beta_1$ isotopic to $\beta_0$ rel $\PPP_G$ for some integer $p\ge 1$. We may assume that both $\beta_0$ and $\beta_1$ are disjoint from $\varphi_{0}(N_1)$.

By Lemma \ref{connecting0}, we may further assume that $\beta_1$ coincides with $\beta_0$ in a neighborhood of the endpoints $\{y^{\pm}\}$. Then $\beta_0$ can be divided into three arcs: $\beta^+_0$, $\beta^-_0$ and the middle piece $\beta^0_0$, such that $\beta^{\pm}_1:=(G^{p}|_{\beta_1})^{-1}
(\beta^{\pm}_0)\subset\beta^{\pm}_0$. Thus $\{\beta^{\pm}_k:=(G^{p}|_{\beta_1})^{-k}(\beta^{\pm}_0)\}$ converges to the endpoints $\{y^{\pm}\}$ as $k\to\infty$.

Set $\delta^{\pm}_k=\varphi_{0}^{-1}(\beta^{\pm}_k)$ for all $k\ge 0$. Then $\delta^{\pm}_{k+1}\subset\delta^{\pm}_k$ and $f^{p}(\delta^{\pm}_{k+1})=\delta^{\pm}_k$. Thus $\delta^{\pm}_0$ lands on a repelling or parabolic periodic point $a^{\pm}$ of $f$ by Lemma \ref{arc} and $\{\delta^{\pm}_k\}$ converges to $a^{\pm}$ as $k\to\infty$. Note that $y^{\pm}=\varphi_0(a^{\pm})$. Thus $a^{\pm}\in\PPP'_f\cup\wh\SSS_0$.

Let $\delta_0=\varphi_{0}^{-1}(\beta_0)$. Then $f^{-p}(\delta_0)$ has a unique component $\delta_1$ isotopic to $\delta_0$ rel $\PPP_f\cup\{a^{\pm}\}$ since $\wt\varphi_0$ is isotopic to $\varphi_0$ rel $\PPP_f$. So $\delta^{\pm}_{1}\subset\delta_1$. For each $k\ge 2$, $f^{-kp}(\delta_0)$ has a unique component $\delta_k$ isotopic to $\delta_0$ rel $\PPP_f\cup\{a^{\pm}\}$. Therefore $\delta^{\pm}_k\subset\delta_k$. Set $\delta^0_k=\delta_k\smm (\delta^+_k\cup\delta^-_k)$. Then its spherical diameter converges to zero as $k\to\infty$ by Lemma \ref{shrinking}. Thus $a^-=a^+$ and $\delta_k$ converges to the point $a^+$. In particular, the Jordan curve $\delta_0\cup\{a^+\}$ does not separate the points of $\PPP_f$. Thus $y^{-}=y^{+}$ and one component of $\cbar\smm(\beta_0\cup\{y^{+}\})$ is disjoint from $\PPP_G$. This is a contradiction.
\qed

\REFLEM{f2g}
There exists a rational map $g$, a sequence of normalized quotient maps $\{\xi_n\}$ $(n\ge 0)$ of $\cbar$, a fundamental set $\UUU_f$ of $f$ and a fundamental set $\UUU_g$ of $g$ such that the following conditions hold:

(a) $\xi_0(\wh\SSS_0)$ consists of some parabolic cycles of $g$.

(b) $\III_0\cup\PPP^s_F\subset\UUU_f$ and $\xi_0(\UUU_f)=\UUU_g\cup\xi_0(\wh\SSS_0)$.

(c) $g\circ\xi_{n+1}=\xi_n\circ f$.

(d) $\xi_{n+1}$ is isotopic to $\xi_{n}$ rel $f^{-n}(\overline{\UUU_f}\cup\PPP_f)$.

(e) $\varphi_n\circ\xi_n^{-1}$ is quasiconformal on $\cbar$ and holomorphic in $g^{-n}(\UUU_g)$.
\ENDLEM

\beginp
Let $Y_G=\varphi_0(\wh\SSS_0)$ and $\WWW_G=\varphi_0(\III_0\smm\wh\SSS_0)$. Then $Y_G\subset\PPP'_G$ and $\varphi_0(\III_0)=\WWW_G\cup Y_G$. Each cycle in $Y_G$ is parabolic and each component of $\WWW_G$ is a sepal of $G$ at a point in $Y_G$.

Applying Theorem \ref{existence} and Lemma \ref{no-obs-2} to $G$, we get a normalized c-equivalence $(\psi_0, \psi_1)$ from $G$ to a rational map $g$. By Lemma \ref{c2q}, we may choose the c-equivalence such that $\psi_0$ is quasiconformal on $\cbar$ and $\psi_1=\psi_0$ in a neighborhood $N$ of $\PPP'_G$.

Since $\overline{\WWW}_G\cap\PPP_G=Y_G$, using a similar argument as in the proof of Lemma \ref{c2q}, we may require further that $\WWW_G\Subset N$, $\psi_0$ is holomorphic in $\WWW_G$ and $\psi_1=\psi_0$ on $\WWW_G$. Thus there exists a fundamental set $\UUU_G$ of $G$ with $\WWW_G\subset\UUU_G$ such that $\psi_0$ is quasiconformal on $\cbar$ and holomorphic in $\UUU_G$, and $\psi_1$ is isotopic to $\psi_0$ rel $\UUU_G\cup\PPP_G$ by Lemma \ref{isotopy}. We may also assume that $\PPP^s_G\subset\UUU_G$ since $G$ is holomorphic in a neighborhood of $\PPP_G$.

Set $\UUU_f=\varphi_0^{-1}(\UUU_G)\cup\III_0$. Then $\varphi_0(\UUU_f)=\UUU_G\cup Y_G$. Thus $\UUU_f$ is a fundamental set of $f$. Set $\xi_0=\psi_0\circ\varphi_0$ and $\xi_1=\psi_1\circ\wt\varphi_0$. Then $g\circ\xi_1=\xi_0\circ f$ and $\xi_1$ is isotopic to $\xi_0$ rel $\overline{\UUU_f}\cup\PPP_f$. Set $\UUU_g=\psi_0(\UUU_G)$. This is a fundamental set of $g$ by Lemma \ref{equivalent}. We have $\xi_0(\wh\SSS_0)=\psi_0(Y_G)$ and $\xi_0(\UUU_f)=\UUU_g\cup\psi_0(Y_G)$. Obviously, $\varphi_0\circ\xi_0^{-1}=\psi_0^{-1}$ is quasiconformal on $\cbar$ and holomorphic in $\UUU_g$.

Let $\xi_{n+1}$ be the lift of $\xi_{n}$ for $n\ge 1$. Then $g\circ\xi_{n+1}=\xi_n\circ f$ and $\xi_{n+1}$ is isotopic to $\xi_n$ rel $f^{-n}(\overline{\UUU_f}\cup\PPP_f)$. Recall that
$$
g_{n+1}\circ(\varphi_{n+1}\circ\xi_{n+1}^{-1})=(\varphi_n\circ\xi_n^{-1})\circ g.
$$
Thus $\varphi_n\circ\xi_n^{-1}$ is quasiconformal on $\cbar$ and holomorphic in $g^{-n}(\UUU_g)$. See the following commutative diagrams:
\[
\xymatrix{
\cbar\ar[d]_{f}\ar[r]^{\wt \varphi_0} & \cbar\ar[d]_{G}\ar[r]^{\psi_1} & \cbar\ar[d]^{g}
& & \cbar\ar[d]_{g_{n+1}} & \cbar\ar[d]_{f}\ar[l]_{\varphi_{n+1}}\ar[r]^{\xi_{n+1}} & \cbar\ar[d]^{g} \\
\cbar\ar[r]^{\varphi_0} & \cbar\ar[r]^{\psi_0} & \cbar
& & \cbar & \cbar\ar[l]_{\varphi_n}\ar[r]^{\xi_n} & \cbar}
\]
\qed

\vskip 0.24cm

Set $\xi_{t,n}:=\xi_n\circ \phi_t^{-1}$ for $t,n\ge 0$. For each periodic filled-in skeleton $\wh S$, $\xi_0(\wh S)$ is a parabolic periodic point in $\PPP'_g$. Let $Y_0=\xi_0(\wh\SSS_0)$. Then $g(Y_0)=Y_0$. Let $Y=\cup_{n\ge 0}Y_n$. We will prove the following lemma in \S9.3 and \S9.4.

\REFLEM{ft2g}
(1) For fixed $t\ge 0$, $\{\xi_{t,n}\}$ converges uniformly to a quotient map $q_t$ of $\cbar$ as $n\to
\infty$, and $q_t\circ f_t=g\circ q_t$.

(2) For each filled-in skeleton $\wh S$, $q_t\circ\phi_t(\wh S)$ is a single point. Conversely, for each point $w\in\cbar$, $(q_t\circ\phi_t)^{-1}(w)$ is a single point if $w\notin Y$, or else is a filled-in skeleton.

(3) The sequence $\{q_t\}$ converges uniformly to the identity as $t\to\infty$.
\ENDLEM

{\noindent\it Proof of Theorem \ref{pinching}}.
By Lemma \ref{ft2g}, $q_t\circ f_t=g\circ q_t$ and $\{q_t\}$ converges uniformly to the identity as $t\to\infty$. It follows that $\{f_t\}$ converges uniformly to $g$ as $t\to\infty$.

Because $\xi_{t,n}=\xi_n\circ \phi_t^{-1}$ and $\{\xi_{t,n}\}$ converges uniformly to $q_t$ as $n\to\infty$, we get $q_0=q_t\circ\phi_t$. Thus $\{\phi_t\}$ converges uniformly to $q_0$ as $t\to\infty$ since $\{q_t\}$ converges uniformly to the identity as $t\to\infty$.

Let $X_0=\wh\SSS_0\cap\JJJ_f$ and $X=\cup_{n\ge 0}f^{-n}(X_0)$. Then $q_0(X)=Y$ by Lemma \ref{ft2g} (2). Thus $q_0(\JJJ_f)=\JJJ_g$ since $X$ is dense in $\JJJ_f$ and $Y$ is dense in $\JJJ_g$.

\vskip 0.24cm

Denote by $\MMM(g)$ the space of $g$-invariant Beltrami differentials $\wt\mu$ on $\FFF_g$ with $\|\wt\mu\|_{\infty}<1$. For each $\wt\mu\in\MMM(g)$, let $\phi_{\wt\mu}$ be a quasiconformal map of $\cbar$ whose Beltrami differential is $\wt\mu$. Then $g_{\wt\mu}:=\phi_{\wt\mu}\circ g\circ\phi_{\wt\mu}^{-1}$ is a rational map whose holomorphic conjugate class $[g_{\wt\mu}]$ is contained in $\fm[g]$. Conversely, since $\JJJ_g$ has zero area \cite{U}, each element of $\fm[g]$ is represented by such a rational map $g_{\wt\mu}$.

Define an equivalence relation on $\MMM(g)$ by $\wt\mu_1\sim\wt\mu_2$ if $\phi_{\wt\mu_2}\circ\phi_{\wt\mu_1}^{-1}$ is isotopic to a conformal map of $\cbar$ rel $\PPP_{g_{\wt\mu_1}}$. Then $\wt\mu_1\sim\wt\mu_2$ if and only if there exists a holomorphic conjugacy between $g_{\wt\mu_1}$ and $g_{\wt\mu_2}$ in the corresponding isotopy class (refer to Lemma \ref{qc-conj} or \cite{McS}).

Recall that $\wt\sR_g$ is the set of wandering points in the attracting and parabolic basins of $g$ whose orbits contains no points of $\PPP_g$, and $\pi_g: \wt\sR_g\to\sR_g$ is the natural projection. Each $\wt\mu\in\MMM(g)$, restricted to $\sR_g$, can be pushed forward to be a Beltrami differential $\mu$ on $\sR_g$, which is exactly the Beltrami differential of the pushforward $\Phi_{\mu}: \sR_g\to\sR_{g_{\wt\mu}}$ of the quasiconformal conjugacy $\phi_{\wt\mu}$. If $\wt\mu_1\sim\wt\mu_2$, then $\Phi_{\mu_2}\circ\Phi_{\mu_1}^{-1}$ is isotopic to a conformal map. Conversely, if $\Phi_{\mu_2}\circ\Phi_{\mu_1}^{-1}$ is isotopic to a conformal map and $\wt\mu_1=\wt\mu_2$ in all super-attracting basins of $g$, then there exists a holomorphic conjugacy between $g_{\wt\mu_1}$ and $g_{\wt\mu_2}$ in the corresponding isotopy class.

Let $\WWW=\xi_0(\III_0)\smm Y_0$. This is a disjoint union of calyxes of $g$, and $\pi_g(\WWW)\subset\sR_g$ is a finite disjoint union of once-punctured disks whose closures are pairwise disjoint.  We claim that for each Beltrami differential $\mu$ on $\sR_g$ with $\|\mu\|_{\infty}<1$, there exists a Beltrami differential $\nu$ on $\sR_g$ with $\|\nu\|_{\infty}<1$, such that $\nu=0$ on $\pi_g(\WWW)$ and $\Phi_{\nu_2}\circ\Phi_{\mu}^{-1}$ is isotopic to a conformal map.

Pick a large once-punctured quasi-disk for each component of $\pi_g(\WWW)$ such that they also have disjoint closures and $\overline{\pi_g(\WWW)}\subset\sQ$, where $\sQ$ is their union. Then there exists a quasiconformal map $\Psi:\, \sR_g\to\Phi_{\mu}(\sR_g)$ such that $\Psi=\Phi_{\mu}$ on $\sR_g\smm\sQ$ and
$\Psi$ is holomorphic in $\pi_g(\WWW)$. Thus $\Psi$ is isotopic to $\Phi_{\mu}$. Let $\nu$ be the Beltrami differential of $\Psi$; then $\nu$ satisfies the conditions in the claim.

For each $\wt\mu\in\MMM(g)$, by the claim, there exists a Beltrami differential $\wt\nu\in\MMM(g)$ such that $\wt\nu=\wt\mu$ on $\cbar\smm\wt\sR_g$, $\wt\nu=0$ on the grand orbit of $\WWW$ and $g_{\wt\nu}=g_{\wt\mu}$. Recall that $q_0^{-1}$ is a univalent map from $\cbar\smm\overline{\cup g^{-n}(\WWW)}$ into $\cbar$. Let $\wt\nu^*$ be the pullback of $\nu$ by $q_0$; then $\wt\nu^*$ is $f$-invariant. It is easy to verify that $g_{\wt\nu}$ is the limit of the pinching path starting from $f_{\wt\nu^*}$. Therefore $\fm[g]\subset\partial \fm[f]$.
\qed

\subsection{Factoring of quotient maps by surgery}

We will prove Lemma \ref{ft2g} in the following two sub-sections. By Lemma \ref{f2g} (c), the sequence $\{\xi_n\}$ converges uniformly on any compact subset of $\FFF_g$. However, we have to take a long detour to prove its global convergence.

The quotient map $\xi_n$ is uniquely determined on band-trees up to step $n$, which contains all its non-trivial fibers. In the remaining part, $\xi_n$ is only determined up to isotopy and asymptotic conformal rigidity near the band-trees. This observation suggests factoring the quotient maps into hard factors which are uniquely determined quotient maps, and soft factors which are uniformly quasiconformal maps with a certain sort of asymptotically conformal rigidity. We will achieve this factoring through a surgery.

\vskip 0.24cm

Denote by $\TTT_n(t)\subset\TTT_n$ the union of all band-trees of $\sA(t)$ with level $k\le n$. Let $\WWW_n$ and $\WWW_n(t)$ be the interiors of $\xi_n(\TTT_n)$ and $\xi_n(\TTT_n(t))$, respectively. Then $\WWW_n=\xi_n(\III_n\smm\wh\SSS_n)$,
$$
\xi_{t,n}:\ \phi_t(\III_n\smm\TTT_n(t))\to \WWW_n\smm\overline{\WWW_n(t)}
$$
is a conformal map and $\xi_{t,n+1}=\xi_{t,n}$ on $\phi_t(\TTT_n)$, by Lemma \ref{f2g}.

Given any $n\ge 0$, glue $\phi_t(\III_n))$ with $\cbar\smm\overline{\WWW_n(t)}$ by $\xi_{t,n}$. The space obtained is a punctured sphere. Thus there exist a unique normalized univalent map
$$
p_{t,n}:\, \cbar\smm\overline{\WWW_n(t)}\to\cbar
$$
and a univalent map
$$
j_{t,n}:\, \phi_t(\III_n))\to\cbar
$$
such that $j_{t,n}(z_1)=p_{t,n}(z_2)$ if and only if $z_1\in\phi_t(\III_n\smm\TTT_n(t))$ and $z_2=\xi_{t,n}(z_1)$. Thus there exists a unique rational map $f_{t,n}$ of $\cbar$ such that:
$$
\begin{cases}
f_{t,n+1}\circ p_{t,n+1}=p_{t,n}\circ g & \text{on }\cbar\smm\overline{\WWW_{n+1}(t)},\text{ and} \\
f_{t,n+1}\circ j_{t,n+1}=j_{t,n}\circ f_t & \text{on }\phi_t(\III_{n+1}).
\end{cases}
$$
Define $j_{t,n}=p_{t,n}\circ\xi_{t,n}$ on $\cbar\smm\phi_t(\III_n)$. Then $j_{t,n}$ is a normalized homeomorphism of $\cbar$ and
$$
f_{t,n+1}\circ j_{t,n+1}=j_{t,n}\circ f_t.
$$
By Lemma \ref{f2g} (e), $\varphi_n\circ\xi_n^{-1}$ is a quasiconformal map of $\cbar$ and is holomorphic in $g^{-n}(\UUU_g)$. Recall that $\varphi_n\circ\phi_t^{-1}$ is holomorphic in $\cbar\smm\phi_t(\TTT_n(t))$. Thus $\xi_{t,n}$ is quasiconformal in $\cbar\smm\phi_t(\TTT_n(t))$. So is $j_{t,n}$. By the definition, $j_{t,n}$ is holomorphic in $\phi_t(\III_n)$. Note that $\TTT_n(t)\smm\III_n$ consists of finitely many points. Hence $j_{t,n}$ is quasiconformal in $\cbar$. Meanwhile, $\xi_n$ is holomorphic in $\UUU_f\smm\TTT_n$ since $\xi_n(\UUU_f)=\xi_0(\UUU_f)=\UUU_g\cup Y_0$. Thus $\xi_{t,n}$ is holomorphic in $\phi_t(\UUU\smm\TTT_n(t))$ and so is $j_{t,n}$. Therefore $j_{t,n}$ is holomorphic in $\phi_t(\UUU\smm(\TTT_n(t)\smm\III_n))$ and hence in $\phi_t(\UUU_f)$.

\vskip 0.24cm

The map $p_{t,n}^{-1}$ can be extended to a quotient map of $\cbar$ by defining:
$$
q_{t,n}=\begin{cases}
p_{t,n}^{-1} & \text{ on } p_{t,n}(\cbar\smm\overline{\WWW_n(t)}), \\
\xi_{t,n}\circ j_{t,n}^{-1} & \text{ on }j_{t,n}\circ\phi_t(\TTT_n).
\end{cases}
$$
By definition, $\xi_{t,n}=q_{t,n}\circ j_{t,n}$ on $\cbar$ for $t,n\ge 0$.  We have proved:

\REFLEM{factor} For any $t, n\ge 0$, there exist a fundamental set $\UUU_f$ of $f$ with $(\PPP^s_f\cup\III_0)\subset\UUU_f$, a normalized quasiconformal map $j_{t,n}$ of $\cbar$ and a normalized quotient map $q_{t,n}$ of $\cbar$ such that the following conditions hold:

(1) $\xi_{t, n}=q_{t,n}\circ j_{t,n}$.

(2) $j_{t,n}$ is holomorphic in $\phi_t(\UUU_f)$.

(3) $q_{t,n}^{-1}$ is univalent in $\cbar\smm\overline{\WWW_n(t)}$.

(4) There exists a rational map $f_{t,n+1}$ such that the following diagram commutes:
\[
\xymatrix{
\cbar\ \ar[d]_{f_t}\ar[r]^{j_{t,n+1}} & \ \cbar\ \ar[d]_{f_{t,n+1}\!\!}\ar[r]^{q_{t,n+1}} & \ \cbar\ \ar[d]^{g} \\
\cbar\ \ar[r]^{j_{t,n}} & \ \cbar\ \ar[r]^{q_{t,n}} & \ \cbar\ }
\]
\ENDLEM

Since $j_{t,n}$ is quasiconformal in $\cbar$ and holomorphic in $\phi_t(\UUU_f)$, applying Lemma \ref{algorithm2} and Lemma \ref{factor} (4), we have:

\REFLEM{g2F}
For a fixed $t\ge 0$, $\{f_{t,n}\}$ converges uniformly to the rational map $f_t$ and $\{j_{t,n}\}$ converges uniformly to the identity as $n\to\infty$.
\ENDLEM

Recall that $Y_0=\xi_0(\wh\SSS_0)$ consists of some parabolic cycles of $g$ and $Y$ is the grand orbit of $Y_0$. Thus $\overline{Y}=\JJJ_g$. Let $\{U_y\}$ be a pullback system of $Y$ as defined in \S8.3. Set $n(y)=0$ if $y\in Y_0$ and $n(y)=n$ if $g^n(y)\in Y_0$ but $g^{n-1}(y)\not\in Y_0$. Since $\WWW_0(t)$ converges to $Y_0$ as $t\to\infty$, for any $r\in (0,1)$, there exists a constant $t_0\ge 0$ such that $\WWW_0(t)\subset\cup_{n(y)=0}U_y(r)$. It follows that $\WWW_n(t)\subset\cup_{n(y)\le n}U_y(r)$. Applying Lemma \ref{factor} (3) and Theorem \ref{distortion1}, we have:

\begin{lemma}\label{distortion2}
For any non-empty open set $V\Subset\FFF_g$, there exists a constant $t_0\ge 0$ and a decreasing function $C(t)$ on $(t_0, \infty)$ with $C(t)\to 0$ as $t\to\infty$ such that $\overline{V}\cap\WWW_n(t_0)=\emptyset$ for $n\ge 0$ and
$$
\sD_0(q_{t,n}^{-1}, V)\le C(t)\quad\text{for $t>t_0$ and $n\ge 0$}.
$$
\end{lemma}

\subsection{Convergence of the hard factors}

We want to start by comparing the maps $\{q_{t,n}\}$ if they share same parameter $t$ or same step $n$. Note that $q_{t,n}\circ j_{t,n}\circ\phi_t=\xi_n$. For $t_1, t_2, n\ge 0$, set
$$
\wp_{t_1, t_2;\, n}=j_{t_1, n}\circ\phi_{t_1}\circ\phi_{t_2}^{-1}\circ j_{t_2,n}^{-1}.
$$
Then $\wp_{t_1, t_2;\, n}$ is a homeomorphism of $\cbar$ and $q_{t_1,n}\circ\wp_{t_1, t_2; n}=q_{t_2,n}$. Note that
$$
\wp_{t_1, t_2;\, n}:\, q_{t_2,n}^{-1}(\cbar\smm\overline{\WWW_n})\to q_{t_1,n}^{-1}(\cbar\smm\overline{\WWW_n})
$$
is holomorphic since $\wp_{t_1, t_2; n}=q_{t_1,n}^{-1}\circ q_{t_2, n}$. On the other hand,
$$
\wp_{t_1, t_2;\, n}=j_{t_1, n}\circ\phi_{t_1}\circ\phi_{t_2}^{-1}\circ j_{t_2,n}:\, j_{t_2,n}\circ\phi_{t_2}(\TTT_n)\to j_{t_1,n}\circ\phi_{t_1}(\TTT_n)
$$
is quasiconformal and its maximal dilatation is equal to the maximal dilatation of $\phi_{t_1}\circ\phi_{t_2}^{-1}$ since both $j_{t_1, n}$ and $j_{t_2,n}$ are holomorphic.

\REFPROP{compare1}
For $t_1, t_2\ge 0$ and $n\ge 0$, there exists a quasiconformal map $\wp_{t_1, t_2;\, n}$ of $\cbar$ whose maximal dilatation is equal to the maximal dilatation of $\phi_{t_1}\circ\phi_{t_2}^{-1}$, such that $q_{t_1,n}\circ\wp_{t_1, t_2; n}=q_{t_2,n}$.
\ENDPROP

By Lemma \ref{f2g} (d), there exists a quotient map $\zeta_{n}$ of $\cbar$ such that $\zeta_n$ is isotopic to the identity rel $f^{-n}(\PPP_f\cup\UUU_f)$ and $\xi_{n}=\xi_{n-1}\circ\zeta_{n}$. It follows that $f\circ\zeta_{n+1}=\zeta_n\circ f$.  By Lemma \ref{algorithm1}, $\{\zeta_{n}\}$ converges uniformly to the identity as $n\to\infty$. Set
$$
\eta_{t, n}=(j_{t,n-1}\circ\phi_t)\circ\zeta_{n}\circ (j_{t,n}\circ\phi_t)^{-1}.
$$
Then for a fixed $t\ge 0$, $\{\eta_{t,n}\}$ also converges uniformly to the identity as $n\to\infty$ by Lemma \ref{g2F}. From $q_{t,n}\circ j_{t,n}\circ\phi_t=\xi_n$ and $\xi_n=\xi_{n-1}\circ\zeta_n$, we get:
$$
q_{t,n}\circ j_{t,n}\circ\phi_t=q_{t,n-1}\circ j_{t,n-1}\circ\phi_t\circ\zeta_n=q_{t,n-1}\circ\eta_{t,n}\circ j_{t,n}\circ\phi_t.
$$
Thus $q_{t,n}=q_{t,n-1}\circ\eta_{t,n}$. We have proved:

\REFPROP{compare2}
For $t, n\ge 0$, there exists a quotient map $\eta_{t, n}$ of $\cbar$ such that $q_{t,n}=q_{t,n-1}\circ\eta_{t,n}$.  For a fixed $t\ge 0$, $\{\eta_{t,n}\}$ converges uniformly to the identity as $n\to\infty$.
\ENDPROP

Now we want to control the distortion of $\{q_{t,n}\}$ using the results obtained in \S8.

\begin{proposition}\label{modulus-q}
Let $w\in\cbar$ be a point and let $U\ni w$ be a disk. Then for any $M>0$, there exist a constant $t_0\ge 0$ and a disk $V\Subset U$ with $w\in V$ such that
$$
\md q_{t,n}^{-1}(U\smm\overline{V})\ge M
$$
for $n\ge 0$ and $t\ge t_0$. Moreover, the constant $t_0$ can be chosen to be $t_0=0$ if $w\not\in Y$.
\end{proposition}

\beginp
If $w\in\FFF_g$, then $\{q_{t,n}\}$ is uniformly quasiconformal in a neighborhood of $w$ by Lemmas \ref{factor} and \ref{compare1}. The lemma is trivial in this case. Now we assume $w\in\JJJ_g$.

{\it Case 1}. Assume that the $\omega$-limit set $\omega(w)\cap Y=\emptyset$. By Proposition \ref{disk}, there exists a pullback system $\{U_y\}$ of $Y$ such that $w\notin\overline{U_y}$ for $y\in Y$. By Lemma \ref{prepare}, there exists a constant $r_0\in (0,1)$ such that if

(a) $A_1\subset A_0$ are annuli in $\cbar$ such that $A_1$ is contained essentially in $A_0$ and for any point $y\in Y$, $U_y\cap\partial A_0\neq\emptyset$ implies that $U_y\cap A_1=\emptyset$, and

(b) $q$ is a quotient map of $\cbar$ such that
$$
q^{-1}:\, \cbar\smm\bigcup_{n(x)\le n}\overline{U_y(r_0)}\to\cbar
$$
is univalent for some integer $n\ge 0$, then $\md q^{-1}(A_0)\ge(\md A_1)/2$.

Let $t_0\ge 0$ be a constant such that $\WWW(t_0)\subset\cup_{n(y)=0}U_y(r_0)$. By Proposition \ref{compare1}, for any $0\le t\le t_0$, there is a quasiconformal map $\wp_{t_0,t;\,n}$ such that $q_{t,n}=\wp_{t_0,t;\, n}\circ q_{t_0,n}$, whose maximal dilatation is bounded by a constant $K>1$ depending only on $t_0$.

For any constant $M>0$, since the diameter of $U_y$ tends to zero as $n(y)\to\infty$, there exist a disk $V\Subset U$ with $w\in V$ and an annulus $A_1$ contained essentially in $A:=U\smm\overline{V}$ such that $\md A_1\ge 2KM$ and for any point $y\in Y$, $U_y\cap\partial A\neq\emptyset$ implies that $U_y\cap A_1=\emptyset$. Now applying Lemmas \ref{prepare} and \ref{factor} (3), we have
$$
\md q_{t,n}^{-1}(A_0)\ge\ds\frac{1}{2}\md A_1\ge KM\quad\text{for $n\ge 0$ and $t\ge t_0$}.
$$
It follows that
$$
\md q_{t,n}^{-1}(U\smm\overline{V})\ge M\quad\text{for $t, n\ge 0$}.
$$

{\it Case 2}. Assume $w=y_k\in Y_k$. We may assume further that $y_k=0$ for simplicity. Set $r_1=\inf\{|z|:\, z\in\partial U\}$. For any constant $M>0$, let $\epsilon\in (0, r_1)$ be a constant such that $2\epsilon e^{4\pi M}\le r_1-\epsilon$. Let
$$
V=\{z:\, |z|<\epsilon\}\quad\text{and}\quad A_1=\{z:\, 2\epsilon<|z|<r_1-\epsilon\}.
$$
Then $A_1\subset U\smm\overline{V}$ and $\md A_1\ge 2M$.

Let $\{U_y\}$ be a pullback system of $Y$ such that the Euclidean diameter $\text{\rm diam} U_y\le\epsilon$ for any $y\in Y$. Applying Lemma \ref{prepare} for this pullback system, we obtain a constant $t_0\in (0,1)$ such that for any quotient map $q$ of $\cbar$, if $q^{-1}: \cbar\smm\overline{\WWW_n(t_0)}\to\cbar$
is univalent for some integer $n\ge 0$, then $\md q^{-1}(U\smm\overline{V})\ge(\md A_1)/2$. In particular,
$$
\md q_{t,n}^{-1}(U\smm\overline{V})\ge\ds\frac{1}{2}\md A_1\ge M\quad\text{for $n\ge 0$ and $t\ge t_0$}.
$$

{\it Case 3}. Assume $\omega(w)\cap Y\neq\emptyset$. Then there exists a point $x\in Y_0$ such that $x\in\omega(w)$. Without loss of generality, we may assume that $g(x)=x$. Since $w\in\JJJ_g\smm Y$, $\omega(w)\neq\{x\}$. Thus we may choose a  pullback system $\{U'_y\}$ of $Y$ such that $\omega(w)\smm U'_x\neq\emptyset$.

By Proposition \ref{nested}, there exists a small pullback system $\{U_y\}$ of $Y$ with $U_y\Subset U'_y$ for $y\in Y$ such that for each point $y\in Y_n$, if $\partial U'_{y'}\cap U'_y\neq\emptyset$ and $n(y')>n(y)$, then $U'_{y'}\cap U_y=\emptyset$.

Applying case 2 for the point $x$ and the domain $U_x$, we obtain a constant $r_0\in (0,1)$ and a disk $V_x\Subset U_x$ with $x\in V_x$ such that $\overline{\WWW(t_0)}\subset\cup_{n(y)=0}U_y$ and for any quotient map $q$ of $\cbar$, if
$$
q^{-1}: \cbar\smm\overline{\WWW_n(t_0)}\to\cbar
$$
is univalent for some integer $n\ge 0$, then $\md q^{-1}(U_x\smm\overline{V_x})\ge M$.

Since $x$ is a parabolic fixed point, we may assume that $W\cap\JJJ_g\subset V_x\cap\JJJ_g$, where $W$ is the component of $g^{-1}(V_x)$ containing the point $x$.

Let $\{g^{n_k}(w)\}$ $(k\ge 0)$ be the first returns to $V_x$ of the orbit $\{g^n(w)\}_{n\ge 1}$. It is defined as follows: $n_0\ge 1$ is the minimal positive integer such that $g^{n_0}(w)\in V_x$. Since $\omega(w)\smm U_x\neq\emptyset$, there exists an integer $n>n_0$ such that $g^n(w)\not\in V_x$. Let $n_1>n_0$ be the minimal integer such that $g^{n_1}(w)\in V_x$ but $g^{n_1-1}(w)\not\in V_x$. Inductively, $n_k>n_{k-1}$ is the minimal integer such that $g^{n_k}(w)\in V_x$ but $g^{n_k-1}(w)\not\in V_x$.

Recall that $w\in\JJJ_g$ and $W\cap\JJJ_g\subset V_x\cap\JJJ_g$. Thus for $k\ge 1$, we have $g^{n_k-1}(w)\notin W$ and hence $g^{n_k-1}(w)\in U_{y_1}$ for a certain point $y_1\in Y_1$ with $g(y_1)=x$.

For each $k\ge 1$, let $V_k$ and $U_{y_k}$ be the components of $g^{-n_k}(V_x)$ and $g^{-n_k}(U_x)$ containing the point $w$, respectively. Then
$g^{n_k}:\ U_{y_k}\to U_x$ is conformal by the assumption at the beginning of this section.

Given any integer $n\ge 0$, let $q$ be a quotient map of $\cbar$ such that $q^{-1}$ is univalent in $\cbar\smm\overline{\WWW_n(t_0)}$. For each component $W$ of $\overline{\WWW_n(t_0)}$ with $W\cap\partial U_{y_k}=\emptyset$, $g^{n_k}(W)$ is also a component of $\overline{\WWW_n(t_0)}$ which is not periodic since $g^{n_k}$ is the first return map. Thus there exists a disk $U_k\Subset U'_{y_k}$ with $U_{y_k}\subset U_k$ such that if $W$ is a component of $\overline{\WWW_n(t_0)}$ with $W\cap\partial U_{y_k}=\emptyset$, then $W\subset U_k$. Define $h_k$ to be equal to $g^{n_k}\circ q$ on $U_k$ and quasiconformal otherwise. Then $h_k$ is also a quotient map of $\cbar$. By the Riemann Mapping Theorem, there exists a quasiconformal map $\phi$ of $\cbar$ which is holomorphic in $U_k$ such that $\phi\circ h_k^{-1}$ is univalent in
$$
\cbar\smm g^{n_k}(\overline{\WWW_n(t_0)}\cap U_k)\subset\cbar\smm\overline{\WWW_n(t_0)}.
$$
Thus
$$
\md \phi\circ h_k^{-1}(U_x\smm\overline{V_x})=\md\phi\circ q^{-1}(U_{y_k}\smm\overline{V_k})\ge M.
$$
Note that $\phi$ is conformal in $U_k\supset U_{y_k}$. Therefore $\md q^{-1}(U_{y_k}\smm\overline{V_k})\ge M$. In particular,
$$
\md q_{t,n}^{-1}(U_{y_k}\smm\overline{V_k})\ge M\quad\text{for $n\ge 0$ and $t\ge t_0$}.
$$

By Proposition \ref{compare1}, there exists a constant $K>1$ such that
$$
\md q_{t,n}^{-1}(U_{y_k}\smm\overline{V_k})\ge\frac{M}{K}
$$
for $t, n\ge 0$. Choose an integer $k_0\ge 1$ such that $k_0/K>1$. Since the diameter of $U_{y_k}$ converges to zero as $k\to\infty$, there are $k_0$
annuli $U_{y_k}\smm\overline{V_k}$ contained in $U$ that are pairwise disjoint. Set $V$ to be the smallest one of these $k_0$ domains $V_k$. Then
$$
\md q_{t,n}^{-1}(U\smm\overline{V})\ge M\quad\text{for $t,n\ge 0$.}
$$
\qed

\REFLEM{normal-q}
The family $\{q_{t,n}\}$ $(t, n\ge 0)$ is equicontinuous.
\ENDLEM

\beginp We claim that for any point $w_0\in\cbar$ and any disk $U\subset\cbar$ with $w_0\in U$, there exists a constant $\delta(w_0)>0$ and a disk $V\ni w_0$ with $\overline{V}\subset U$, such that for any $q_{t,n}$, the spherical distance
$$
\text{\rm dist}_s(q_{t,n}^{-1}(\partial U), q^{-1}(V))>\delta(w_0).
$$

The claim holds for $w_0\in\FFF_g$ by Lemma \ref{distortion2}, Proposition \ref{compare1} and the normalization condition. Now we assume $w_0\in\JJJ_g$.

Assume $\infty\in\FFF_g$ for simplicity. Choose a constant $M>5\log 2/(2\pi)$. Then there exists a constant $t_0\ge 0$ and a disk $V\ni w_0$ with $\overline{V}\subset U$, such that
$$
\md q_{t,n}^{-1}(U\smm\overline{V})\ge M
$$
for $n\ge 0$ and $t\ge t_0$, by Proposition \ref{modulus-q}. From Lemma \ref{mod1}, we have
$$
\text{\rm dist}(q_{t,n}^{-1}(\partial U), q_{t,n}^{-1}(V))>C(M)\text{diam}\,q_{t,n}^{-1}({\overline V})
$$
for all $n\ge 0$ and $t\ge t_0$. On the other hand, there exists a disk $D\subset V$ such that $D\Subset\FFF_g$. Thus there exists a constant $t_1\ge t_0$ such that $D$ is disjoint from $g^{-n}(\WWW_t)$ for $n\ge 0$ and $t\ge t_1$. So there exists a constant $C>0$ such that
$$
\text{diam}\,q_{t,n}^{-1}({\overline D})\ge C\text{ for $n\ge 0$ and $t\ge t_1$}.
$$
Take $\delta_1(w_0)=C\cdot C(M)$. Then the claim holds for $t\ge t_1$.

By Proposition \ref{compare1}, there exists a normalized quasiconformal map $\wp_{t_1, t;\, n}$ of $\cbar$ whose maximal dilatation is equal to the maximal dilatation of $\phi_{t_1}\circ\phi_t$, such that $q_{t_1, n}\circ\wp_{t_1, t;\, n}=q_{t,n}$. Thus $\wp_{t_1, t;\, n}^{-1}$ is uniformly H\"{o}lder continuous for all $n\ge 0$ and $t\in (0,t_1)$. In particular, there exists a constant $\delta(w_0)\le\delta_1(w_0)$ such that
$$
\text{\rm dist}(\wp_{t_1, t;\, n}^{-1}(z_1), \wp_{t_1, t;\, n}^{-1}(z_2))\ge\delta(w_0)\text{ if }
\text{\rm dist}(z_1, z_2)\ge\delta_1(w_0)
$$
for all $n\ge 0$ and $t\in (0,t_1)$. Now the claim is proved. By Lemma \ref{equi-p}, the family $\{q_{t,n}\}$ is equicontinuous. \qed

\REFLEM{limit-q}
For fixed $t\ge 0$, $\{q_{t,n}\}$ converges uniformly to a quotient map $q_t$ of $\cbar$ as $n\to\infty$, and $q_t\circ f_t=g\circ q_t$. For each filled-in skeleton $\wh S$, $q_t\circ\phi_t(\wh S)$ is a single point. Conversely, for each point $w\in\cbar$, $(q_t\circ\phi_t)^{-1}(w)$ is a single point if $w\notin Y$, or else is a filled-in skeleton of level $n$ if $w\in Y_n$.
\ENDLEM

\beginp Let $q_t$ be the limit of a subsequence $\{q_{t,n_k}\}$ as $n_k\to\infty$. Then the subsequence $\{\xi_{t,n_k}=q_{t,n_k}\circ j_{t,n_k}\}$ also converges uniformly to $q_t$ as $n_k\to\infty$ since the sequence $\{j_{t,n}\}$ converges uniformly to the identity as $n\to\infty$ by Lemma \ref{g2F}. By Proposition \ref{compare2}, $q_{t,n}=q_{t,n-1}\circ\eta_{t,n}$ where $\{\eta_{t,n}\}$ converges uniformly to the identity. Thus the sequence $\{\xi_{t,n_k-1}=q_{t,n_k-1}\circ j_{t,n_k-1}\}$ also converges uniformly to $q_t$ as $n_k\to\infty$. From $\xi_{t,n-1}\circ f_t=g\circ\xi_{t,n}$, we get $q_t\circ f_t=g\circ q_t$.

Let $w\in Y_0$ be a periodic point with period $p\ge 1$. Then $\xi_{t,n}^{-1}(w)$ is the $\phi_t$-image of a periodic filled-in skeleton $\wh S$, by Lemma \ref{f2g}. Thus $q_t(\phi_t(\wh S))=w$ since $\{j_{t,n}\}$ converges uniformly to the identity as $n\to\infty$. It follows that $\phi_t(\wh S)\subset q_t^{-1}(w)$.

Consider the continuum $q_t^{-1}(w)$. For any point $z\in\FFF_{f_t}$ with $q_t(z)=w$, there exists an integer $N>0$ such that $f_t^n(z)\in\phi_t(\UUU)$ for $n\ge N$, where $\UUU$ is a fundamental set of $f$ as defined in Lemma \ref{f2g}. Thus $\xi_{t,n}(z)=q_t(z)=w$ for $n\ge N$ by Lemma \ref{f2g}. Therefore  $q_t^{-1}(w)\cap\FFF_{f_t}=\phi_t(\wh S)\cap\FFF_{f_t}$.

Note that $E:=\phi_t(\wh S)\cap\JJJ_{f_t}$ contains only finitely many points, each of them periodic. There are finitely many disjoint disks in $\FFF_{f_t}$ such that their union $D$ contains $\phi_t(\wh S)\cap\FFF_{f_t}$ and $\partial D\cap\phi_t(\wh S)=E$.  Thus $\partial D\cap q_t^{-1}(w)=E$. We claim that each component of $K:=q_t^{-1}(w)\cap\JJJ_{f_t}$ contains a point of $E$. Otherwise, assume that $W_0, W_1\subset\cbar$ are disjoint domains that both contain points of $K$, $K\subset W_0\cup W_1$ and $W_0\cap E=\emptyset$. Then $W_0\smm\overline{D}$ has a component $W'$ such that $W'\cap K\neq\emptyset$. Moreover, $W'\cap K$ is compact since $\partial W'\subset\partial W_0\cup\partial D$ has positive distance from $W_0\cap K$. On the other hand, $W'\cap q_t^{-1}(w)=W'\cap K$ since $q_t^{-1}(w)\cap\FFF_{f_t}=\phi_t(\wh S)\cap\FFF_{f_t}\subset D$ and $W'\cap D=\emptyset$. It follows that $q_t^{-1}(w)$ is disconnected. This is a contradiction.

By the claim, $q_t^{-1}(w)\cap\JJJ_{f_t}$ has only finitely many components and each of them is eventually periodic. From the equation $q_t\circ f_t=g\circ q_t$, we see that $q_t^{-1}(w)$ is a component of $f_t^{-p}(q_t^{-1}(w))$ and $f_t^p$ is univalent in a neighborhood of $q_t^{-1}(w)$ since $g^p(w)=w$ and $\deg_w g^p=1$.
Thus each component of $q_t^{-1}(w)\cap\JJJ_{f_t}$ must be a single point by Lemma \ref{shrinking}. Therefore $q_t^{-1}(w)\cap\JJJ_{f_t}=\phi_t(\wh S)\cap\JJJ_{f_t}$ and hence $q_t^{-1}(w)=\phi_t(\wh S)$.

Let $w\in\cbar\smm Y$ be a point. For contradiction we assume that $q_{t}^{-1}(w)$ is not a single point. Let $U\subset\cbar$ be a disk with $w\in U$. Then
$$
\md q_t^{-1}(U\smm\{w\})=M<\infty.
$$
Pick a disk $U_1\ni w$ such that $\overline{U_1}\subset U$. From Proposition \ref{modulus-q}, there exists a disk $V\ni w$ with $\overline{V}\subset U_1$, such that
$$
\md q_{t,n}^{-1}(U_1\smm\overline{V})> M\quad\text{for }t, n\ge 0.
$$
Since $\{q_{t,n_k}\}$ converges uniformly to $q_t$ as $n_k\to\infty$, there exists an integer $N>0$ such that $q_{t,n_k}^{-1}(U_1)\subset q_t^{-1}(U)$ and $q_{t,n_k}^{-1}(V)\supset q_t^{-1}(w)$ for $n_k\ge N$. Thus
$$
\md q_{t,n}^{-1}(U_1\smm\overline{V})> M=\md q^{-1}(U\smm\{w\})>\md q_{t,n}^{-1}(U_1\smm\overline{V}).
$$
This is a contradiction.

If $\wt q_t$ is another limit of the sequence $\{q_{t,n}\}$, then $\theta:=\wt q_t\circ q_t^{-1}$ is a well-defined homeomorphism of $\cbar$ and $g\circ\theta=\theta\circ g$. Moreover, $\theta$ is holomorphic in $\FFF_g$. Thus $\theta$ is holomorphic on $\cbar$ by Theorem \ref{unicity} and hence is the identity by the normalization condition. Therefore $\wt q_t=q_t$ and hence the whole sequence $\{q_{t,n}\}$ converges uniformly to $q_t$.
\qed

\vskip 0.24cm

{\noindent\it Proof of Lemma \ref{ft2g}}. Parts (1) and (2) are direct consequences of Lemmas \ref{factor}, \ref{g2F} and \ref{limit-q}. Part (3) comes from Lemma \ref{distortion2} and Corollary \ref{identity}.
\qed

\subsection{Proof of Theorem \ref{boundedness}}

The proof of Theorem \ref{boundedness} is a small modification of the proof of Theorem \ref{pinching}. We continue to use the notation from the proof of Theorem \ref{pinching}.

\vskip 0.24cm

{\noindent\it Proof of Theorem \ref{boundedness}}.
Let $g$ be the limit of the pinching path $\{f_t=\phi_t\circ f\circ\phi_t^{-1}\}$ supported on $\sA$. Let $\varphi$ be the limit of $\{\phi_t\}$. Denote by $\WWW$ the interior of $\varphi(\BBB_0)$.

Let $D\Subset\cbar\smm\overline{\cup_{n\ge 0}\WWW_n}$ be a disk. Then $g^{-1}(D)\Subset \cbar\smm\overline{\cup_{n\ge 0}\WWW_n}$. From Theorem \ref{distortion1} and Lemma \ref{factor}, there exist constants $\delta<\infty$ and $t_0\ge 0$ such that for any $n\ge 0$ and any univalent map
$\phi:\, \cbar\smm\overline{\WWW_n(t_0)}\to\cbar$,
$$
\sD_0(\phi, D\cup g^{-1}(D))\le\delta<\infty.
$$

Define $\fm[g,D,\delta]\subset\fm_d$ ($d=\deg g$) by $[h]\in \fm[g,D,\delta]$ if there exists a univalent map $\phi: D\cup g^{-1}(D)\to\cbar$ such that $h\circ\phi=\phi\circ g$ on $g^{-1}(D)$ and $\sD_0(\phi, D\cup g^{-1}(D))\le\delta$. Since $\deg(g: g^{-1}(D)\to D)=d$, $\fm[g,D, \delta]$ is compactly contained in $\fm_d$ by Lemma \ref{rational-s}.

Let $\Theta: \sR_f\to\sR_{f_{t_0}}$ be the quasiconformal map whose lift is $\phi_{t_0}$. Let $\nu$ be a Beltrami differential supported on $\Theta(\sA(t_0))$ with $\|\nu\|_{\infty}<1$. Let $\phi$ be a quasiconformal map of $\cbar$ whose Beltrami differential is the lift of $\nu$. Then $h=\phi\circ f_{t_0}\circ\phi^{-1}$ is a rational map. Repeat the proof of Theorem \ref{pinching} with $\xi_{t,n}$ replaced by $\xi_{t_0,n}\circ\phi$. We again obtain a sequence of quotient maps $\{q_n\}$ and a Thurston sequence $\{h_n\}$ of $h$ such that $q_n\circ h_n=g\circ q_{n+1}$ and $q_n^{-1}$ is holomorphic and injective in $\cbar\smm\overline{\WWW_n(t_0)}$. Thus $[h]\in\fm[g,D,\delta]$.

From the definition of a multi-annulus, there exists a quasiconformal map $\Phi$ from $\sR_f$ to itself such that $\Phi(\sA(t_0))\subset\sA$. Then $\Theta\circ\Phi^{-1}(\sA)\subset\Theta(\sA(t_0))$.

For any Beltrami differential $\mu$ on $\sR_f$ supported on $\sA$, let $\Phi_{\mu}: \sR(f)\to\sR(f_{\tilde\mu})$ be the quasiconformal map with Beltrami differential $\mu$. Let
$$
\Psi=\Phi_{\mu}\circ\Phi\circ\Theta^{-1}:\ \sR(f_{t_0})\to\sR(f_{\tilde\mu}).
$$
Then $\mu(\Psi)=\mu(\Phi\circ\Theta^{-1})$ on $\sR_{f_{t_0}}\smm\Theta(\sA(t_0))$ and hence $\mu(\Psi)$ is independent of the choice of $\mu$. Let $\Psi_0$ be a quasiconformal map on $\sR(f_{t_0})$ such that $\mu(\Psi_0)=\mu(\Phi)$ on $\Theta(\sA(t_0))$ and $\mu(\Psi_0)=0$ otherwise. Set $\Psi_1=\Psi\circ\Psi_0^{-1}$. Then the maximal dilatation of $\Psi_1$ depends only on $\Phi\circ\Theta^{-1}$ and hence is bounded by a constant $K<\infty$.

Let $h$ be the quasiconformal deformation of $f_{t_0}$ with Beltrami differential $\mu(\Psi_0)$; then $h\in \fm[g,D,\delta]$ and $h$ is $K$-quasiconformally conjugate to $f_{\tilde\mu}$. Define $\fm[g,D,\delta;K]\subset\fm_d$ by $[h_1]\in\fm[g,D,\delta; K]$ if $h_1$ is $K$-quasiconformally conjugate to a rational map $h$ with $[h]\in\fm[g,D,\delta]$. Then $[f_{\tilde\mu}]\in\fm[g,D,\delta;K]$. Obviously $\fm[g,D,\delta;K]$ is compactly contained in $\fm_d$. It is easy to check that each rational map in the closure of $\fm[f,\sA]$ is geometrically finite. \qed

\section{Parabolic-hyperbolic deformation}

Let $g$ be a geometrically finite rational map with parabolic cycles. Let $\fw$ be a collection of a pairwise disjoint sepals such that the closure of their union $\WWW$ is disjoint from $\PPP_g\smm\PPP'_g$ and $g(\WWW)=\WWW$. A bijection $\sigma: \fw\to\fw$ is called a {\bf plumbing correspondence}\/ if it satisfies the following conditions:

(1) $\sigma^2=\text{id}$ but $\sigma$ has no fixed element.

(2) $\sigma$ is compatible with $g$, i.e. $\sigma (g(W_i))=g(\sigma(W_i))$ for each sepal $W_i$.

(3) If $\sigma(W_i)=W_j$, then $W_i$ and $W_j$ touch each other and $W_j$ is a left sepal if and
only if $W_i$ is a right sepal.

(4) (non-crossing condition) If $W_i$ touches $W_j$ at $y\in\PPP'_g$ but $\sigma(W_i)\neq W_j$, then both $W_j$ and $\sigma(W_j)$ lie on the same side of $W_i\cup\{y\}\cup\sigma(W_i)$.

Two plumbing correspondences $\sigma: \fw\to\fw$ and $\sigma': \fw'\to\fw'$ are called {\bf equivalent}\ if there is a plumbing correspondence $\sigma'': \fw''\to\fw''$ such that all of them have same number of sepals and for each sepal $W''\in\fw''$, there are sepals $W\in\fw$ and $W'\in\fw'$ such that $W''\in W\cap W'$ and
$\sigma''(W'')\subset\sigma(W)\cap\sigma'(W')$. A {\bf plumbing combinatorics} is an equivalence class of plumbing correspondences.

\vskip 0.24cm

In this section, we will prove the following theorem -- a precise version of Theorem \ref{plumbing}.

\REFTHM{plumbing-1}
Let $g$ be a geometrically finite rational map and let $Y$ be a set of parabolic cycles of $g$. Let $\fw=\{W_1, \cdots, W_{2m}\}$ be a collection of pairwise disjoint sepals of cycles in $Y$ such that the closure of their union $\WWW$ is disjoint from $\PPP_g\smm Y$ and $g(\WWW)=\WWW$.  Let $\sigma: \fw\to\fw$ be a plumbing correspondence. Then there exist a geometrically finite rational map $f$ and a non-separating multi-annulus $\sA\subset\sR_f$ such that the following conditions hold:

(1) The pinching path $f_t=\phi_t\circ f\circ\phi_t^{-1}$ ($t\ge 0$) starting from $f=f_0$ supported on $\sA$ converges uniformly to the rational map $g$.

(2) Let $\varphi$ be the limit of the conjugacy $\phi_t$ as $t\to\infty$. Let $\BBB_0$ and $\SSS_0$ be the unions of all periodic bands and skeletons, respectively; then $\varphi(\BBB_0\smm{\SSS_0})=\WWW$ and for any two distinct components $B_1, B_2$ of $\BBB_0\smm{\SSS_0}$, $\varphi(B_1)=\sigma(\varphi(B_2))$ if and only if $B_1$ and $B_2$ are contained in the same periodic band.
\ENDTHM

{\bf Step 1. Plumbing surgery}.
Note that the quotient space $\WWW/\langle g\rangle$ is a finite disjoint union of once-punctured disks. Thus there is a natural holomorphic projection $\pi: \WWW\to\DD^*$ such that for each sepal $W$ of $\WWW$ with period $p\ge 1$, $\pi: W\to\DD^*$ is a universal covering and $\pi(z_1)=z_2$ if and only if $z_1=g^{kp}(z_2)$ for some integer $k\in\ZZ$.

Given any $0<r<1$, let $\WWW(r)=\pi^{-1}(\DD^*(r))$ and $\RRR(r)=\WWW\smm\overline{\WWW(r)}$. Then $g(\RRR(r))=\RRR(r)$. Thus there is a conformal map $\tau:\ \RRR(r^2)\to\RRR(r^2)$ such that:

$\bullet$ $\tau(z)\in\sigma(W_i)$ if $z\in W_i$, and

$\bullet$ $\tau^2=\text{id}$ and $g\circ\tau=\tau\circ g$.

Define an equivalence relation on $\cbar\smm\overline{\WWW(r^2)}$ by $z_1\sim z_2$ if $\tau(z_1)=z_2$. Then the quotient space is a punctured sphere with finitely many punctures. Thus there exists a finite set $X\subset\cbar$ and a holomorphic onto map
$$
p:\, \cbar\smm\overline{\WWW(r^2)}\to\cbar\smm X
$$
such that $p(z_1)=p(z_2)$ if and only if $z_1=\tau(z_2)$.

Let $\SSS=p(\partial\WWW(r)\smm Y)\cup X$. This is a finite disjoint union of trees whose vertex set is $X$. Let $\BBB=p(\RRR(r^2))$. This is a finite disjoint union of disks.

\vskip 0.24cm

{\bf Step 2. The induced map after surgery}.
Let $\WWW_1=g^{-1}(\WWW)\smm\WWW$. Since $g\circ\tau=\tau\circ g$, there is a unique holomorphic map $F_0:\,\cbar\smm p(\overline{\WWW_1})\to\cbar$, such that $F_0\circ p=p\circ g$ on $\cbar\smm g^{-1}(\overline\WWW)$. Obviously, all the sets $X$, $\SSS$ and $\BBB$ are fixed by $F_0$.

Define $X_c\subset X$ by $x\in X_c$ if $x$ is an accumulation point of $p(\PPP_g\smm Y)$. Then $F_0(X_c)=X_c$. For each component $B$ of $\BBB$, $\overline{B}\cap X$ contains exactly two points. Let $p\ge 1$ be the period of $B$. Then $\{F_0^{kp}(z)\}$ converges to a point in $\overline{B}\cap X$ for any point $z\in B$.  Denote the point by $a(B)$. Since each attracting petal of $G$ at a point $y\in Y$ contains infinitely many points of $\PPP_G$, we have $a(B)\in X_c$. Denote the other point of $\overline{B}\cap X$ by $r(B)$.

Pick an attracting flower of $g$ at each point $y\in Y$ whose union $\VVV$ satisfies the following conditions:

(1) $\overline{g(\VVV)}\subset\VVV\cup Y\subset\cbar\smm\overline{\WWW_1}$, and

(2) each component $R$ of $\RRR(r^2)$ intersects $\partial\VVV$ at exactly two open arcs and $\tau(\RRR(r^2)\cap\partial\VVV)=\RRR(r^2)\cap\partial\VVV$.

Then each component of $p(\partial\VVV\smm\overline{\WWW(r^2)})$ is either a Jordan curve, or an open arc whose two endpoints land on the same point in $X$. Denote by $\UUU$ the union of disks enclosed by these closed curves and open arcs together with their endpoints. Then $\overline{F_0(\UUU)}\subset\UUU\cup X$.

Let $x\in X_c$ be a point. If $x\in\UUU$, then $x$ is an attracting point of $F_0$. Otherwise, let $\VVV_x$ be the union of the components of $\UUU$ touching the point $x$; then $\VVV_x$ satisfies the conditions of Lemma \ref{parabolic}. Thus the point $x$ is a parabolic point of $F_0$ and $\VVV_x$ is an attracting flower of $F_0$ at the point $x$. Obviously, each component of $\VVV_x$ contains infinitely many points of $p(\PPP_g\smm Y)$. The following proposition is easy to verify.

\begin{proposition}\label{infty-plumbing} Let $S$ be a component of $\SSS$ and let $x\in S\cap X_c$ be a point. Let $k\ge 1$ be the number of components of $S\smm\{x\}$. Let $D_x\ni x$ be a sufficiently small disk such that $D_x\smm S$ has $k$ components $U_i$ whose closures contain the point $x$. Then $U_i$ contains infinitely many points of $p(\PPP_g)$ if there exists a component $B$ of $\BBB$ such that $a(B)=x$ and $U_i\cap B\neq\emptyset$.
\end{proposition}

\vskip 0.24cm

{\bf Step 3. Quotient extension of the inverse map of the projection}. Note that the map $p: \cbar\smm\overline{\WWW}\to\cbar\smm\overline{\BBB}$ is a conformal map. We want to extend its inverse map to be a quotient map of $\cbar$ as follows. Let $w:\, \AA(r,1)\to\DD^*$ be the homeomorphism defined in Proposition \ref{model} (4). Then there exists a unique homeomorphism $\wt w:\, \RRR(r)\to\WWW$ such that $\pi\circ\wt w=w\circ\pi$, $g\circ\wt w=\wt w\circ g$ and the continuous extension of $\wt w$ to $\partial\WWW\subset\partial\RRR(r)\to\partial\WWW$ is the identity. Define
$$
q=\begin{cases}
p^{-1}:\, \cbar\smm\overline{\BBB}\to\cbar\smm\overline{\WWW} \\
{\wt w}\circ p^{-1}:\, \BBB\smm\SSS\to\RRR(r)\to\WWW.
\end{cases}
$$
Then $q: \cbar\smm\SSS\to\cbar\smm Y$ is a homeomorphism and hence can be extended to a quotient map of $\cbar$ with $q(S)=Y$. Since $F_0\circ p=p\circ g$ on $\cbar\smm\overline{\WWW_1}$ and $g\circ\wt w=\wt w\circ g$, we have $g\circ q=q\circ F_0$ on $\cbar\smm q^{-1}(\overline{\WWW_1}))$.

\vskip 0.24cm

{\bf Step 4. Construction of a marked semi-rational map}. Each component $E$ of $\overline{\WWW}$ is a full continuum which contains exactly one point of $\PPP_g$. Pick a disk $U(E)\supset E$ such that $U(E)\smm E$ contains no critical values of $g$ and $\partial U(E)$ is disjoint from $\PPP_g$. We may assume that all these domains $U(E)$ have disjoint closures. Denote by $U_0$ their union. Note that $U(E)$ contains at most one critical value of $g$. Each component of $g^{-1}(U_0)$ is a disk containing exactly one component of $g^{-1}(\overline{\WWW})$.

Let $U_1$ be the union of components of $g^{-1}(U_0)$ which contain a component of $\overline{\WWW_1}$. Since $\overline{\WWW_1}$ is disjoint from $\PPP'_g$, Once $U_0$ is close enough to $\WWW$, we may assume $U_1\smm\overline{\WWW_1}$ is disjoint from $\PPP_g$.

\vskip 0.24cm

Define a branched covering $F$ of $\cbar$ such that:

(a) $F(z)=F_0(z)$ on $\cbar\smm q^{-1}(U_1)$ and hence $g\circ q=q\circ F$ on $\cbar\smm q^{-1}(U_1)$.

(b) $F: q^{-1}(U_1)\to q^{-1}(U_0)$ is a branched covering with at most one critical point in $q^{-1}(Y_1)$ and $F(q^{-1}(Y_1))\subset X$, where $Y_1=g^{-1}(Y)\smm Y$.

It follows that:
$$
q^{-1}(\PPP_g\smm Y)\cup X_a\subset\PPP_F\subset q^{-1}(\PPP_g\smm Y)\cup X,
$$
and $\PPP'_F=q^{-1}(\PPP'_g\smm Y)\cup X_a$. In particular, $q(\PPP_F)=\PPP_g$. Set $\PPP=\PPP_f\cup X$. Then $q(\PPP)=\PPP_g$, $F(\PPP)=\PPP$ and $\#(\PPP\smm\PPP_F)<\infty$.

Since $F$ is holomorphic in $\cbar\smm q^{-1}(\overline{U_1})$ and each component of $q^{-1}(U_1)$ contains at most one point of $\PPP_F$,  which is an isolated point of $\PPP_F$,  we know that $F$ is holomorphic in a neighborhood of $\PPP'_F$. By Step 2, $(F,\PPP)$ is a marked semi-rational map.

\vskip 0.24cm

{\bf Step 5. Lift of the quotient map}.
For each component $D$ of $U_1$, $g: D\to g(D)$ is proper with at most one critical value in $Y$. On the other hand, $F: q^{-1}(D)\to F(q^{-1}(D))$ is also a branched covering with at most one critical value in $\SSS$, and $\deg g|_D=\deg F|_{q^{-1}(D)}$. Since $q(\SSS)=Y$, there exists a quotient map $\wt q: q^{-1}(D)\to D$ that coincides with $q$ on the boundary such that $g\circ\wt q=q\circ F$ on $q^{-1}(D)$.

Define $\wt q=q$ on $\cbar\smm q^{-1}(U_1)$. Then $\wt q$ is a quotient map of $\cbar$ isotopic to $q$ rel $(\cbar\smm q^{-1}(U_1))\cup\PPP$ and $g\circ\wt q=q\circ F$ on $\cbar$.

\vskip 0.24cm

\REFLEM{fundamental} If $\UUU$ is a fundamental set of $F$, then $q(\UUU)$ contains a fundamental set of $g$.
\ENDLEM

\beginp We only need to prove that $q(\UUU)$ contains an attracting flower of $g$ at each point $y\in Y$. Let $S=q^{-1}(y)$. For simplicity of notation we assume that each point in $S\cap X$ is fixed by $F$.

Let $x\in S\cap X_a$ be a point. If it is attracting, then $x\in\UUU$. Thus there exists a disk $D_x\subset\UUU$ with $x\in D$ such that $f$ is injective on $D_x$ and $f(D_x)\Subset D_x$. Moreover, we may require that $\partial D_x$ intersects each component of $S\smm\{x\}$ at a single point. Then for each component $U$ of $D_x\smm S$,  $V=q(U)$ is a disk and $\overline{g(V)}\subset V\cup\{y\}$.

Now suppose that $x$ is parabolic. Then there exists an attracting flower $\VVV_x$ of $F$ at $x$ such that $\VVV_x\subset\UUU$. We may also require that each component of $\partial\VVV_x\smm\{x\}$ is either disjoint from $S$ or intersects with each component of $S\smm\{x\}$ at a single point. Then for each component $U$ of $\VVV_x\smm S$,  $V=q(U)$ is a disk and $\overline{g(V)}\subset V\cup\{y\}$.

Denote by $\VVV_1$ the union of $V=q(U)$ for all components $U$ of $D_x\smm S$ if $x\in S\cap X_c$ is attracting and for all components $U$ of $\VVV_x\smm S$ if $x\in S\cap X_c$ is parabolic. If $\{w_n=g^{n}(w)\}$ is an orbit in $\cbar\smm U_1$ converging to the point $y$ as $n\to\infty$ but $w_n\neq y$ for all $n\ge 1$, then $\{z_n=q^{-1}(w_n)\}$ converges to a point $x\in X_a$ and $F(z_n)=z_{n+1}$. Thus once $n$ is large enough, the point $z_n$ is contained in either $D_x$ if $x$ is attracting or $\VVV_x$ if $x$ is parabolic. So $w_n\in\VVV_1$ once $n$ is large enough. Therefore $\VVV_1$ is an attracting flower of $g$ at $y$.
\qed

\REFLEM{no-obs-3}
The marked semi-rational map $(F, \PPP)$ has neither Thurston obstructions nor connecting arcs.
\ENDLEM

\beginp
The proof of this theorem is similar to the proof of Theorem \ref{no-obs1}.

Assume for contradiction that $\Gamma$ is an irreducible multicurve of $(F, \PPP)$ with $\lambda(\Gamma)\ge 1$. Assume futher that for each $\gamma\in\Gamma$, $\#(\gamma\cap\SSS)$ is minimal in its isotopy class. Since $F: \SSS\to\SSS$ is bijective, $k=\#(\gamma\cap\SSS)<\infty$ is a constant for $\gamma\in\Gamma$.

If $k=0$, then $\Gamma_1=\{q(\gamma):\, \gamma\in\Gamma\}$ is a multicurve of $g$ since $q: \PPP\smm X\to\PPP_g\smm Y$ is injective and $q(X)=Y\in\PPP'_g$. Noticing that $g\circ\wt q=q\circ F$ on $\cbar$ and $\wt q$ is a quotient map of $\cbar$ isotopic to $q$ rel $\PPP$, we have $\lambda(\Gamma)=\lambda(\Gamma_1)<1$. This is a contradiction.

Now we assume that $k>0$. Then there exists at most one component of $F^{-1}(\gamma)$ isotopic to a curve in $\Gamma$ rel $\PPP$ for each $\gamma\in\Gamma$ since $F: \SSS\to\SSS$ is bijective. Thus for each $\gamma\in\Gamma$, there is exactly one curve $\beta\in\Gamma$ such that $F^{-1}(\beta)$ has a component isotopic to $\gamma$ rel $\PPP$, since $\Gamma$ is irreducible. Therefore each entry of the transition matrix $M(\Gamma)$ is less than or equal to $1$. Because $\lambda(\Gamma)\ge 1$, there is a curve $\gamma\in\Gamma$ such that $\gamma$ is isotopic to a component $\delta$ of $F^{-p}(\gamma)$ rel $\PPP$ for some integer $p\ge 1$, and $F^p$ is injective on $\delta$.

Let $\UUU$ be a fundamental set of $F$ that is disjoint from every curve in $\Gamma$. Since $q(\UUU)\smm Y$ contains a fundamental set of $g$ and $q$ is injective on $\cbar\smm\SSS$, $q(\gamma)$ is disjoint from a fundamental set of $g$.

Suppose $\gamma$ intersects at least two components of $\SSS$. Let $\beta$ be a component of $q(\gamma)\smm Y$ such that $\beta$ joins two distinct points in $Y$. Then $\beta$ is isotopic to a component of $g^{-kp}(\beta)$ rel $\PPP_g$ for some integer $k>0$ since $\#(\gamma\cap\SSS)$ is minimal in its isotopy class and $\gamma$ is isotopic to a component of $F^{-p}(\gamma)$ rel $\PPP$. Thus $\beta$ is a connecting arc of $g$. This is a contradiction.

Suppose that $\gamma$ intersects exactly one component $S$ of $\SSS$. We claim that at least two components of $\cbar\smm(\gamma\cup S)$ contain points of $\PPP\smm X$.

Let $V_0, V_1$ be the components of $\cbar\smm\gamma$. If both of them contain points of $\PPP\smm X$, then both of them contain a component of $\cbar\smm(\gamma\cup S)$ which contains points of $\PPP\smm X$. The claim is proved. Now we assume that one of them, say $V_0$, contains no points of $\PPP\smm X$. Since $\gamma$ is non-peripheral, $V_0$ contains at least two points of $\PPP$. Thus $V_0$ contains two distinct points $x_0, x_1\in X$. Since $X_c\subset\PPP'_F$, we have $x_1, x_2\in X\smm X_c$.

As $S$ is a tree, there exists a unique arc $l\subset S$ whose endpoints are $(x_1, x_2)$. Let $B_0,  B_1$ be the components of $\BBB$ intersecting $l$ such that $x_i\in\overline{B_i}$. Then $x_0=r(B_0)$ and $x_1=r(B_1)$ since $a(B)\in X_c$ for each component $B$ of $\BBB$. Consequently, there exists a point $x_2\in X_c\cap l$. Now $V_1\smm l$ has exactly two components $U_0, U_1$ whose closures contain the point $x_2$. By Proposition \ref{infty-plumbing}, each of them contains infinitely many points of $\PPP$. Therefore, $V_1\smm S$ has at least two components $U'_0, U'_1$ which contain infinitely many points of $\PPP$. The claim is proved.

Let $y=q(S)$. Then there exists a component $\beta$ of $q(\gamma)\smm\{y\}$ such that $\beta\cup\{y\}$ separates $q(U'_0)$ from $q(U'_1)$. In other words, each component of $\cbar\smm(\beta\cup\{y\})$ contains at least one point of $\QQQ$ since $q(\PPP)=\QQQ$. As above, $\beta$ is isotopic to a component of $G^{-kp}(\beta)$ rel $\QQQ$ for some integer $k>0$. Thus $\beta$ is a connecting arc of $(G, \QQQ)$. This is a contradiction. Thus $(F, \PPP)$ has no Thurston obstructions.

Suppose that $\beta$ is a connecting arc of $(F, \PPP)$. We may assume that $k=\#(\beta\cap\SSS)$ is minimal in the isotopy class of $\beta$. If $k=0$, then $q(\beta)$ is a connecting arc of $G$. This ia a contradiction. Otherwise, by discussion similar to the one as above, there exists a component $\delta$ of $\beta\smm\SSS$ such that $q(\delta)$ is a connecting arc of $G$. This ia a contradiction. Therefore $(F,\PPP)$ has no connecting arc.
\qed

\vskip 0.24cm

{\noindent\it Proof of Theorems \ref{plumbing-1} and \ref{plumbing}}. With the previous preparation, the proof is the same as the proof of Theorem \ref{s-plumbing-1}.  We omit it here.
\qed

\noindent
Guizhen Cui \\
Academy of Mathematics and Systems Science \\
Chinese Academy of Sciences \\
Beijing 100190, P. R. China \\
Email: gzcui@math.ac.cn

\vskip 0.24cm

\noindent
Lei TAN \\
LAREMA, UMR 6093 CNRS \\
Universit\'e d'Angers, 2 bd Lavoisier \\
49045 Angers Cedex, France \\
Email: tanlei@math.uinv-angers.fr

\begin{thebibliography}{CPT}

\bibitem{A} L. Ahlfors, {\em Conformal Invariants: Topics in Geometric Function Theory}, McGraw-Hill Book Company, 1973.

\bibitem{B} A. Beardon, {\em Iteration of Rational Functions}, Graduate Texts in Math. 132, Springer-Verlag, 1991.

\bibitem{BCT} X. Buff, G. Cui and Tan L., Teichm\"{u}ller spaces and holomorphic dynamics, in {\em Handbook of Teichm\"{u}ller theory, Vol. IV}, ed. A. Papadopoulos, European Mathematical Society (2014), 717-756.

\bibitem{BD} B. Farb and D. Margalit, {\em A Primer on Mapping Class Groups} (PMS-49), Princeton University Press, 2011.

\bibitem{BP} A. Beardon and C. Pommerenke, The Poincar\'e metric of plane domains. {\em J. London Math. Soc.} 18 (1978), 475-483.

\bibitem{Bi} B. Bielefeld (editor), Conformal Dynamics Problem List. {\em Stony Brook IMS Preprint 1990/1}.

\bibitem{C} G. Cui, Conjugacies between rational maps and extremal quasiconformal maps, {\em Proc. Amer. Math. Soc.} 129, no. 7 (2001), 1949-1953.

\bibitem{CJ} G. Cui and Y. Jiang,  Geometrically finite and semi-rational branched coverings of the two-sphere, {\em Trans. Amer. Math. Soc.} 363 (2011), 2701-2714.

\bibitem{CPT} G. Cui, W. Peng and Tan Lei, On a theorem of Rees-Shishikura, {\em Annales de la Facult\'{e}s des Sciences de Toulouse}, S\'{e}r. 6 Vol. 21 no. 5 (2012), 981-993.

\bibitem{CT1} G. Cui and Tan Lei, A characterization of hyperbolic rational maps, {\em Invent. Math}., vol. 183 (2011), 451-516.

\bibitem{CT2} G. Cui and Tan Lei, Distortion control of conjugacies between quadratic polynomials, {\em Science China Mathematics}, vol. 53, no. 3 (2010), 625-634.

\bibitem{DD} R. and A. Douady, {\em  Alg\`ebre et th\'eories galoisiennes}, Cassini, 2005.

\bibitem{DH1} A. Douady and J. Hubbard, \'{E}tude dynamique des polyn\^{o}mes complexes I and II, Publication math\'ematiques d'Orsay, 84-02 and 85-04.

\bibitem{DH2} A. Douady and J. Hubbard, A proof of Thurston's topological characterization of  rational functions, {\em Acta Math.}, 171 (1993), 263-297.

\bibitem{GM} L. Goldberg and J. Milnor, Fixed points of polynomial maps, II. Fixed point portraits, {\em Ann. Sci. \'Ecole Norm. Sup.} 26 (1993), 51-98.

\bibitem{H1} P. Ha\"{\i}ssinsky, Chirurgie parabolique, {\em Comptes Rendus de l'Acad. des Sc. S\'er. I Math.,} { 327} (1998),  195-198.

\bibitem{H2} P. Ha\"{\i}ssinsky,  D\'eformation $J$-\'equivalente de polyn\^omes g\'eom\'etriquement finis, {\em Fund. Math.} {163} (2000), no. 2, 131-141.

\bibitem{H3} P. Ha\"{\i}ssinsky,  Pincements de polyn\^omes,  {\em Commentarii Helvetici Mathematici}, 77 (2002), no. 1, 1-23.

\bibitem{HT} P. Ha\"{\i}ssinsky and Tan L.,  Convergence of pinching deformations and matings of geometrically finite polynomials, {\em Fundamenta Mathematicae}, 181 (2004), 143-188.

\bibitem{JZ} Y. Jiang and G. Zhang, Combinatorial characterization of sub-hyperbolic rational maps, {\em Adv. Math.}, 221 (2009), 1990-2018.

\bibitem{Ki} J. Kiwi, Rational Rays and Critical Portraits of Complex Polynomials, Ph.D. thesis, Institute for Mathematical Sciences, Preprint ims97-15.

\bibitem{K1} T. Kawahira, Semiconjugacies between the Julia sets of geometrically finite rational maps, {\em Erg. Th. \& Dyn. Sys.} 23 (2003), 1125-1152.

\bibitem{K2} T. Kawahira, Semiconjugacies between the Julia sets of geometrically finite rational maps II, in {\em Dynamics on the Riemann sphere}, Bodilfest, ed. P. Hjorn and C. Petersen.

\bibitem{L} O. Lehto, {\em Univalent Functions and Teichm\"uller Spaces}, Springer-Verlag, 1987.

\bibitem{LM} M. Lyubich and Y. Minsky, Laminations in holomorphic dynamics, {\em J. Diff. Geom.} 47 (1997), 17-94.

\bibitem{PM} C. Petersen and D. Meyer, On the notions of mating, {\em Annales de la Facult\'{e}s des Sciences de Toulouse}, S\'{e}r. 6 Vol. 21 no. 5 (2012), 839-876.

\bibitem{PT} K. Pilgrim and Tan Lei, Combining rational maps and controlling obstructions, {\em Erg. Th. \& Dyn. Sys.}, 18 (1998), pp. 221-246.

\bibitem{R} J. Rivera-Letelier, A connecting lemma for rational maps satisfying a no growth condition, arXiv 2006.

\bibitem{Mak} P. Makienko, Unbounded components in parameter space of rational maps, {\em Conformal geometry and dynamics}, 4 (2000), 1-21.

\bibitem{Mas} B. Maskit, Parabolic elements in Kleinian groups, {\em Ann. of Math.} 117 (1983), 659-668.

\bibitem{Mc1} C. McMullen, {\em Complex Dynamics and Renormalization}, Annals of Mathematics Studies 135, Princeton University Press, 1994.

\bibitem{Mc2} C. McMullen, Self-similarity of Siegel disks and the Hausdorff dimension of Julia sets, {\em Acta Math.} 180 (1998), 247-292.

\bibitem{Mc3} C. McMullen, Cusps are dense, {\em Ann. of Math.} 133 (1991), 217-247.

\bibitem{Mc4} C. McMullen, Rational maps and Teichm\"{u}ller space, in {\em Linear and Complex Analysis Problem Book 3} (eds. V.P. Havin and N.K. Nikolskii), Lecture Notes in Mathematics Volume 1574 (1994), Springer-Verlag.

\bibitem{McS} C. McMullen and D. Sullivan, Quasiconformal homeomorphisms and dynamics III. The Teichm\"{u}ller space of a holomorphic dynamical system. {\em Adv. Math.} 135 (1998), no. 2, 351-395.

\bibitem{Mi} J. Milnor, {\em Dynamics in One Complex Variable,} Third Edition, Annals of Mathematics Studies 160, Princeton University Press, 2006.

\bibitem{O} R. Oudkerk, The parabolic implosion for $f_0(z)=z+z^{\nu+1}+\OOO(z^{\nu+2})$, Ph.D. thesis, University of Warwick, 1999.

\bibitem{P} Ch. Pommerenke, {\em Univalent Functions}, Vandenhoeck and Ruprecht in G\"{o}ttingen, 1975.

\bibitem{Re} M. Rees, Components of degree two hyperbolic rational maps, {\em Invent. math.}, Vol. 100 (1990), 357-382.

\bibitem{Ri} S. Rickman, Removability theorem for quasiconformal mappings, {\em Ann. Ac. Scient. Fenn.}, 449 (1969), 1-8.

\bibitem{RS} E. Reich and K. Strebel, Extremal quasiconformal mapping with given boundary values, {\em Bull. Amer. Math. Soc.}, 79 (1973), 488-490.

\bibitem{Sh} M. Shishikura, On a theorem of Mary Rees for matings of polynomials, in {\em The Mandelbrot set, Theme and variations}, ed. Tan Lei, LMS Lecture Note Series 274, Cambridge University Press, 2000, 289-305.

\bibitem{ST} M. Shishikura and Tan Lei,  A family of cubic rational maps and matings of cubic polynomials, {\em Experi. Math.}, 9 (2000), 29-53.

\bibitem{Str} K. Strebel, On quasiconformal mappings of open Riemann surfaces, {\em Comment. Math. Helv.}, 53 (1978), 301-321.

\bibitem{T1} Tan Lei, Matings of quadratic polynomials, {\em Ergodic Theory Dynam. Systems} 12 (1992),  no. 3, 589-620.

\bibitem{T2} Tan Lei, On pinching deformations of rational maps, {\em Ann. Scient. \'Ec. Norm. Sup.}, $4^e$ s\'erie, t.35 (2002), 353-370.

\bibitem{T3} Tan Lei, Stretching rays and their accumulations, following Pia Willumsen,  in {\em Dynamics on the Riemann sphere}, European Mathematical Society (2006) 183-208.

\bibitem{U} M. Urb\'{a}nski, Rational functions with no recurrent critical points, {\em Erg. Th. and Dyn. Sys.} 14 (1994), 1-29.

\bibitem{Wi} P. Willumsen, Holomorphic Dynamics: On accumulation of Stretching rays, Ph.D. thesis, Danmarks Tekniske Universitet, February 1997.

\end{thebibliography}
\end{document}